\documentclass[11pt]{article}

\usepackage{amssymb,amsthm,amsmath,hyperref}

\numberwithin{equation}{section}

\newtheorem{thm}{Theorem}[section]
\newtheorem{lem}{Lemma}[section]
\newtheorem{cor}{Corollary}[section]
\newtheorem{prop}{Proposition}[section]
\theoremstyle{definition}
\newtheorem{defn}{Definition}[section]
\theoremstyle{remark}
\newtheorem{rem}{Remark}[section]

\def\sumetage#1#2{
\sum_{\scriptstyle {#1}\atop\scriptstyle {#2}} }

\allowdisplaybreaks

\begin{document}
\title{Global wellposed problem for the 3-D incompressible anisotropic Navier-Stokes equations\thanks{
This work is supported by NSFC 10571158  and China Postdoctoral
Science Foundation 20060400335} }
\author{Ting Zhang\thanks{E-mail: zhangting79@hotmail.com}, Daoyuan Fang\thanks{E-mail:
dyf@zju.edu.cn}\\
\textit{\small Department of Mathematics, Zhejiang University,
Hangzhou 310027, China} }
\date{}
\maketitle
\begin{abstract}In this paper, we consider a global wellposed problem for the 3-D
incompressible anisotropic Navier-Stokes equations (\textit{ANS}).
In order to do so, we first introduce the scaling invariant
Besov-Sobolev type spaces, $B^{-1+\frac{2}{p},\frac{1}{2}}_{p}$
and $B^{-1+\frac{2}{p},\frac{1}{2}}_{p}(T)$, $p\geq2$. Then,  we
prove the global wellposedness for (\textit{ANS}) provided the
initial data are sufficient small
 compared to the horizontal viscosity in some suitable sense,
which is stronger than $B^{-1+\frac{2}{p},\frac{1}{2}}_{p}$ norm.
In particular, our results imply the global wellposedness of
(\textit{ANS}) with high oscillatory initial data.
\end{abstract}
\section{Introduction}
\subsection{Introduction to the anisotropic Navier-Stokes equations}
 In this paper, we are going to study the 3-D incompressible
anisotropic Navier-Stokes equations (\textit{ANS}), namely,
    \begin{equation}
    \left\{\begin{array}{l}
    u_t+u\cdot \nabla u-\nu_h\Delta_h u-\nu_3 \partial_{x_3}^2u=-\nabla
    P,\\
        \mathrm{div} u=0,\\
            u|_{t=0}=u_0,
    \end{array}
    \right.\label{ANS-E1.1}
    \end{equation}
where $u(t,x)$ and $P(t,x)$ denote the fluid velocity and the
pressure, respectively, the viscosity coefficients $\nu_h$ and
$\nu_3$ are two constants satisfying
    $$
    \nu_h>0,\ \nu_3\geq0,
    $$
$x=(x_h,x_3)\in\mathbb{R}^3$ and
$\Delta_h=\partial_{x_1}^2+\partial_{x_2}^2$. When
$\nu_h=\nu_3=\nu$, such system is the classical (isotropic)
Navier-Stokes system (\textit{NS}). It is
 appeared in geophysical fluids (see for instance
\cite{Chemin06}). In fact, instead of putting the classical
viscosity $-\nu\Delta$ in (\textit{NS}), meteorologists often
simulate the turbulent diffusion by putting a viscosity of the
form $-\nu_h\Delta_h-\nu_3\partial_{ x_3}^2$, where $\nu_h$ and
$\nu_3$ are empiric constants, and $\nu_3$ usually is  much
smaller than $\nu_h$. We refer to the book of J. Pedlosky
\cite{Pedlosky}, Chapter 4 for a more complete discussion. In
particular, in the studying of Ekman boundary layers
 for rotating fluids \cite{Chemin06,Ekman05,Grenier97}, it makes sense to consider anisotropic viscosities of
 the type
   $-\nu_h\Delta_h-\varepsilon\beta\partial_{x_3}^2$,
 where $\varepsilon$ is a very small parameter.
The system  ($ANS$) has been studied first by J.Y. Chemin, B.
Desjardins, I. Gallagher and E. Grenier in \cite{Chemin00} and D.
Iftimie in \cite{Iftimie02}, where the authors proved that such
system  is locally wellposed for initial data in the anisotropic
Sobolev space
    $$
    H^{0,\frac{1}{2}+\varepsilon}=\left\{
    u\in L^2(\mathbb{R}^3);
    \|u\|^2_{\dot{H}^{\frac{1}{2}+\varepsilon}}=
    \int_{\mathbb{R}^3}|\xi_3|^{1+2\varepsilon}|\hat{u}(\xi_h,\xi_3)|^2d\xi<+\infty
    \right\},
    $$
for some $\varepsilon>0$. Moreover, it has also been proved that
if the initial data $u_0$ is small enough in the sense of that
        \begin{equation}
        \|u_0\|^\varepsilon_{L^2}\|u_0\|^{1-\varepsilon}_{\dot{H}^{0,\frac{1}{2}+\varepsilon}}\leq
        c\nu_h\label{ANS-E1.2}
        \end{equation}
for some sufficiently small constant $c$, then the system
(\ref{ANS-E1.1}) is global wellposed.

Similar to  the classical Navier-Stokes equations, the system
(\textit{ANS}) has a scaling invariance. Indeed, if $u$ is a
solution of (\textit{ANS}) on a time interval $[0,T]$ with initial
data $u_0$, then the vector field $u_\lambda$ defined by
    $$
    u_\lambda(t,x)=\lambda u(\lambda^2t,\lambda x)
    $$
is also a solution of (\textit{ANS}) on the time interval
$[0,\lambda^{-2}T]$ with the initial data $\lambda u_0(\lambda
x)$. The smallness condition (\ref{ANS-E1.2}) is of course scaling
invariant, but the norm
$\|\cdot\|_{\dot{H}^{\frac{1}{2}+\varepsilon}}$ is not. M. Paicu
proved in \cite{Paicu05} a similar result for the system
(\textit{ANS}) with $\nu_3=0$ in the case of the initial data
$u_0\in B^{0,\frac{1}{2}}$. This space has a scaling invariant
norm. Then J.Y. Chemin and P. Zhang \cite{Chemin07} obtained  a
similar result in the scaling invariant space
$B^{-\frac{1}{2},\frac{1}{2}}_4$. Considering the periodic
anisotropic Naiver-Stokes equations, Paicu obtained the global
wellposedness in \cite{Paicu05-2}.

On the other hand, the classical (isotropic) Navier-Stokes system
(\textit{NS}) is globally wellposed for small initial data in
Besov norms of negative index. In \cite{Cannone93},     M.
Cannone, Y. Meyer and F. Planchon proved that: if the initial data
satisfy
    $$
    \|u_0\|_{\dot{B}^{-1+\frac{3}{p}}_{p,\infty}}\leq c\nu,
    $$
for  $p>3$ and some constant $c$ small enough, then the classical
Navier-Stokes system (\textit{NS}) is globally wellposed. Then, H.
Koch and D. Tataru generalized this theorem to the $BMO^{-1}$ norm
(see \cite{Koch01}), D. Iftimie in \cite{Iftimie99} obtained the
global wellposedness in anisotropic spaces $H^{s_1,s_2,s_3}$ and
$B^{0,\frac{1}{2}}$. Recently, J.Y. Chemin and I. Gallagher
\cite{Chemin07-0} proved that if a certain nonlinear function of
the initial data is small enough, then there is a global solution
to the Navier-Stokes equations (\textit{NS}).

Let $\phi_0(x_3)$ be a function in the Schwartz space
$\mathcal{S}(\mathbb{R})$ satisfying $\mathrm{supp}
\hat{\phi}_0\subset \mathcal{C}_v$, $\phi_1(x_h)$ be a function in
the Schwartz space $\mathcal{S}(\mathbb{R}^2)$ satisfying
$\mathrm{supp} \hat{\phi}_1\subset \mathcal{C}_h$, where
$\mathcal{C}_h$ (resp. $\mathcal{C}_v$) is a ring of
$\mathbb{R}^2_h$ (resp. $\mathbb{R}_v$). The mentioned results
imply that the system ($NS$) is globally wellposed for the initial
data $u_0^\varepsilon$ defined by
    \begin{equation}
      u_0^\varepsilon(x)=\varepsilon^{-\frac{1}{2}}\sin(\frac{x_1}{\varepsilon})
      (0,-\partial_{x_3}(\phi_0\phi_1),\partial_{x_2}(\phi_0\phi_1))\label{ANS-E1.3}
    \end{equation}
with small enough $\varepsilon$. The goal of our work is to prove
a result of this type for the anisotropic Navier-Stokes system
(\ref{ANS-E1.1}).

\subsection{Statement of the results.}
As in \cite{Chemin07}, let us begin with the definition of the
spaces, which we will be going to work with. It requires an
anisotropic version of dyadic decomposition of the Fourier space,
let us first recall the following operators of localization in
Fourier space, for $(k,l)\in\mathbb{Z}^2$,
    $$
        \Delta^h_k a=\mathcal{F}^{-1}(\varphi(2^{-k}|\xi_h|)\hat{a}),
         \ \Delta^v_l a=\mathcal{F}^{-1}(\varphi(2^{-l}|\xi_3|)\hat{a}),
    $$
        $$
            S^h_ka=\sum_{k'\leq k-1}\Delta^h_{k'}a,
            \ S^v_la=\sum_{l'\leq l-1}\Delta^v_{l'}a,
        $$
    where $\mathcal{F}a$ or $\hat{a}$ denotes the Fourier transform
    of the function $a$, and $\varphi$ is a function in
    $\mathcal{D}((\frac{3}{4},\frac{8}{3}))$ satisfying
        $$
            \sum_{j\in\mathbb{Z}}\varphi(2^{-j}\tau)=1,
            \ \forall\ \tau>0.
        $$

 Our main motivation to introduce the following
spaces is to find a scaling invariant Besov-Sobolve type space
such that  $u_0^\varepsilon$ can be small. According to the
definitions of $B^{0,\frac{1}{2}}$  (in \cite{Iftimie99,Paicu05})
and $B^{-\frac{1}{2},\frac{1}{2}}_4$ (in \cite{Chemin07}), we give
the definition of $B^{-1+\frac{2}{p},\frac{1}{2}}_p$,
$p\in[2,\infty)$, as follows.

\begin{defn}
  We denote by $B^{-1+\frac{2}{p},\frac{1}{2}}_p$, $p\in[2,\infty)$, the space of distributions, which
  is the completion of $\mathcal{S}(\mathbb{R}^3)$ by the following
  norm:
    $$
        \|a\|_{B^{-1+\frac{2}{p},\frac{1}{2}}_p}=\sum_{l\in\mathbb{Z}}2^{\frac{l}{2}}
        \left(\sum_{k=l-1}^\infty 2^{(-2+\frac{4}{p})k}\|\Delta^h_k\Delta^v_l
        a\|_{L^p_h(L^2_v)}^2\right)^\frac{1}{2}+\sum_{j\in\mathbb{Z}}2^{\frac{j}{2}}\|S^h_{j-1}\Delta^v_ja\|_{
        L^2(\mathbb{R}^3)}.
    $$
    Let $B^{0,\frac{1}{2}}:=B^{0,\frac{1}{2}}_2$ and $
    \|a\|_{B^{0,\frac{1}{2}}}\simeq \sum_{j\in\mathbb{Z}}2^{\frac{j}{2}}\|\Delta^v_ja\|_{
        L^2(\mathbb{R}^3)}$.
\end{defn}

To study the evolution of (\ref{ANS-E1.1}) with initial data in
$B^{-1+\frac{2}{p},\frac{1}{2}}_p$, we need also to introduce the
following space.

\begin{defn}
  We denote by $B^{-1+\frac{2}{p},\frac{1}{2}}_p(T)$ the space of distributions, which
  is the completion of $C^\infty([0,T];\mathcal{S}(\mathbb{R}^3))$ by the following
  norm:
            \begin{eqnarray*}
        &&\|a\|_{B^{-1+\frac{2}{p},\frac{1}{2}}_p(T)}\\
        &=&\sum_{l\in\mathbb{Z}}2^{\frac{l}{2}}\left[
        \sum_{k=l-1}^\infty\left( 2^{(-2+\frac{4}{p})k}\|\Delta^h_k\Delta^v_l
        a\|_{L^\infty_T(L^p_h(L^2_v))}^2
            +\nu_h
             2^{\frac{4}{p}k}\|\Delta^h_k\Delta^v_l
        a\|_{L^2_T(L^p_h(L^2_v))}^2\right.\right.\\
                &&\left.\left. +\nu_3
             2^{(-2+\frac{4}{p})k+2l}\|\Delta^h_k\Delta^v_l
        a\|_{L^2_T(L^p_h(L^2_v))}^2\right)\right]^\frac{1}{2}
        +\sum_{j\in\mathbb{Z}}2^{\frac{j}{2}}\left(
        \|S^h_{j-1}\Delta^v_ja\|_{
        L^\infty_T(L^2(\mathbb{R}^3))}\right. \\
                &&\left.+\nu_h^\frac{1}{2}\|\nabla_hS^h_{j-1}\Delta^v_ja\|_{
        L^2_T(L^2(\mathbb{R}^3))}+\nu_3^\frac{1}{2}\|\partial_3S^h_{j-1}\Delta^v_ja\|_{
        L^2_T(L^2(\mathbb{R}^3))}
        \right).
  \end{eqnarray*}
    Let $B^{0,\frac{1}{2}}(T):=B^{0,\frac{1}{2}}_2(T)$ and
        \begin{equation*}
        \|a\|_{B^{0,\frac{1}{2}}(T)}\simeq\sum_{j\in\mathbb{Z}}2^{\frac{j}{2}}\left(
        \|\Delta^v_ja\|_{
        L^\infty_T(L^2(\mathbb{R}^3))}+\nu_h^\frac{1}{2}\|\nabla_h\Delta^v_ja\|_{
        L^2_T(L^2(\mathbb{R}^3))} +\nu_3^\frac{1}{2}\|\partial_3\Delta^v_ja\|_{
        L^2_T(L^2(\mathbb{R}^3))}
        \right).
        \end{equation*}
\end{defn}

In our global result, we need that the initial data $u_0$ and a
certain nonlinear function of the initial data $u_F\cdot\nabla
u_F$ are small enough in some suitable sense, where
    \begin{equation}
            u_F:=e^{\nu_ht\Delta_h+\nu_3t\partial_3^2}u_{0hh}
            \ \textrm{ and }\             u_{0hh}=\sum_{k\geq l-1}\Delta^h_k\Delta^v_l
            u_0.           \label{ANS-E1.9}
    \end{equation}

\begin{defn}\label{ANS-D2}
  We define  $[\cdot]_{E^p_T}$ by:
    $$
        [a]_{E^p_T}:=
            \|a\|_{B^{-1+\frac{2}{p},\frac{1}{2}}_p}+
        \|a_{F}\cdot\nabla
         a_{F}\|_{{L}^1_T(B^{0,\frac{1}{2}})},
    $$
where
    $$
        \|f\|_{{L}^1_T(B^{0,\frac{1}{2}})}
        :=\int^T_0\sum_{j\in\mathbb{Z}}2^{\frac{j}{2}}\|\Delta^v_jf\|_{L^2(\mathbb{R}^3)}dt,
    $$
    $$
    a_F:=
            e^{\nu_ht\Delta_h+\nu_3t\partial_3^2}a_{hh},\
            a_{hh}:=\sum_{k\geq l-1}\Delta^h_k\Delta^v_l a.
    $$
\end{defn}

Now, we present the main results of this paper, which cover the
results in \cite{Chemin07,Paicu05} and partial result in
\cite{Iftimie99}.

\begin{thm}\label{ANS-T1}
  A constant $c$ exists such that, if the divergence free vector field $u_0\in
  B^{-1+\frac{2}{p},\frac{1}{2}}_p$, $p\geq2$   and $[u_0]_{E^p_\infty}\leq c\nu_h$,
  then,   the system (\ref{ANS-E1.1}) with initial data $u_0$ has a
  unique global solution $u\in
  B^{-1+\frac{2}{p},\frac{1}{2}}_p(\infty)\cap C([0,\infty);B^{-1+\frac{2}{p},\frac{1}{2}}_p)$, and $\|u\|_{
  B^{-1+\frac{2}{p},\frac{1}{2}}_p(\infty)}$ is independent of $\nu_3$.

  Furthermore, let $u_i\in B^{-1+\frac{2}{p},\frac{1}{2}}_p(T)$ be the solution for
  the system (\ref{ANS-E1.1}) with initial data $u_{0i}$, $p\geq2$,
  $i=1,2$. If  $\nu_3>0$ and $u_{01}-u_{02}\in L^2$, then we have
    \begin{eqnarray}
        &&\|u_1-u_2\|_{L^\infty_T(L^2(\mathbb{R}^3))}\label{ANS-E1.5}\\
    &\leq& \|u_{01}-
    u_{02}\|_{L^2(\mathbb{R}^3)}\exp\left
        \{C\nu_h^{-1}(\nu_h^{-\frac{p+1}{p-1}}+\nu_3^{-\frac{p+1}{p-1}})
        \left(
        \sum_{i=1}^2\|u_i\|_{B^{-1+\frac{2}{p},\frac{1}{2}}_p(T)
        }        \right)^{\frac{2p}{p-1}}\right\}.\nonumber
    \end{eqnarray}
\end{thm}

In what follows, we always use $C$ to denote a generic positive
constant independent of $\nu_3$. Repeating the proof of Theorem
\ref{ANS-T1}, we may conclude the following theorem concerning
local wellposedness for large data.

\begin{thm}\label{ANS-T2}
If the divergence free vector field $u_0\in
  B^{-1+\frac{2}{p},\frac{1}{2}}_p$, $p\geq 2$ and $[u_0]_{E^p_T}<\infty$,
  then there exists a constant $T_0\in(0,T]$  such that  the system (\ref{ANS-E1.1}) has a
  unique  solution $u\in
  B^{-1+\frac{2}{p},\frac{1}{2}}_p(T_0)\cap C([0,T_0];B^{-1+\frac{2}{p},\frac{1}{2}}_p)$, and $\|u\|_{
  B^{-1+\frac{2}{p},\frac{1}{2}}_p(T_0)}$ is independent of $\nu_3$.
\end{thm}

\begin{rem}
  These theorems imply that the third
viscosity coefficient $\nu_3$ do not play a role except the
continuous dependence (\ref{ANS-E1.5}).
\end{rem}

\begin{prop}\label{ANS-P2}
If $p\in[2,4]$, we have
    $$
    \sum_{j\in\mathbb{Z}}2^{\frac{j}{2}}\int^\infty_0\|\Delta^v_j(a_F\cdot\nabla a_F)\|_{L^2(\mathbb{R}^3)}dt
         \lesssim \nu_h^{-1}\|a\|_{B^{-1+\frac{2}{p},\frac{1}{2}}_p}^2
         $$
    and
    $$[a]_{E^p_\infty}\lesssim\|a\|_{B^{-1+\frac{2}{p},\frac{1}{2}}_p}+\nu_h^{-1}\|a\|_{B^{-1+\frac{2}{p},\frac{1}{2}}_p}^2.$$
\end{prop}

\begin{rem}
This proposition will be proved in Section \ref{ANS-Sec5}. It
implies that if $p\in[2,4]$, then the condition
$[u_0]_{E^p_\infty}\leq
  c\nu_h$ in Theorem \ref{ANS-T1} can be replaced by
    $$
    \|u_0\|_{B^{-1+\frac{2}{p},\frac{1}{2}}_p}\leq
  c\nu_h,
    $$
    and the condition $[u_0]_{E^p_T}<\infty
  $ in Theorem \ref{ANS-T2} can be omitted.
\end{rem}

\begin{rem}
  Using the similar argument to that in the proof of Proposition \ref{ANS-P2}, we
     obtain
     $$
    2^{\frac{j}{2}}
    \|\Delta^v_j(a_F\cdot\nabla
    a_F)\|_{L^q_TL^\frac{p}{2}(\mathbb{R}^3)}
         \lesssim d_j\nu_h^{-\frac{1}{q}}2^{(3-\frac{4}{p}-\frac{2}{q})j}\|a\|_{B^{-1+\frac{2}{p},\frac{1}{2}}_p}^2,
         \ q<\frac{p}{p-2},
         $$
where $(d_k)_{k\in\mathbb{Z}}$ denotes a generic element of the
sphere of  $l^1(\mathbb{Z})$. From $3-\frac{4}{p}-\frac{2}{q}=0$
and $q\geq1$, we have $p\leq4$.  Thus, we think $p=4$ seems a
special point.
\end{rem}

The following proposition, which will be proved in Section
\ref{ANS-Sec6-0}, shows that Theorem \ref{ANS-T1} can be applied
to initial data given by (\ref{ANS-E1.3}).

\begin{prop}\label{ANS-P1}
 Let  $\phi(x_h)$ and $\psi(x_h)$  be in the Schwartz space
$\mathcal{S}(\mathbb{R}^2)$, $\mathrm{supp} \hat{\phi}$,
$\mathrm{supp} \hat{\psi}\subset \mathcal{C}_h$, where
$\mathcal{C}_h$ is a ring of $\mathbb{R}^2$. Denote
$\phi_\varepsilon(x_h)=e^{ix_1/\varepsilon}\phi(x_h)$ and
$\psi_\varepsilon(x_h)=e^{ix_1/\varepsilon}\psi(x_h)$, we have,
    for any positive $\varepsilon$, $q\in(1,\infty]$,
$\sigma>0$,
    $\alpha\in(0,2(1-\frac{1}{q}))$,
                $$
                              \|\phi_\varepsilon\|_{\tilde{B}^{-\sigma}_{q,1}}\leq C_\phi
               \varepsilon^{\sigma},
               $$
               \begin{equation}
                              \|\phi_\varepsilon\|_{\dot{B}^{-\alpha}_{q,1}}\leq C_\phi
               \varepsilon^{\alpha},\label{ANS-E1.4}
               \end{equation}
                    $$
                    \|\phi_\varepsilon\|_{\dot{B}^{-\sigma}_{q,\infty}}\geq C_{\phi}
                    \varepsilon^{\sigma},
                    $$
               $$
        \int^\infty_0\|\phi_{\varepsilon,F'}\psi_{\varepsilon,F'}\|_{L^2(\mathbb{R}^2)}dt\leq C_{\phi,\psi,\nu_h}
                    \varepsilon^2,
               $$
    where
        $$
            \phi_{\varepsilon,F'}(t)=e^{\nu_ht\Delta_h}\sum_{k\geq 0}\Delta^h_k
            \phi_{\varepsilon},
            \        \psi_{\varepsilon,F'}(t)=e^{\nu_ht\Delta_h}\sum_{k\geq 0}\Delta^h_k
            \psi_{\varepsilon},
        $$
    $$\|a\|_{\tilde{B}^{-\sigma}_{q,1}}=\|S_0^ha\|_{L^2}+
    \sum_{k=0}^\infty 2^{-\sigma k}\|\Delta_k^ha\|_{L^q},\ \|a\|_{\dot{B}^{-\alpha}_{q,1}}=
    \sum_{k\in\mathbb{Z}} 2^{-\alpha k}\|\Delta_k^ha\|_{L^q}
    $$
    and
    $\|a\|_{\dot{B}^{-\sigma}_{q,\infty}}=\sup_{k\in \mathbb{Z}}
     2^{-\sigma k}\|\Delta_k^ha\|_{L^q}$.
\end{prop}

\begin{rem}
  From Proposition \ref{ANS-P1}, we get
        $$\|u_0^\varepsilon\|_{B^{-\frac{1}{2},\frac{1}{2}}_4}\simeq C_{\phi_0,\phi_1},
        \ \textrm{ and }\ [u_0^\varepsilon]_{E^p_\infty}\leq C_{\phi_0,\phi_1,\nu_h} \varepsilon^{\frac{1}{2}-\frac{2}{p}},
                \ \textrm{ for all }\ p>4.$$
Theorem \ref{ANS-T1} and Proposition \ref{ANS-P1} imply that the
anisotropic Navier-Stokes system (\ref{ANS-E1.1}) with initial
data $u_0^{\varepsilon,q}$, which defined by
   $$
      u_0^{\varepsilon,q}(x)=\varepsilon^{-1+\frac{2}{q}}\sin(\frac{x_1}{\varepsilon})
      (0,-\partial_{x_3}(\phi_0\phi_1),\partial_{x_2}(\phi_0\phi_1)),\ \forall
      \ q\geq 2,
   $$
    is globally wellposed  when $\varepsilon$ is small enough.
\end{rem}

At last, we give an imbedding result in the following proposition,
which will be proved in Section \ref{ANS-Sec6}.
\begin{prop}\label{ANS-P3}
For $p\geq2$, we have
  $$B^{0,\frac{1}{2}}\subset B^{-1+\frac{2}{p},\frac{1}{2}}_p\subset \dot{B}^{-1}_{\infty,2}\subset BMO^{-1}\subset
  \dot{B}^{-1}_{\infty,\infty}=C^{-1},$$
where
$\|f\|_{\dot{B}^{-1}_{\infty,q}}=\left\|2^{-k}\|\Delta_kf\|_{L^\infty}\right\|_{l^q_k}$,
$\Delta_k a=\mathcal{F}^{-1}(\varphi(2^{-k}|\xi|)\hat{a}),$
        $$\|f\|_{BMO^{-1}}:=\|f\|_{\dot{B}^{-1}_{\infty,\infty}}
            +\sup_{x\in\mathbb{R}^3,R>0}R^{-\frac{3}{2}}\left(\int_{P(x,R)}|e^{t\Delta}f(y)|^2dydt
            \right)^\frac{1}{2},$$
 $P(x,R)=[0,R^2]\times B(x,R)$ and
$B(x,R):=\{y\in\mathbb{R}^3;|x-y|\leq R\}$ (see \cite{Koch01}).
\end{prop}

\subsection{Structure of the proof of Theorem \ref{ANS-T1}.}

The purpose of Section \ref{ANS-Sec2} is to establish some results
about anisotropic Littlewood-Paley theory, which will be of
constant use in what follows.

Section \ref{ANS-Sec3} will be devoted to the proof of the
existence of a solution of (\ref{ANS-E1.1}). In order to do it, we
shall search for a solution of the form, (following the idea in
\cite{Chemin07})
    $$
    u=u_F+w, \ \textrm{ and }\ w\in B^{0,\frac{1}{2}}(\infty).
    $$

    In Section \ref{ANS-Sec4}, we shall prove the uniqueness in the
    following way. First, we shall establish a regularity theorem,
   which claim that if $u\in B^{-1+\frac{2}{p},\frac{1}{2}}_{p}(T)$ is a
    solution of (\ref{ANS-E1.1}) with
    $[u_0]_{E^p_T}<\infty$, then there exists $T_1\in(0,T]$ such that $w=u-u_F\in B^{0,\frac{1}{2}}(T_1)$.
    Therefore, looking at the equation of $w$, we shall prove the
    uniqueness of the solution $u$ in the space
    $u_F+B^{0,\frac{1}{2}}(T_1)$. Since $u\in  C([0,T];B^{-1+\frac{2}{p},\frac{1}{2}}_p)$, one can
    easily obtain the uniqueness of the solution $u$ on $[0,T]$.

In Section \ref{ANS-Sec5-0}, we shall prove that if  $\nu_3>0$,
then the continuous dependence of the solution on the initial data
holds.

     We should mention that the methods
    introduced by  Chemin-Gallagher in \cite{Chemin07-0}, Chemin-Zhang in
    \cite{Chemin07}, Koch-Tataru in \cite{Koch01} and  Paicu in \cite{Paicu05}
   will play a crucial role
in our proof here.

\section{Anisotropic Littlewood-Paley theory}\label{ANS-Sec2}
At first, we list anisotropic Berstein inequalities in the
following, (please see the detail in \cite{Chemin07,Paicu05}).

\begin{lem}\label{ANS-L2.1}
Let $\mathcal{B}_h$ (resp. $\mathcal{B}_v$) be a ball  of
$\mathbb{R}^2_h$ (resp. $\mathbb{R}_v$), and $\mathcal{C}_h$
(resp. $\mathcal{C}_v$) be a ring of $\mathbb{R}^2_h$ (resp.
$\mathbb{R}_v$).  Then, for $1\leq p_2\leq p_1\leq \infty$ and
$1\leq q_2\leq q_1\leq \infty$, there holds:
\begin{enumerate}
  \item If the support of $\hat{a}$ is included in
  $2^k\mathcal{B}_h$, then
        $$
        \|\partial_{h}^\beta a\|_{L^{p_1}_h(L^{q_1}_v)}\lesssim
        2^{k(|\beta|+2(\frac{1}{p_2}-\frac{1}{p_1}))}
                \| a\|_{L^{p_2}_h(L^{q_1}_v)},
        $$
  where $\partial_{h}:=\partial_{x_h}$.
  \item If the support of $\hat{a}$ is included in
  $2^l\mathcal{B}_v$, then
        $$
        \|\partial_{3}^N a\|_{L^{p_1}_h(L^{q_1}_v)}\lesssim
        2^{l(N+\frac{1}{q_2}-\frac{1}{q_1})}
                \| a\|_{L^{p_1}_h(L^{q_2}_v)},
        $$
         where $\partial_{3}:=\partial_{x_3}$.
  \item If the support of $\hat{a}$ is included in
  $2^k\mathcal{C}_h$, then
        $$
        \| a\|_{L^{p_1}_h(L^{q_1}_v)}\lesssim
        2^{-kN}
               \sup_{|\beta|=N} \|\partial_h^\beta a\|_{L^{p_1}_h(L^{q_1}_v)}.
        $$
  \item If the support of $\hat{a}$ is included in
  $2^l\mathcal{C}_v$, then
        $$
        \| a\|_{L^{p_1}_h(L^{q_1}_v)}\lesssim
        2^{-lN}
               \|\partial_3^N a\|_{L^{p_1}_h(L^{q_1}_v)}.
        $$
\end{enumerate}
  \end{lem}

Let us state two corollaries of this lemma, the proofs of which
are obvious and thus omitted.

\begin{cor}\label{ANS-C1.1}
  The space $B^{0,\frac{1}{2}}$ is continuously embedded in the
  space $B^{-1+\frac{2}{p},\frac{1}{2}}_p$ and so is
  $B^{0,\frac{1}{2}}(T)$ in $B^{-1+\frac{2}{p},\frac{1}{2}}_p(T)$, $p\geq2$.
  Moreover, the space $B^{0,\frac{1}{2}}(T)$ is continuously
  embedded in the space $L^\infty_T(L^2_h(L^\infty_v))$.
  Furthermore,  The space $B^{-1+\frac{2}{p},\frac{1}{2}}_p$ is continuously embedded in the
  space $B^{-1+\frac{1}{p},\frac{1}{2}}_{2p}$ and so is
  $B^{-1+\frac{2}{p},\frac{1}{2}}_p(T)$ in $B^{-1+\frac{1}{p},\frac{1}{2}}_{2p}(T)$, $p\geq2$.
\end{cor}

\begin{cor}\label{ANS-C1.2}
  If a belongs to $B^{-1+\frac{2}{p},\frac{1}{2}}_p(T)$, $p\geq2$, then we
  have
    $$
    \sum_{l\in\mathbb{Z}}2^{\frac{l}{2}}\left(
    \sum_{k\in\mathbb{Z}}2^{2k(-1+\frac{2}{p})}\|\Delta^h_k\Delta^v_l
    a(0)\|^2_{L^p_h(L^2_v)}
    \right)^\frac{1}{2}\lesssim
    \|a(0)\|_{B^{-1+\frac{2}{p},\frac{1}{2}}_p}
    $$
  and
        $$
        \sum_{l\in\mathbb{Z}}2^{\frac{l}{2}}\left[
    \sum_{k\in\mathbb{Z}}\left(2^{k(-2+\frac{4}{p})}\|\Delta^h_k\Delta^v_l
    a\|^2_{L^\infty_T(L^p_h(L^2_v))}+
    \nu_h2^{\frac{4k}{p}}\|\Delta^h_k\Delta^v_l
    a\|^2_{L^2_T(L^p_h(L^2_v))}\right)
    \right]^\frac{1}{2}\lesssim
    \|a\|_{B^{-1+\frac{2}{p},\frac{1}{2}}_p(T)}.
        $$
\end{cor}

\noindent\textit{Notations.} In what follows, as in
\cite{Chemin07}, we make the convention that
$(c_k)_{k\in\mathbb{Z}}$ (resp. $(d_k)_{k\in\mathbb{Z}}$) denotes
a generic element of the sphere of $l^2(\mathbb{Z})$ (resp.
$l^1(\mathbb{Z})$). Moreover, $(c_{k,l})_{(k,l)\in\mathbb{Z}^2}$
denotes a generic element of the sphere of $l^2(\mathbb{Z}^2)$ and
$(d_{k,l})_{(k,l)\in\mathbb{Z}^2}$ denotes a generic sequence such
that
            $$
            \sum_{l\in\mathbb{Z}}\left(
            \sum_{k\in\mathbb{Z}}d^2_{k,l}
            \right)^\frac{1}{2}=1.
            $$

The following lemma will be of frequent use in this work. It
describes some estimates of dyadic parts of functions in
$B^{-1+\frac{2}{p},\frac{1}{2}}_p(T)$.

\begin{lem}\label{ANS-L2.2}
  For any $a\in B^{-1+\frac{2}{p'},\frac{1}{2}}_{p'}(T)$, $p'\geq2$, we have
    $$
    S_{k,l}(a)
                \lesssim \frac{p'}{\beta}d_{k,l}2^{\frac{\beta k}{p'}}2^{-\frac{l}{2}}
    \|a\|_{B^{-1+\frac{2}{p'},\frac{1}{2}}_{p'}(T)}
    $$
    and
       $$
    S_{k}(a)
                \lesssim\frac{p'}{\beta} c_{k}2^{\frac{\beta k}{p'}}
    \|a\|_{B^{-1+\frac{2}{p'},\frac{1}{2}}_{p'}(T)},
       $$
    where $\beta>0$, and
    $$S_{k,l}(a)=\sum_{k'\leq k-1} 2^{(-1+\frac{2+\beta}{p'})k'}
            \left(\|\Delta^h_{k'}\Delta^v_la\|_{L^\infty_T(L^{p'}_h(L^2_v))}+\nu_h^\frac{1}{2}\|\nabla_h
            \Delta^h_{k'}\Delta^v_la\|_{
    L^2_T(L^{p'}_h(L^2_v))}
            \right),$$
    $$S_k(a):= \sum_{k'\leq k-1}2^{(-1+\frac{2+\beta}{p'})k'}
            \left(\|\Delta^h_{k'}a\|_{L^\infty_T(L^{p'}_h(L^\infty_v))}+\nu_h^\frac{1}{2}\|\nabla_h\Delta^h_{k'} a\|_{
    L^2_T(L^{p'}_h(L^\infty_v))} \right).$$
\end{lem}
\noindent\textbf{Proof.}
  Since
    \begin{eqnarray*}
            2^{\frac{l}{2}}2^{-\frac{\beta k}{p'}}S_{k,l}&\leq&2^\frac{l}{2}\sum_{k'\leq k-1}
            2^{-\frac{\beta}{p'}(k-k')}2^{(-1+\frac{2}{p'})k'}
            \left(\|\Delta^h_{k'}\Delta^v_la\|_{L^\infty_T(L^{p'}_h(L^2_v))}\right.\\
                    &&\left.+\nu_h^\frac{1}{2}\|\nabla_h
            \Delta^h_{k'}\Delta^v_la\|_{
    L^2_T(L^{p'}_h(L^2_v))}
            \right),
    \end{eqnarray*}
using Young's inequality, we obtain
             \begin{eqnarray*}
            2^{\frac{l}{2}}\left(\sum_{k\in\mathbb{Z}}
            2^{-\frac{2\beta k}{p'}}S_{k,l}^2\right)^\frac{1}{2}&\lesssim&2^\frac{l}{2}\frac{p'}{\beta}
            \left(\sum_{k'\in\mathbb{Z}}
            2^{(-2+\frac{4}{p'})k'}
            \left(\|\Delta^h_{k'}\Delta^v_la\|_{L^\infty_T(L^{p'}_h(L^2_v))}\right.\right.\\
                    &&\left.\left.+\nu_h^\frac{1}{2}\|\nabla_h
            \Delta^h_{k'}\Delta^v_la\|_{
    L^2_T(L^{p'}_h(L^2_v))}
            \right)^2\right)^\frac{1}{2}.
            \end{eqnarray*}
Combining it with Corollary \ref{ANS-C1.2}, we can easily obtain
the first inequality.

To get the second inequality, we shall prove that, for any
$(c_k)_{k\in\mathbb{Z}}$, we have
    \begin{equation}
    I(a):=\sum_{k\in\mathbb{Z}}2^{-\frac{\beta k}{p'}}S_k(a)c_k\lesssim
    \frac{p'}{\beta}\|a\|_{B^{-1+\frac{2}{p'},\frac{1}{2}}_{p'}(T)}.
    \label{ANS-E2.2}
    \end{equation}
 Using Lemma  \ref{ANS-L2.1}, we get
        $$
    S_k(a)\lesssim \sum_{k'\leq k-1}\sum_{l\in\mathbb{Z}}
    2^{\frac{l}{2}+(-1+\frac{2+\beta}{p'})k'}\left(
        \|\Delta^h_{k'}\Delta^v_la\|_{L^\infty_T(L^{p'}_h(L^2_v))}
        +\nu_h^\frac{1}{2}\|\nabla_h
            \Delta^h_{k'}\Delta^v_la\|_{
    L^2_T(L^{p'}_h(L^2_v))}
    \right),
        $$
and
    \begin{eqnarray*}
    I(a)&\lesssim& \sum_{l\in\mathbb{Z}}
    2^{\frac{l}{2}}\sumetage{(k,k')\in\mathbb{Z}^2}{k'\leq k-1}2^{-\frac{\beta}{p'}(k-k')}
    2^{(-1+\frac{2}{p'})k'}c_k\left(
        \|\Delta^h_{k'}\Delta^v_la\|_{L^\infty_T(L^{p'}_h(L^2_v))}\right.\\
                &&\left.  +\nu_h^\frac{1}{2}\|\nabla_h
            \Delta^h_{k'}\Delta^v_la\|_{
    L^2_T(L^{p'}_h(L^2_v))}
    \right).
    \end{eqnarray*}
Using the Cauchy-Schwarz inequality and Young's inequality, we
have
        \begin{eqnarray*}
    &&I(a)\\
    &\lesssim& \big(\sumetage{(k,k')\in\mathbb{Z}^2}{k'\leq k-1}2^{-\frac{\beta}{p'}(k-k')}
    c_k^2\big)^\frac{1}{2}\sum_{l\in\mathbb{Z}}
    2^{\frac{l}{2}}\big(\sumetage{(k,k')\in\mathbb{Z}^2}{k'\leq k-1}2^{-\frac{\beta}{p'}(k-k')}
    2^{(-2+\frac{4}{p'})k'}\\
            &&\left.\times\left(
        \|\Delta^h_{k'}\Delta^v_la\|_{L^\infty_T(L^{p'}_h(L^2_v))} +\nu_h^\frac{1}{2}\|\nabla_h
            \Delta^h_{k'}\Delta^v_la\|_{
    L^2_T(L^{p'}_h(L^2_v))}
    \right)^2\right)^\frac{1}{2}\\
                    &\lesssim& \sqrt{\frac{p'}{\beta}}\sum_{l\in\mathbb{Z}}
    2^{\frac{l}{2}}\big(\sumetage{(k,k')\in\mathbb{Z}^2}{k'\leq k-1}2^{-\frac{\beta}{p'}(k-k')}
    2^{(-2+\frac{4}{p'})k'}\left(
        \|\Delta^h_{k'}\Delta^v_la\|_{L^\infty_T(L^{p'}_h(L^2_v))}\right.\\
            &&\left. +\nu_h^\frac{1}{2}\|\nabla_h
            \Delta^h_{k'}\Delta^v_la\|_{
    L^2_T(L^{p'}_h(L^2_v))}
    \right)^2\big)^\frac{1}{2}\\
            &\lesssim& \frac{p'}{\beta}\sum_{l\in\mathbb{Z}}
    2^{\frac{l}{2}}\left(\sum_{k'\in\mathbb{Z}}
    2^{(-2+\frac{4}{p'})k'}\left(
        \|\Delta^h_{k'}\Delta^v_la\|_{L^\infty_T(L^{p'}_h(L^2_v))}
         +\nu_h^\frac{1}{2}\|\nabla_h
            \Delta^h_{k'}\Delta^v_la\|_{
    L^2_T(L^{p'}_h(L^2_v))}
    \right)^2\right)^\frac{1}{2}\\
            &\lesssim& \frac{p'}{\beta}\|a\|_{B^{-1+\frac{2}{p'},\frac{1}{2}}_{p'}(T)},
    \end{eqnarray*}
which proves (\ref{ANS-E2.2}) and thus Lemma \ref{ANS-L2.2}.
{\hfill $\square$\medskip}

With Lemma \ref{ANS-L2.2}, we are going to state two lemmas, which
is very closed to Sobolev's embedding Theorem and will be of
constant use in the proof of Theorem \ref{ANS-T1}.

\begin{lem}\label{ANS-L2.3}
  The space $B^{-1+\frac{2}{p'},\frac{1}{2}}_{p'}(T)$ is included in
  $L^\frac{2p'}{p'-2}_T(L^{p'}_h(L^\infty_v))$, $p'>2$. Moreover precisely, if $a\in
  B^{-1+\frac{2}{p'},\frac{1}{2}}_{p'}(T)$,
   then, we have
    \begin{equation}
   \|\Delta^v_ja\|_{L^\frac{2p'}{p'-2}_T(L^{p'}_h(L^2_v))}
   \lesssim \sqrt{\frac{p'}{\beta}} d_j\nu_h^{\frac{1}{p'}-\frac{1}{2}}2^{-\frac{j}{2}}
   \|a\|_{B^{-1+\frac{2}{p'},\frac{1}{2}}_{p'}(T)}\label{ANS-E2.2--1}
   \end{equation}
  and
    \begin{equation}
    \|a\|_{L^\frac{2p'}{p'-2}_T(L^{p'}_h(L^\infty_v))}
    \lesssim  \|a\|_{\widetilde{L}^\frac{2p'}{p'-2}_T(\widetilde{L}^{p'}_h(B^\frac{1}{2}_v))}
    \lesssim \sqrt{\frac{p'}{\beta}}\nu_h^{\frac{1}{p'}-\frac{1}{2}}
   \|a\|_{B^{-1+\frac{2}{p'},\frac{1}{2}}_{p'}(T)},\label{ANS-E2.2--2}
    \end{equation}
   where
 $\|a\|_{\widetilde{L}^\frac{2p'}{p'-2}_T(\widetilde{L}^{p'}_h(B^\frac{1}{2}_v))}=
 \displaystyle{\sum_{j\in\mathbb{Z}}}2^{\frac{j}{2}}
\|\Delta^v_ja\|_{{L}^\frac{2p'}{p'-2}_T({L}^{p'}_h(L^2_v))}$ and
 $\beta\in(0,\min\{p'-2,2\}]$.
Furthermore, for $p\geq2$, we have
    \begin{equation}
   \|\Delta^v_ja\|_{L^\frac{2p}{p-1}_T(L^{2p}_h(L^2_v))}
   \lesssim d_j\nu_h^{\frac{1}{2p}-\frac{1}{2}}2^{-\frac{j}{2}}
   \|a\|_{B^{-1+\frac{1}{p},\frac{1}{2}}_{2p}(T)}
   \lesssim d_j\nu_h^{\frac{1}{2p}-\frac{1}{2}}2^{-\frac{j}{2}}
   \|a\|_{B^{-1+\frac{2}{p},\frac{1}{2}}_{p}(T)}\label{ANS-E2.2--3}
   \end{equation}
  and
   \begin{equation}
    \|a\|_{L^\frac{2p}{p-1}_T(L^{2p}_h(L^\infty_v))}
    \lesssim  \|a\|_{\widetilde{L}^\frac{2p}{p-1}_T(\widetilde{L}^{2p}_h(B^\frac{1}{2}_v))}
    \lesssim\nu_h^{\frac{1}{2p}-\frac{1}{2}}
   \|a\|_{B^{-1+\frac{2}{p},\frac{1}{2}}_{p}(T)},\label{ANS-E2.2--4}
    \end{equation}
      where
 $\|a\|_{\widetilde{L}^\frac{2p}{p-1}_T(\widetilde{L}^{2p}_h(B^\frac{1}{2}_v))}=
 \displaystyle{\sum_{j\in\mathbb{Z}}}2^{\frac{j}{2}}
\|\Delta^v_ja\|_{{L}^\frac{2p}{p-1}_T({L}^{2p}_h(L^2_v))}$.
\end{lem}
\begin{rem}
  From now on, $A\lesssim B$ means $A\leq C(p) B$, where $C(p)$ is a
  constant depending on $p$.
\end{rem}
\noindent\textbf{Proof.} Let us first notice that
    $$
        \|\Delta^v_ja\|_{L^\frac{2p'}{p'-2}_T(L^{p'}_h(L^2_v))}^2=
        \|(\Delta^v_ja)^2\|_{L^\frac{p'}{p'-2}_T(L^\frac{p'}{2}_h(L^1_v))}.
    $$
Using Bony's decomposition in the horizontal variable, we have
    $$
        (\Delta^v_ja)^2=\sum_{k\in\mathbb{Z}}S^h_{k-1}\Delta^v_ja\Delta^h_{k}\Delta^v_ja
        +\sum_{k\in\mathbb{Z}}S^h_{k+2}\Delta^v_ja\Delta^h_{k}\Delta^v_ja.
    $$
These two terms are estimated exactly in the same way. Applying
H\"{o}lder's inequality and Lemma \ref{ANS-L2.1}, we obtain
    \begin{eqnarray*}
        &&\|S^h_{k-1}\Delta^v_ja\Delta^h_{k}\Delta^v_ja\|_{L^\frac{p'}{p'-2}_T(L^\frac{p'}{2}_h(L^1_v))}\\
        &\lesssim&\|S^h_{k-1}\Delta^v_ja\|_{L^\frac{2p'}{p'-2-\beta}_T(L^{p'}_h(L^2_v))}
        \|\Delta^h_{k}\Delta^v_ja\|_{L^\frac{2p'}{p'-2+\beta}_T(L^{p'}_h(L^2_v))}\\
                &\lesssim&\sum_{k'\leq k-2}\|\Delta^h_{k'}\Delta^v_ja\|_{L^\frac{2p'}{p'-2-\beta}_T(L^{p'}_h(L^2_v))}
                \|\Delta^h_{k}\Delta^v_ja\|_{L^\frac{2p'}{p'-2+\beta}_T(L^{p'}_h(L^2_v))}\\
        &\lesssim&\sum_{k'\leq
        k-2}\|\Delta^h_{k'}\Delta^v_ja\|_{L^2_T(L^{p'}_h(L^2_v))}^{1-\frac{2+\beta}{p'}}
        \|\Delta^h_{k'}\Delta^v_ja\|_{L^\infty_T(L^{p'}_h(L^2_v))}^{\frac{2+\beta}{p'}}
         \|\Delta^h_{k}\Delta^v_ja\|_{L^2_T(L^{p'}_h(L^2_v))}^{1-\frac{2-\beta}{p'}}
         \|\Delta^h_{k}\Delta^v_ja\|_{L^\infty_T(L^{p'}_h(L^2_v))}^\frac{2-\beta}{p'}\\
                           &\lesssim&\sum_{k'\leq
                            k-2}\left(2^{\frac{2+\beta}{p'}k'}
                            \|\Delta^h_{k'}\Delta^v_ja\|_{L^2_T(L^{p'}_h(L^2_v))}\right)^{1-\frac{2+\beta}{p'}}
                            \left(2^{(-1+\frac{2+\beta}{p'})k'}
                            \|\Delta^h_{k'}\Delta^v_ja\|_{L^\infty_T(L^{p'}_h(L^2_v))}\right)^{\frac{2+\beta}{p'}}\\
         &&\times\left(2^{\frac{2-\beta}{p'}k}
                             \|\Delta^h_{k}\Delta^v_ja\|_{L^2_T(L^{p'}_h(L^2_v))}\right)^{1-\frac{2-\beta}{p'}}
                       \left(2^{(-1+\frac{2-\beta}{p'})k}
                         \|\Delta^h_{k}\Delta^v_ja\|_{L^\infty_T(L^{p'}_h(L^2_v))}\right)^\frac{2-\beta}{p'}\\
        &\lesssim&\nu_h^{\frac{2}{p'}-1}\sum_{k'\leq
                            k-2}2^{(-1+\frac{2+\beta}{p'})k'}\left(\nu_h^\frac{1}{2}
                            \|\nabla_h\Delta^h_{k'}\Delta^v_ja\|_{L^2_T(L^{p'}_h(L^2_v))}+
                            \|\Delta^h_{k'}\Delta^v_ja\|_{L^\infty_T(L^{p'}_h(L^2_v))}\right)\\
         &&\times2^{(-1+\frac{2-\beta}{p'})k}\left(\nu_h^\frac{1}{2}
                             \|\nabla_h\Delta^h_{k}\Delta^v_ja\|_{L^2_T(L^{p'}_h(L^2_v))}+
                                                \|\Delta^h_{k}\Delta^v_ja\|_{L^\infty_T(L^{p'}_h(L^2_v))}\right).
    \end{eqnarray*}
From Corollary \ref{ANS-C1.2} and Lemma \ref{ANS-L2.2}, we get
        $$
            \|S^h_{k-1}\Delta^v_ja\Delta^h_{k}\Delta^v_ja\|_{L^\frac{p'}{p'-2}_T(L^\frac{p'}{2}_h(L^1_v))}
        \lesssim \frac{p'}{\beta}
        d^2_{k,j}\nu_h^{\frac{2}{p'}-1}2^{-j}\|a\|^2_{B^{-1+\frac{2}{p'},\frac{1}{2}}_{p'}(T)}.
        $$
Taking the sum over $k$, we obtain
        $$
        \|(\Delta^v_ja)^2\|_{L^\frac{p'}{p'-2}_T(L^\frac{p'}{2}_h(L^1_v))}
        \lesssim \frac{p'}{\beta}
        d^2_{j}\nu_h^{\frac{2}{p'}-1}2^{-j}\|a\|^2_{B^{-1+\frac{2}{p'},\frac{1}{2}}_{p'}(T)},
        $$
which is exactly the first inequality of this lemma. Combining it
with Lemma \ref{ANS-L2.1} and Corollary \ref{ANS-C1.1}, we can
immediately obtain (\ref{ANS-E2.2--2})-(\ref{ANS-E2.2--4}).
{\hfill $\square$\medskip}

\begin{lem}\label{ANS-L4.0}
 Let $a$ be in $B^{0,\frac{1}{2}}(T)$.
   Then, we have
        \begin{equation}
   \|\Delta^v_ja\|_{L^{q_1}_T(L^{q_2}_h(L^2_v))}
   \lesssim d_j\nu_h^{-\frac{1}{q_1} }2^{-\frac{j}{2}}
   \|a\|_{B^{0,\frac{1}{2}}(T)},
    \ \frac{1}{q_1}+\frac{1}{q_2}=\frac{1}{2},\
    q_2\in[2,4],\label{ANS-E2.6-00}
            \end{equation}
    $$
   \|\Delta^v_ja\|_{L^{2p}_T(L^\frac{2p}{p-1}_h(L^2_v))}
   \lesssim d_j\nu_h^{-\frac{1}{2p} }2^{-\frac{j}{2}}
   \|a\|_{B^{0,\frac{1}{2}}(T)}
   $$
  and
     $$
    \|a\|_{L^{2p}_T(L^\frac{2p}{p-1}_h(L^\infty_v))}
    \lesssim  \|a\|_{\widetilde{L}^{2p}_T(\widetilde{L}^\frac{2p}{p-1}_h(B^\frac{1}{2}_v))}
    \lesssim \nu_h^{-\frac{1}{2p}}
   \|a\|_{B^{0,\frac{1}{2}}(T)},
  $$
 where
 $\|a\|_{\widetilde{L}^{2p}_T(\widetilde{L}^\frac{2p}{p-1}_h(B^\frac{1}{2}_v))}=
 \displaystyle{\sum_{j\in\mathbb{Z}}}2^{\frac{j}{2}}\|\Delta^v_ja\|_{L^{2p}_T(L^\frac{2p}{p-1}_h(L^2_v))}$
 and $p\geq2$.
\end{lem}
\noindent\textbf{Proof.} From Corollary \ref{ANS-C1.1}, we have
    $$
    \|a\|_{B^{-\frac{1}{2},\frac{1}{2}}_4(T)}\lesssim
    \|a\|_{B^{0,\frac{1}{2}}(T)}.
    $$
From (\ref{ANS-E2.2--1}), we get
     $$
   \|\Delta^v_ja\|_{L^{4}_T(L^4_h(L^2_v))}
   \lesssim d_j\nu_h^{-\frac{1}{4} }2^{-\frac{j}{2}}
   \|a\|_{B^{-\frac{1}{2},\frac{1}{2}}_4(T)}
   \lesssim d_j\nu_h^{-\frac{1}{4} }2^{-\frac{j}{2}}
   \|a\|_{B^{0,\frac{1}{2}}(T)}.
   $$
Combining it with
    $$
    \|\Delta^v_ja\|_{L^{\infty}_T(L^2_h(L^2_v))}
   \lesssim d_j2^{-\frac{j}{2}}
   \|a\|_{B^{0,\frac{1}{2}}(T)},
    $$
 using interpolation, we obtain (\ref{ANS-E2.6-00}).
   Choosing $q_1=2p$, we can finish the proof of this lemma.
{\hfill $\square$\medskip}

Using Lemma \ref{ANS-L2.1}, we can obtain some estimates of $u_F$
in the following lemma.

\begin{lem}\label{ANS-L2.4}
  Let $u_0\in B^{-1+\frac{2}{p},\frac{1}{2}}_p$, $p\geq2$, and $u_F$ be as in
  (\ref{ANS-E1.9}), $1\leq q\leq \infty$. Then, there holds
            \begin{eqnarray}
            &&\|\Delta^h_k\Delta^v_l u_F\|_{L^q_T(L^p_h(L^2_v))}\nonumber\\
                    &\lesssim&
               d_{k,l}2^{(1-\frac{2}{p})k}
                2^{-\frac{l}{2}}
                \min(\nu_h^{-\frac{1}{q}}2^{-\frac{2k}{q}},
                \nu_3^{-\frac{1}{q}}2^{-\frac{2l}{q}})\|u_0\|_{B^{-1+\frac{2}{p},\frac{1}{2}}_p},
                \textrm{for}
                \ k\geq l-1,      \label{ANS-E2.2-2}
            \end{eqnarray}
    \begin{equation}
        \|\Delta^h_k\Delta^v_l u_F\|_{L^q_T(L^p_h(L^2_v))}=0,
              \  \textrm{for}
                \ k< l-1.
    \end{equation}
Moreover, $u_F$ belongs to
$B^{-1+\frac{2}{p},\frac{1}{2}}_p(\infty)\cap
C([0,\infty);B^{-1+\frac{2}{p},\frac{1}{2}}_p)$, and we have
    \begin{equation}
      \|u_F\|_{B^{-1+\frac{2}{p},\frac{1}{2}}_p(\infty)}\lesssim\|u_0\|_{B^{-1+\frac{2}{p},\frac{1}{2}}_p}.
      \label{ANS-E2.2-20}
    \end{equation}
\end{lem}
\noindent\textbf{Proof.}
  The relation (2.5) in \cite{Chemin07} tell us
        \begin{equation}
          \Delta^h_k\Delta^v_lu_F(t,x)
          =2^{2k+l}\int_{\mathbb{R}^3}g_h(t,2^k(x_h-y_h))
          g_v(t,2^l(x_3-y_3))\Delta^h_k\Delta^v_l u_0(y)dy.
            \label{ANS-E2.4}
        \end{equation}
with
    $$
          \|g_h(t,\cdot)\|_{L^1(\mathbb{R}^2)}\lesssim
          e^{-c\nu_ht2^{2k}},
        \  \|g_v(t,\cdot)\|_{L^1(\mathbb{R})}\lesssim
          e^{-c\nu_3t2^{2l}}.
    $$
  From Corollary \ref{ANS-C1.2} and (\ref{ANS-E2.4}), we have
        \begin{eqnarray*}
         \|\Delta^h_k\Delta^v_lu_F(t)\|_{L^p_h(L^2_v)}
        &\lesssim& e^{-c\nu_ht2^{2k}-c\nu_3t2^{2l}}
        \|\Delta^h_k\Delta^v_lu_0\|_{L^p_h(L^2_v)}\\
                &\lesssim&
                e^{-c\nu_ht2^{2k}-c\nu_3t2^{2l}}d_{k,l}2^{(1-\frac{2}{p})k}2^{-\frac{l}{2}}
                \|u_0\|_{B^{-1+\frac{2}{p},\frac{1}{2}}_p}.
        \end{eqnarray*}
By integration, we can obtain
(\ref{ANS-E2.2-2})-(\ref{ANS-E2.2-20}). The proof of  $u_F\in
C([0,\infty);B^{-1+\frac{2}{p},\frac{1}{2}}_p)$ is simple, and we
omit the details. {\hfill $\square$\medskip}

From Lemmas \ref{ANS-L2.1} and \ref{ANS-L2.4}, we can immediately
deduce the following corollary.

\begin{cor}\label{ANS-C2.3}
   For any $(q,p')\in
  [1,\infty]\times[p,\infty]$, $p\geq2$, we have
        $$
    \|\Delta^h_k u_F\|_{L^q(\mathbb{R}^+;L^{p'}_h(L^\infty_v))}
    \lesssim \nu_h^{-\frac{1}{q}}c_k2^{-k(\frac{2}{q}+\frac{2}{p'}-1)}
    \|u_0\|_{B^{-1+\frac{2}{p},\frac{1}{2}}_p}.
        $$
If $p'\in[p,\infty]$ and $q\in[1,\frac{2p'}{p'-2})$, we have
            $$
    \|\Delta^v_j u_F\|_{L^q(\mathbb{R}^+;L^{p'}_h(L^2_v))}
    \lesssim \nu_h^{-\frac{1}{q}}d_j2^{-j(\frac{2}{q}+\frac{2}{p'}-\frac{1}{2})}
    \|u_0\|_{B^{-1+\frac{2}{p},\frac{1}{2}}_p}.
        $$
\end{cor}

The following lemma is the end point of the second estimate of
Corollary \ref{ANS-C2.3}.

\begin{lem}\label{ANS-L2.5}
  Under the assumptions of Lemma \ref{ANS-L2.4}, we have
        $$
        \|\Delta^v_ju_F\|_{L^2(\mathbb{R}^+;L^\infty_h(L^2_v))}
        \lesssim d_j\nu_h^{-\frac{1}{2}}2^{-\frac{j}{2}}
        \|u_0\|_{B^{-1+\frac{2}{p},\frac{1}{2}}_p}
        $$
  and
        $$
        \|u_F\|_{L^2(\mathbb{R}^+;L^\infty(\mathbb{R}^3))}\lesssim
        \nu_h^{-\frac{1}{2}}
         \|u_0\|_{B^{-1+\frac{2}{p},\frac{1}{2}}_p}.
        $$
\end{lem}
\noindent\textbf{Proof.} Trivially, there holds
    $$
        \|\Delta^v_j u_F\|_{L^2_T(L^\infty_h(L^2_v))}^2
        =\|(\Delta^v_j u_F)^2\|_{L^1_T(L^\infty_h(L^1_v))}.
    $$
Using Bony's paradifferential decomposition in the horizontal
variables, we have
        \begin{equation}
        (\Delta^v_j u_F)^2
        =\sum_{k\in\mathbb{Z}}S^h_{k-1}\Delta^v_ju_F
        \Delta^h_{k}\Delta^v_ju_F
        +\sum_{k\in\mathbb{Z}}S^h_{k+2}\Delta^v_ju_F
        \Delta^h_{k}\Delta^v_ju_F.
        \label{ANS-E2.6}
        \end{equation}
Using Lemma \ref{ANS-L2.1} and H\"{o}lder's inequality, we get
    $$
    \|S^h_{k-1}\Delta^v_ju_F
        \Delta^h_{k}\Delta^v_ju_F\|_{L^1_T(L^\infty_h(L^1_v))}
        \lesssim
        2^{\frac{2}{p}k}\|S^h_{k-1}\Delta^v_ju_F\|_{L^\infty_T(L^{2p}_h(L^2_v))}
        \|\Delta^h_{k}\Delta^v_ju_F\|_{L^1_T(L^{2p}_h(L^2_v))}
    $$
By (\ref{ANS-E2.2-2}) and the proof of Lemma \ref{ANS-L2.2}, we
obtain
        $$
    \|S^h_{k-1}\Delta^v_ju_F\|_{L^\infty_T(L^{2p}_h(L^2_v))}\lesssim
    d_{k,j}2^{(1-\frac{1}{p})k}2^{-\frac{j}{2}}\|u_0\|_{B^{-1+\frac{2}{2p},\frac{1}{2}}_{2p}}.
        $$
Therefore, using (\ref{ANS-E2.2-2}) once again, we get
    $$
    \|\sum_{k\in\mathbb{Z}}S^h_{k-1}\Delta^v_ju_F
        \Delta^h_{k}\Delta^v_ju_F\|_{L^1_T(L^\infty_h(L^1_v))}\lesssim
        2^{-j}\nu_h^{-1}\left(
        \sum_{k\in\mathbb{Z}}d_{k,j}^2
        \right)\|u_0\|_{B^{-1+\frac{2}{2p},\frac{1}{2}}_{2p}}^2.
    $$
A similar argument yields a similar estimate for the other term in
(\ref{ANS-E2.6}). Then we deduce that
          $$
        \|\Delta^v_ju_F\|_{L^2(\mathbb{R}^+;L^\infty_h(L^2_v))}
        \lesssim d_j\nu_h^{-\frac{1}{2}}2^{-\frac{j}{2}}
        \|u_0\|_{B^{-1+\frac{2}{2p},\frac{1}{2}}_{2p}}
        \lesssim d_j\nu_h^{-\frac{1}{2}}2^{-\frac{j}{2}}
        \|u_0\|_{B^{-1+\frac{2}{p},\frac{1}{2}}_p}.
        $$
From Lemma \ref{ANS-L2.1}, we conclude that
    $$
    \|u_F\|_{L^2_T(L^\infty(\mathbb{R}^3))}\lesssim
    \sum_{j\in\mathbb{Z}}2^{\frac{j}{2}}\|\Delta^v_{j}u_F\|_{L^2_T(L^\infty_h(L^2_v))}
    \lesssim \nu_h^{-\frac{1}{2}}\|u_0\|_{B^{-1+\frac{2}{p},\frac{1}{2}}_p}.
    $$
{\hfill $\square$\medskip}

\section{The proof of an existence theorem}\label{ANS-Sec3}
The purpose of this section is to prove the following existence
theorem.

\begin{thm}\label{ANS-T3.1}
 A sufficiently small constant $c$ exists which satisfies the following property: if the divergence free vector
field $u_0\in
  B^{-1+\frac{2}{p},\frac{1}{2}}_p$, $p\geq2$ and $[u_0]_{E^p_\infty}\leq c\nu_h$,
  then the system (\ref{ANS-E1.1})  with initial data $u_0$ has a
  global solution in the space $\{u_F+
  B^{0,\frac{1}{2}}(\infty)\}\cap C([0,\infty);B^{-1+\frac{2}{p},\frac{1}{2}}_p)$.
\end{thm}
\noindent\textbf{Proof.}
  As announced in the introduction, we shall look for a solution of
   the form
            $$
        u=u_F+w.
            $$
 Actually, by
  substituting the above formula  to (\ref{ANS-E1.1}), we get
        \begin{equation}
            \left\{
        \begin{array}{l}
        w_t+w\cdot\nabla w-\nu_h\Delta_h w-\nu_3\partial_3^2w
        +w\cdot\nabla u_F+u_F\cdot\nabla w=-u_F\cdot\nabla
        u_F-\nabla P,\\
                \mathrm{div}w=0,\\
                    w|_{t=0}=u_{0ll}=u_0-u_{0hh},
        \end{array}
            \right.\label{ANS-E3.1}
        \end{equation}
where
    \begin{equation}
      u_{0ll}=\sum_{j\in\mathbb{Z}}S^h_{j-1}\Delta^v_ju_0.
    \end{equation}
Moreover, we obtain
        $
        \|u_{0ll}\|_{B^{0,\frac{1}{2}}}\lesssim
        \|u_0\|_{B^{-1+\frac{2}{p},\frac{1}{2}}_p}.
        $
We shall use the classical Friedrichs' regularization method to
construct the approximate solutions to (\ref{ANS-E3.1}). For
simplicity, we just outline it here (for the details, see
\cite{Chemin07,Chemin06,Paicu05}). In order to do so, let us
define the sequence of operators $(P_n)_{n\in\mathbb{N}}$,
$(P_{1n})_{n\in\mathbb{N}}$ and $(P_{2n})_{n\in\mathbb{N}}$  by
    $$
        P_na:=\mathcal{F}^{-1}(1_{B(0,n)}\hat{a}),
        \ P_{1n}a:=\mathcal{F}^{-1}
        (1_{\{|\xi|\leq n,\ |\xi_3|
    \geq\frac{1}{n}\}}\hat{a}),
    \ P_{2n}a:=\mathcal{F}^{-1}
        (1_{\{|\xi_3|<\frac{1}{n}\}}\hat{a}),
    $$
and we define the following approximate system:
      \begin{equation}
            \left\{
        \begin{array}{l}
        \partial_tw_n+P_n(w_n\cdot\nabla w_n)-\nu_h\Delta_h
        w_n-\nu_3\partial_3^2w_n
        +P_n(w_n\cdot\nabla u_F)+P_n(u_F\cdot\nabla w_n)\\
            =-P_{1n}(u_F\cdot\nabla
        u_F)-P_n\nabla(-\Delta)^{-1}\partial_j\partial_k\left(
        (u^j_F+w^j_n)(u^k_F+w^k_n)-P_{2n}(u_F^j
        u_F^k)
        \right),\\
                \mathrm{div}w_n=0,\\
                    w_n|_{t=0}=P_n(u_{0ll}),
        \end{array}
            \right.\label{ANS-E3.3}
        \end{equation}
where $(-\Delta)^{-1}\partial_j\partial_k$ is defined precisely by
    $$(-\Delta)^{-1}\partial_j\partial_ka
    :=\mathcal{F}^{-1}(|\xi|^{-2}\xi_j\xi_k\hat{a}).$$
Then, the system (\ref{ANS-E3.1}) appears to be an ordinary
differential equation in the space
    $$
        L^2_n:=\{
        a\in L^2(\mathbb{R}^3)\big|\mbox{div}a=0,\ \mathrm{Supp} \hat{a}\subset
        B(0,n)\}.
    $$
Such system is globally wellposed because
    $$
        \frac{d}{dt}\|w_n\|_{L^2}^2\leq
        C_n\|u_F\|_{L^\infty}\|w_n\|_{L^2}^2
    +C_n\|u_F\cdot\nabla u_F\|_{B^{0,\frac{1}{2}}}\|w_n\|_{L^2},
    $$
and $u_F$ belongs to $L^2(\mathbb{R}^+;L^\infty(\mathbb{R}^3))$.

Now, the proof of Theorem \ref{ANS-T3.1} reduces to the following
three propositions, which we shall admit for the time begin.

\begin{prop}\label{ANS-P4.1}
  Let $a$  be a divergence free vector filed in
  $B^{0,\frac{1}{2}}(T)$ and $u$ be a divergence free vector field in $B^{-1+\frac{2}{p},
        \frac{1}{2}}_p(T)$. Then, for any $j\in\mathbb{Z}$, we have
        $$
            F_j(T):=\int^T_0
            \left|\int_{\mathbb{R}^3}
            \Delta_j^v(u\cdot\nabla a)
            \Delta^v_jadx
            \right|dt\lesssim
            d_j^2\nu_h^{-\frac{1}{2}-\frac{1}{2p}}2^{-j}\|a\|^2_{B^{0,\frac{1}{2}}(T)}
            \|u \|_{\widetilde{L}^\frac{2p}{p-1}_T(\widetilde{L}^{2p}_h(B^\frac{1}{2}_v))}.
        $$
\end{prop}

\begin{prop}\label{ANS-P4.2}
  Let $a$ and $b$ be two divergence free vector fields in
  $B^{0,\frac{1}{2}}(T)$. Then, for any $j\in\mathbb{Z}$, we have
        \begin{eqnarray*}
            &&G_j(T):=\int^T_0
            \left|\int_{\mathbb{R}^3}
            \Delta_j^v(a\cdot\nabla u_F)
            \Delta^v_jbdx
            \right|dt\\
                    &\lesssim&
            d_j^22^{-j}\|a\|_{B^{0,\frac{1}{2}}(T)}
            \|b\|_{B^{0,\frac{1}{2}}(T)}\left(\nu_h^{-\frac{1}{2}-\frac{1}{2p}}
            \|u_F\|_{\widetilde{L}^\frac{2p}{p-1}_T(\widetilde{L}^{2p}_h(B^\frac{1}{2}_v))}
            +\|u_F\|_{\widetilde{L}^1_T(\widetilde{L}^\infty_h(B^\frac{3}{2}_v))}\right),
        \end{eqnarray*}
where
$\|u_F\|_{\widetilde{L}^1_T(\widetilde{L}^\infty_h(B^\frac{3}{2}_v))}:=\sum_{j\in\mathbb{Z}}
2^{\frac{3}{2}j}\|\Delta^v_ju_F\|_{L^1_T(L^\infty_h(L^2_v))}$.
\end{prop}

From Corollary \ref{ANS-C1.1}, Lemma \ref{ANS-L2.3} and
Proposition \ref{ANS-P4.1}, we have the following proposition.
\begin{prop}\label{ANS-P3.3}
  Let $a$ and $b$ be two divergence free vector fields in
  $B^{0,\frac{1}{2}}(T)$. Then, for any $j\in\mathbb{Z}$, we have
        $$
            \int^T_0
            \left|\int_{\mathbb{R}^3}
            \Delta_j^v(a\cdot\nabla b)
            \Delta^v_jbdx
            \right|dt\lesssim
            d_j^2\nu_h^{-1}2^{-j}\|a\|_{B^{0,\frac{1}{2}}(T)}\|b\|^2_{B^{0,\frac{1}{2}}(T)}.
        $$
\end{prop}

\noindent\textbf{Conclusion of the proof of Theorem
\ref{ANS-T3.1}.} Applying the operator $\Delta^v_j$ to
(\ref{ANS-E3.3}) and taking the $L^2$ inner product of the
resulting equation with $\Delta^v_j w_n$, we have
    \begin{eqnarray*}
      &&\frac{d}{dt}\|\Delta^v_jw_n\|^2_{L^2}+2\nu_h\|\nabla_h\Delta^v_jw_n\|^2_{L^2}
      +2\nu_3\|\partial_3\Delta^v_jw_n\|^2_{L^2}\\
            &=&-2\int_{\mathbb{R}^3}\Delta^v_j(w_n\cdot\nabla
            w_n)\Delta^v_jw_ndx-2\int_{\mathbb{R}^3}\Delta^v_j(u_F\cdot\nabla
            w_n)\Delta^v_jw_ndx\\
            &&-2\int_{\mathbb{R}^3}\Delta^v_j(w_n\cdot\nabla
            u_F)\Delta^v_jw_ndx-2\int_{\mathbb{R}^3}\Delta^v_jP_{1n}(u_F\cdot\nabla
            u_F)\Delta^v_jw_ndx.
    \end{eqnarray*}
From Lemmas \ref{ANS-L2.3}, \ref{ANS-L2.4}, Corollary
\ref{ANS-C2.3} and Propositions \ref{ANS-P4.1}-\ref{ANS-P3.3}, we
get
        \begin{eqnarray*}
          &&2^j\left(
            \|\Delta^v_jw_n\|_{L^\infty_T(L^2)}^2+{2\nu_h}\|\nabla_h\Delta^v_jw_n\|_{L^2_T(L^2)}^2
            +{2\nu_3}\|\partial_3\Delta^v_jw_n\|_{L^2_T(L^2)}^2
          \right)\\
                &\leq&2^j
            \|\Delta^v_jw_n(0)\|^2_{L^2}+
            Cd_j^2\nu_h^{-1}\| w_n\|_{B^{0,\frac{1}{2}}(T)}^2\left(
            \| w_n\|_{B^{0,\frac{1}{2}}(T)}+\|u_0\|_{B^{-1+\frac{2}{p},\frac{1}{2}}_p}
            \right)\\
                &&+C2^j\|\Delta^v_j(u_F\cdot\nabla
            u_F)\|_{L^1_T(L^2)}\|\Delta^v_jw_n\|_{L^\infty_T(L^2)}
        \end{eqnarray*}
and
          \begin{eqnarray*}
          &&
            \| w_n\|_{B^{0,\frac{1}{2}}(T)}\\
                &\leq&2C_0[ u_0]_{E^p_\infty}+
            C\nu_h^{-\frac{1}{2}}\| w_n\|_{B^{0,\frac{1}{2}}(T)}\left(
            \| w_n\|_{B^{0,\frac{1}{2}}(T)} +[u_0]_{E^p_\infty}
            \right)^\frac{1}{2}.
        \end{eqnarray*}
Then, we have
    $$
        \| w_n\|_{B^{0,\frac{1}{2}}(T)}
                 \leq2C_0[
                 u_0]_{E^p_\infty}+\frac{4CC_0\sqrt{4C_0+1}}{\sqrt{\nu_h}}[
                 u_0]_{E^p_\infty}^\frac{3}{2},
    $$
for all $T<T_n:=\sup\{t>0; \| w_n\|_{B^{0,\frac{1}{2}}(t)}\leq
4C_0[ u_0]_{E^p_\infty}\}$. Then, if $[ u_0]_{E^p_\infty}$ is
small enough with respect to $\nu_h$, we get for any $n$ and for
any $T<T_n$,
        $$
        \| w_n\|_{B^{0,\frac{1}{2}}(T)}
                 \leq \frac{5}{2}C_0[ u_0]_{E^p_\infty}.
        $$
Thus, $T_n=+\infty$. Then, the existence follows from classical
compactness method, the details of which are omitted (see
\cite{Chemin06,Paicu05}).

In order to prove the continuity of the solution $u$, we have to
prove the continuity of $w$. From (\ref{ANS-E3.1}), we have
    \begin{eqnarray*}
    \Delta^v_j w_t&=&\nu_h\Delta^v_j\Delta_h w+\nu_3\Delta^v_j
    \partial_3^2w-\Delta^v_j(w\cdot\nabla w)\\
            &&
        -\Delta^v_j(w\cdot\nabla u_F)-\Delta^v_j(u_F\cdot\nabla w)-\Delta^v_j(u_F\cdot\nabla
        u_F)-\Delta^v_j\nabla P.
    \end{eqnarray*}
We can easily obtain that for any $T>0$ and $j\in \mathbb{Z}$,
    $$
  \nu_3\Delta^v_j
    \partial_3^2w \in L^\infty([0,T];L^2),\  \nu_h\Delta^v_j\Delta_h w\in
    L^2(0,T;L^2_v(\dot{H}^{-1}_h)),
    $$
    $$
    \Delta^v_j(u_F\cdot\nabla
        u_F)\in L^1([0,T];L^2)
    $$
  and
    $$(\nu_h\Delta^v_j\Delta_h w+\nu_3\Delta^v_j
    \partial_3^2w-\Delta^v_j(u_F\cdot\nabla
        u_F)|\Delta_j^vw)_{L^2}\in L^1([0,T]).$$
 From Propositions \ref{ANS-P4.1}-\ref{ANS-P3.3}, we have
    $$(\Delta^v_j(w\cdot\nabla w)        +\Delta^v_j(w\cdot\nabla u_F)
    +\Delta^v_j(u_F\cdot\nabla w)
        |\Delta_j^vw)_{L^2}\in L^1([0,T]).$$
Thus, we have $\frac{d}{dt}\|\Delta^v_j w(t)\|^2_{L^2}\in
L^1([0,T])$, for any $T>0$ and $j\in \mathbb{Z}$. Combining it
with $w\in B^{0,\frac{1}{2}}(\infty)$, we can easily get $w\in
C([0,\infty);B^{0,\frac{1}{2}})$. Then Theorem \ref{ANS-T3.1} is
proved provided of course that we have proved Propositions
\ref{ANS-P4.1}-\ref{ANS-P4.2}. {\hfill $\square$\medskip}

To prove Propositions \ref{ANS-P4.1}-\ref{ANS-P4.2}, we need the
following two lemmas.

\begin{lem}\label{ANS-L4.1}
  Let a be in $B^{0,\frac{1}{2}}(T)$ and $u$ be in $B^{-1+\frac{2}{p},\frac{1}{2}}_p(T)$. We have
    $$
        \|\Delta^v_j(u\partial_h a)\|_{L^\frac{2p}{2p-1}_T(L^\frac{2p}{p+1}_h(L^2_v))}
        \lesssim
        d_j \nu_h^{-\frac{1}{2}}
        2^{-\frac{j}{2}}\|a\|_{B^{0,\frac{1}{2}}(T)}
        \|u\|_{\widetilde{L}^\frac{2p}{p-1}_T(\widetilde{L}^{2p}_h(B^\frac{1}{2}_v))}.
    $$
\end{lem}
\noindent\textbf{Proof.}
  Using Bony's decomposition in the vertical variable, we obtain
        $$
        \Delta^v_j(u\partial_ha)=\sum_{|j-j'|\leq5}\Delta^v_j(S^v_{j'-1}u\partial_h\Delta^v_{j'} a)
        +\sum_{j'\geq j-N_0}\Delta^v_j(\Delta^v_{j'}u\partial_h S^v_{j'+2}a
        ).
        $$
  Using H\"{o}lder's inequality and Lemma \ref{ANS-L2.3}, we get
        \begin{eqnarray*}
          \|\Delta^v_j(S^v_{j'-1}u\partial_h\Delta^v_{j'} a)\|_{L^\frac{2p}{2p-1}_T(L^\frac{2p}{p+1}_h(L^2_v))}
          &\lesssim& \|S^v_{j'-1}u\|_{L^\frac{2p}{p-1}_T(L^{2p}_h(L^\infty_v))}
          \|\Delta^v_{j'} \partial_ha\|_{L^2_T(L^2(\mathbb{R}^3))}\\
                &\lesssim&
        d_{j'}\nu_h^{-\frac{1}{2}}2^{-\frac{j'}{2}}\|a
        \|_{B^{0,\frac{1}{2}}(T)}\|u\|_{\widetilde{L}^\frac{2p}{p-1}_T(\widetilde{L}^{2p}_h(B^\frac{1}{2}_v))}
        \end{eqnarray*}
and
            \begin{eqnarray*}
          \|\Delta^v_j(\Delta^v_{j'}u\partial_h S^v_{j'+2}a
        )\|_{L^\frac{2p}{2p-1}_T(L^\frac{2p}{p+1}_h(L^2_v))}
          &\lesssim& \|S^v_{j'+2}(\partial_ha)\|_{L^2_T(L^2_h(L^\infty_v))}
          \|\Delta^v_{j'} u\|_{L^\frac{2p}{p-1}_T(L^{2p}_h(L^2_v))}\\
                &\lesssim&
        d_{j'}\nu_h^{-\frac{1}{2}}2^{-\frac{j'}{2}}
        \|a\|_{B^{0,\frac{1}{2}}(T)}\|u \|_{\widetilde{L}^\frac{2p}{p-1}_T(\widetilde{L}^{2p}_h(B^\frac{1}{2}_v))}.
        \end{eqnarray*}
Then, we can immediately finish the proof. {\hfill
$\square$\medskip}

\begin{lem}\label{ANS-L4.2}
 Let a be in $B^{0,\frac{1}{2}}(T)$ and $u$ be in $B^{-1+\frac{2}{p},\frac{1}{2}}_p(T)$.  We have
    $$
        \|\Delta^v_j(u a)\|_{L^2_T(L^2(\mathbb{R}^3))}
        \lesssim
        d_j\nu_h^{-\frac{1}{2p}}2^{-\frac{j}{2}}
        \|a\|_{B^{0,\frac{1}{2}}(T)}\|u\|_{\widetilde{L}^\frac{2p}{p-1}_T(\widetilde{L}^{2p}_h(B^\frac{1}{2}_v))}.
    $$
\end{lem}
\noindent\textbf{Proof.}
  Using Bony's decomposition in the vertical variable, we obtain
        $$
        \Delta^v_j(u a)=\sum_{|j-j'|\leq5}\Delta^v_j(S^v_{j'-1}u\Delta^v_{j'}  a)
        +\sum_{j'\geq j-N_0}\Delta^v_j(S^v_{j'+2} a \Delta^v_{j'}u
        ).
        $$
  Using H\"{o}lder's inequality, Lemmas \ref{ANS-L2.3}-\ref{ANS-L4.0}, we get
        \begin{eqnarray*}
          \|\Delta^v_j(S^v_{j'-1}u\Delta^v_{j'}a)\|_{L^2_T(L^2(\mathbb{R}^3))}
          &\lesssim& \|S^v_{j'-1}u\|_{L^\frac{2p}{p-1}_T(L^{2p}_h(L^\infty_v))}
          \|\Delta^v_{j'}  a\|_{L^{2p}_T(L^\frac{2p}{p-1}_h(L^2_v))}\\
                &\lesssim&
        d_{j'}\nu_h^{-\frac{1}{2p}}2^{-\frac{j'}{2}}\|a\|_{B^{0,\frac{1}{2}}(T)}\|u
        \|_{\widetilde{L}^\frac{2p}{p-1}_T(\widetilde{L}^{2p}_h(B^\frac{1}{2}_v))}
        \end{eqnarray*}
and
            \begin{eqnarray*}
          \|\Delta^v_j(S^v_{j'+2}a\Delta^v_{j'}u)\|_{L^2_T(L^2(\mathbb{R}^3))}
          &\lesssim& \|S^v_{j'+2}a\|_{L^{2p}_T(L^\frac{2p}{p-1}_h(L^\infty_v))}
          \|\Delta^v_{j'} u\|_{L^\frac{2p}{p-1}_T(L^{2p}_h(L^2_v))}\\
                &\lesssim&
        d_{j'}\nu_h^{-\frac{1}{2p}}
        2^{-\frac{j'}{2}}\|a\|_{B^{0,\frac{1}{2}}(T)}\|
        u\|_{\widetilde{L}^\frac{2p}{p-1}_T(\widetilde{L}^{2p}_h(B^\frac{1}{2}_v))}.
        \end{eqnarray*}
Then, we can immediately finish the proof. {\hfill
$\square$\medskip}

\noindent\textbf{Proof of Proposition \ref{ANS-P4.1}.} We
distinguish the terms with horizontal derivatives from the terms
with vertical ones, writing
    $$F_j(T)\leq F_j^h(T)+F_j^v(T),$$
where
    $$F_j^h(T):=\int^T_0\left|\int_{\mathbb{R}^3}
    \Delta^v_j(u^h\cdot\nabla_ha)\Delta^v_jadx
    \right|dt$$
and
     $$F_j^v(T):=\int^T_0\left|\int_{\mathbb{R}^3}
    \Delta^v_j(u^3\partial_3a)\Delta^v_jadx
    \right|dt.$$
Using H\"{o}lder's inequality, Lemmas \ref{ANS-L4.0} and
\ref{ANS-L4.1}, we obtain
    \begin{eqnarray*}
      F_j^h(T)&\leq& \|\Delta^v_j(u^h\cdot\nabla_h a)\|_{L^\frac{2p}{2p-1}_T(L^\frac{2p}{p+1}_h(L^2_v))}
      \|\Delta^v_ja\|_{L^{2p}_T(L^\frac{2p}{p-1}_h(L^2_v))}\\
            &\lesssim&d_j^2\nu_h^{-\frac{1}{2}-\frac{1}{2p}}2^{-j}
            \|a\|^2_{B^{0,\frac{1}{2}}(T)}
            \|u\|_{\widetilde{L}^\frac{2p}{p-1}_T(\widetilde{L}^{2p}_h(B^\frac{1}{2}_v))}.
    \end{eqnarray*}
Applying the trick from \cite{Chemin07,Paicu05}, using
paradifferential decomposition in the vertical variable to
$\Delta^v_j(u^3\partial_3a)$ first, then by a commutator process,
one get
        \begin{eqnarray*}
        \Delta^v_j(u^3\partial_3a)&=&S^v_{j-1}u^3\partial_3\Delta^v_ja+\sum_{|j-l|\leq5}[\Delta^v_j;
        S^v_{l-1}u^3]\partial_3\Delta^v_la\\
                    &&+\sum_{|j-l|\leq5}(S^v_{l-1}u^3-S^v_{j-1}u^3)\partial_3\Delta^v_j\Delta_l^va+
                    \sum_{l\geq
                    j-N_0}\Delta^v_j(\Delta^v_lu^3\partial_3S^v_{l+2}a).
    \end{eqnarray*}
Correspondingly, we decompose $F_j^v(T)$ as
    $$
        F_j^v(T):=F_j^{1,v}(T)+F_j^{2,v}(T)+F_j^{3,v}(T)+F_j^{4,v}(T).
    $$
Using integration by parts and the fact that div$u$=0, we have
    $$
        F_j^{1,v}(T)=\frac{1}{2}\int^T_0\left|\int_{\mathbb{R}^3} S^v_{j-1}\mathrm{div}_h
        u^h|\Delta^v_j a|^2dx\right|dt=\int^T_0\left|\int_{\mathbb{R}^3} S^v_{j-1}
        u^h\cdot\nabla_h \Delta^v_j a\Delta^v_j adx\right|dt.
    $$
From Lemmas \ref{ANS-L2.3}-\ref{ANS-L4.0} and H\"{o}lder's
inequality, we get
    \begin{eqnarray*}
            F_j^{1,v}(T)&\leq& \|S^v_{j-1}u\|_{L^\frac{2p}{p-1}_T(L^{2p}_h(L^\infty_v))}
            \|\Delta^v_ja\|_{L^{2p}_T(L^\frac{2p}{p-1}_h(L^2_v))}
            \|\nabla_h\Delta^v_ja\|_{L^2_T(L^2(\mathbb{R}^3))}\\
                    &\lesssim&d_j^2\nu_h^{-\frac{1}{2}-\frac{1}{2p}}2^{-j}\|u\|_{\widetilde{L}^\frac{2p}{p-1}_T(
                    \widetilde{L}^{2p}_h(B^\frac{1}{2}_v))}
                    \|a\|^2_{B^{0,\frac{1}{2}}(T)}.
    \end{eqnarray*}

To deal with the commutator in $F^{2,v}_j$, we first use the
Taylor formula to get
    \begin{eqnarray*}
      F^{2,v}_j(T)&=&\sum_{|j-l|\leq5}\int^T_0\left|\int_{\mathbb{R}^3}
      2^j\int_\mathbb{R}h(2^j(x_3-y_3))\int^1_0S^v_{l-1}\partial_3u^3(x_h,\tau
      y_3+(1-\tau)x_3)d\tau\right.\\
                &&\left.\times
                (y_3-x_3)\partial_3\Delta^v_la(x_h,y_3)dy_3\Delta^v_ja(x)dx
      \right|dt.
    \end{eqnarray*}
Using div$u=0$ and integration by parts, we have
    \begin{eqnarray*}
      F^{2,v}_j(T)&=&\sum_{|j-l|\leq5}\int^T_0\left|\int_{\mathbb{R}^3}
      \int_\mathbb{R}\bar{h}(2^j(x_3-y_3))\int^1_0S^v_{l-1}u^h(x_h,\tau
      y_3+(1-\tau)x_3)d\tau\right.\\
                &&\left.\cdot\nabla_h
                \partial_3\Delta^v_la(x_h,y_3)dy_3\Delta^v_ja(x)dx
      \right|dt\\
                &&+\sum_{|j-l|\leq5}\int^T_0\left|\int_{\mathbb{R}^3}
      \int_\mathbb{R}\bar{h}(2^j(x_3-y_3))\int^1_0S^v_{l-1}u^h(x_h,\tau
      y_3+(1-\tau)x_3)d\tau\right.\\
                &&\left.\times
                \partial_3\Delta^v_la(x_h,y_3)dy_3\nabla_h\Delta^v_ja(x)dx
      \right|dt,
    \end{eqnarray*}
where $\bar{h}(x_3)=x_3h(x_3)$. Using H\"{o}lder's inequality,
Young's inequality and Lemma \ref{ANS-L2.5}, we obtain
    \begin{eqnarray*}
      F^{2,v}_j(T)&\lesssim&\sum_{|j-l|\leq5}2^{l-j}\|S^v_{l-1}
      u^h\|_{L^\frac{2p}{p-1}_T(L^{2p}_h(L^\infty_v))}
      \|\nabla_h\Delta^v_la\|_{L^2_T(L^2(\mathbb{R}^3))}
            \|\Delta^v_ja\|_{L^{2p}_T(L^\frac{2p}{p-1}_h(L^2_v))}\\
                    &&+\sum_{|j-l|\leq5}2^{l-j}\|S^v_{l-1}
      u^h   \|_{L^\frac{2p}{p-1}_T(L^{2p}_h(L^\infty_v))}
      \|\nabla_h\Delta^v_ja\|_{L^2_T(L^2(\mathbb{R}^3))}
            \|\Delta^v_la\|_{L^{2p}_T(L^\frac{2p}{p-1}_h(L^2_v))}\\
      &\lesssim&d_j^2\nu_h^{-\frac{1}{2}-\frac{1}{2p}
      }2^{-j}\|u\|_{\widetilde{L}^\frac{2p}{p-1}_T(\widetilde{L}^{2p}_h(B^\frac{1}{2}_v))}
      \|a\|_{
      B^{0,\frac{1}{2}}(T)}^2.
    \end{eqnarray*}

It is easy to see that
    $$
        F_j^{3,v}(T)\leq\sumetage{|j-l'|\leq5}{|j-l|\leq5}\int^T_0\left|
        \int_{\mathbb{R}^3}\Delta^v_{l'}u^3\partial_3\Delta_j^v\Delta_{l}^va\Delta^v_jadx
        \right|dt.
    $$
We can rewrite $\Delta^v_{l'}u^3$ as following:
   \begin{equation}\label{ANS-E4.1}
        \Delta^v_{l'}u^3
        =\int_{\mathbb{R}}g^v(2^{l'}(x_3-y_3))\partial_3\Delta^v_{l'}u^3(x_h,y_3)dy_3,
    \end{equation}
where $g^v\in\mathcal{S}(\mathbb{R})$ satisfying
$\mathcal{F}(g^v)(\xi_3)=\frac{\tilde{\varphi}(|\xi_3|)}{i\xi_3}$.
Using div$u=0$, integration by parts,  Young's inequality and
Lemma \ref{ANS-L2.5}, we get
    \begin{eqnarray*}
        F^{3,v}_j(T)&\leq&\sumetage{|j-l'|\leq5}{|j-l|\leq5}\int^T_0\left|
        \int_{\mathbb{R}^3}\int_{\mathbb{R}}g^v(2^{l'}(x_3-y_3))
        \Delta^v_{l'}u^h(x_h,y_3)dy_3\cdot\nabla_h\partial_3\Delta_j^v\Delta_{l}^va\Delta^v_jadx
        \right|dt\\
                &&+\sumetage{|j-l'|\leq5}{|j-l|\leq5}\int^T_0\left|
        \int_{\mathbb{R}^3}\int_{\mathbb{R}}g^v(2^{l'}(x_3-y_3))
        \Delta^v_{l'}u^h(x_h,y_3)dy_3\cdot\nabla_h\Delta^v_ja\partial_3\Delta_j^v\Delta_{l}^vadx
        \right|dt\\
                &\lesssim&\sumetage{|j-l'|\leq5}{|j-l|\leq5}2^{l-l'}
                \|\Delta^v_{l'}u\|_{L^\frac{2p}{p-1}_T(L^{2p}_h(L^\infty_v))}
                \|\nabla_h\Delta_j^va\|_{L^2_T(L^2(\mathbb{R}^3))}
                                    \|\Delta_j^va\|_{L^{2p}_T(L^\frac{2p}{p-1}_h(L^2_v))}\\
        &\lesssim&d_j^2\nu_h^{-\frac{1}{2}-\frac{1}{2p}
        }2^{-j}\|u\|_{\widetilde{L}^\frac{2p}{p-1}_T(\widetilde{L}^{2p}_h(B^\frac{1}{2}_v))}\|a\|_{
      B^{0,\frac{1}{2}}(T)}^2.
    \end{eqnarray*}

    Similarly, we have
        \begin{eqnarray*}
       && F^{4,v}_j(T)\\&\leq&\sum_{l\geq
        j-N_0}\int^T_0\left|
        \int_{\mathbb{R}^3}\Delta^v_j\left(\int_{\mathbb{R}}g^v(2^{l }(x_3-y_3))
        \Delta^v_{l }u^h(x_h,y_3)dy_3\cdot\nabla_h\partial_3S_{l+2}^va\right)\Delta^v_jadx
        \right|dt\\
                &&+\sum_{l\geq
        j-N_0}\int^T_0\left|
        \int_{\mathbb{R}^3}\Delta^v_j\left(\int_{\mathbb{R}}g^v(2^{l}(x_3-y_3))
        \Delta^v_{l }u^h(x_h,y_3)dy_3\partial_3S_{l+2}^va\right)\cdot\nabla_h\Delta^v_jadx
        \right|dt\\
                &\lesssim&\sum_{l\geq
        j-N_0}
                \|\Delta^v_{l }u\|_{L^\frac{2p}{p-1}_T(L^{2p}_h(L^2_v))}
                \|\nabla_hS_{l+2}^va\|_{L^2_T(L^2_h(L^\infty_v))}
                                    \|\Delta_j^va\|_{L^{2p}_T(L^\frac{2p}{p-1}_h(L^2_v))}\\
            &&+\sum_{l\geq
        j-N_0}
                \|\Delta^v_{l }u\|_{L^\frac{2p}{p-1}_T(L^{2p}_h(L^2_v))}
                \| S_{l+2}^va\|_{L^{2p}_T(L^\frac{2p}{p-1}_h(L^\infty_v))}
                                    \|\nabla_h\Delta_j^va\|_{L^2_T(L^2_h(L^2_v))}\\
        &\lesssim&d_j^2\nu_h^{-\frac{1}{2}-\frac{1}{2p}
        }2^{-j}\|u\|_{\widetilde{L}^\frac{2p}{p-1}_T(\widetilde{L}^{2p}_h(B^\frac{1}{2}_v))}
        \|a\|_{
      B^{0,\frac{1}{2}}(T)}^2.
    \end{eqnarray*}
 This completes the proof of Proposition \ref{ANS-P4.1}.   {\hfill
$\square$\medskip}

\noindent\textbf{Proof of Proposition \ref{ANS-P4.2}.}
   We
distinguish the terms with horizontal derivatives from the terms
with vertical ones, writing
    $$
        G_j(T)\leq G^h_j(T)+G^v_j(T),
    $$
where
        $$
            G_j^h(T):=\int^T_0
            \left|\int_{\mathbb{R}^3}
            \Delta_j^v(a^h\cdot\nabla_h u_F)
            \Delta^v_jbdx\right|dt
        $$
and
    $$
            G_j^v(T):=\int^T_0
            \left|\int_{\mathbb{R}^3}
            \Delta_j^v(a^3\partial_3 u_F)
            \Delta^v_jbdx\right|dt.
        $$
Using integration by parts, we have
    $$
        \int_{\mathbb{R}^3}
            \Delta_j^v(a^h\cdot\nabla_h u_F)
            \Delta^v_jbdx
            =-\int_{\mathbb{R}^3}
            \Delta_j^v(u_F\mathrm{div}_ha^h)
            \Delta^v_jbdx-\int_{\mathbb{R}^3}
            \Delta_j^v(a^h\otimes u_F):
            \nabla_h\Delta^v_jbdx.
    $$
From Lemmas \ref{ANS-L4.0} and \ref{ANS-L4.1}-\ref{ANS-L4.2}, we
have
    \begin{eqnarray*}
      \int^T_0\left|
      \int_{\mathbb{R}^3}
            \Delta_j^v(u_F\mathrm{div}_ha^h)
            \Delta^v_jbdx
      \right|dt
            &\leq&\|\Delta_j^v(u_F\mathrm{div}_ha^h)\|_{L^\frac{2p}{2p-1}_T(L^\frac{2p}{p+1}_h(L^2_v))}
      \|\Delta^v_jb\|_{L^{2p}_T(L^\frac{2p}{p-1}_h(L^2_v))}\\
            &\lesssim&
            d_j^2\nu_h^{-\frac{1}{2}-\frac{1}{2p}}2^{-j}\|a\|_{B^{0,\frac{1}{2}}(T)}
            \|b\|_{B^{0,\frac{1}{2}}(T)}
            \|u_F\|_{\widetilde{L}^\frac{2p}{p-1}_T(\widetilde{L}^{2p}_h(B^\frac{1}{2}_v))}
    \end{eqnarray*}
and
        \begin{eqnarray*}
        &&\int^T_0\left|
      \int_{\mathbb{R}^3}
           \Delta_j^v(a^h\otimes u_F):
            \nabla_h\Delta^v_jbdx
      \right|dt\\
                &\leq&\|           \Delta_j^v(a^h\otimes u_F)
      \|_{L^2_T(L^2(\mathbb{R}^3))}
      \|\Delta^v_j\nabla_h b\|_{L^2_T(L^2(\mathbb{R}^3))}\\
            &\lesssim&
            d_j^2\nu_h^{-\frac{1}{2}-\frac{1}{2p}}2^{-j}\|a\|_{B^{0,\frac{1}{2}}(T)}
            \|b\|_{B^{0,\frac{1}{2}}(T)}
             \|u_F\|_{\widetilde{L}^\frac{2p}{p-1}_T(\widetilde{L}^{2p}_h(B^\frac{1}{2}_v))}.
    \end{eqnarray*}

On the other hand, using Bony's decomposition in the vertical
variables, we obtain
    \begin{equation}
      \Delta^v_j(a^3\partial_3u_F)=
      \sum_{|j'-j|\leq 5}\Delta^v_j(S^v_{j'-1}a^3\partial_3\Delta^v_{j'}u_F)
      +
            \sum_{j'\geq j- N_0}
            \Delta^v_j(\Delta^v_{j' }a^3\partial_3S^v_{j'+2}u_F).
    \end{equation}
Using H\"{o}lder's inequality, Corollary \ref{ANS-C2.3}, we get
    \begin{eqnarray*}
      \|\Delta^v_j(S^v_{j'-1}a^3\partial_3\Delta^v_{j'}u_F)\|_{L^1_T(L^2(\mathbb{R}^3))}
      &\lesssim&
      2^{j'}\|S^v_{j'-1}a^3\|_{L^\infty_T(L^2_h(L^\infty_v))}
      \|\Delta^v_{j'}u_F\|_{L^1_T(L^\infty_h(L^2_v))}\\
        &\lesssim&d_{j'}2^{-\frac{j'}{2}}\|u_F\|_{\widetilde{L}^1_T(
        \widetilde{L}^\infty_h(B^\frac{3}{2}_v))}
        \|a\|_{B^{0,\frac{1}{2}}(T)}
    \end{eqnarray*}
and
    $$
    \sum_{|j'-j|\leq5}\|\Delta^v_j(S^v_{j'-1}a^3\partial_3\Delta^v_{j'}u_F)\|_{
      L^1_T(L^2(\mathbb{R}^3))}\lesssim
        d_{j }2^{-\frac{j}{2}}\|u_F\|_{\widetilde{L}^1_T(
        \widetilde{L}^\infty_h(B^\frac{3}{2}_v))}
        \|a\|_{B^{0,\frac{1}{2}}(T)}.
    $$

    Using H\"{o}lder's inequality, Lemmas \ref{ANS-L2.1}, \ref{ANS-L2.5} and the fact
    that div$a=0$, we have
        \begin{eqnarray*}
           \|\Delta^v_{j'}a^3S^v_{j'+2}\partial_3u_F\|_{L^\frac{2p}{2p-1}_T(L^\frac{2p}{p+1}_h(L^2_v))}
           &\lesssim&
           2^{j'}\|\Delta^v_{j'}a^3\|_{L^2_T(L^2(\mathbb{R}^3))}
           \|S^v_{j'+2} u_F\|_{L^\frac{2p}{p-1}_T(L^{2p}_h(L^\infty_v))}\\
                &\lesssim&\|\Delta^v_{j'}\partial_3a^3\|_{L^2_T(L^2(\mathbb{R}^3))}
           \|S^v_{j'+2} u_F\|_{L^\frac{2p}{p-1}_T(L^{2p}_h(L^\infty_v))}\\
                    &=&\|\Delta^v_{j'}\mathrm{div}_ha^h\|_{L^2_T(L^2(\mathbb{R}^3))}
           \|S^v_{j'+2} u_F\|_{L^\frac{2p}{p-1}_T(L^{2p}_h(L^\infty_v))}\\
           &\lesssim&d_{j'}\nu_h^{-\frac{1}{2}}2^{-\frac{j'}{2}}\|u_F\|_{\widetilde{L}^\frac{2p}{p-1}_T(
           \widetilde{L}^{2p}_h(B^\frac{1}{2}_v))}
        \|a\|_{B^{0,\frac{1}{2}}(T)}
    \end{eqnarray*}
and
    $$
    \sum_{j'\geq j-N_0}\|\Delta^v_{j'}a^3S^v_{j'+2}\partial_3u_F\|_{L^\frac{2p}{2p-1}_T
(L^\frac{2p}{p+1}_h(L^2_v))}\lesssim
        d_{j }\nu_h^{-\frac{1}{2}}2^{-\frac{j }{2}}\|u_F\|_{\widetilde{L}^\frac{2p}{p-1}_T
        (\widetilde{L}^{2p}_h(B^\frac{1}{2}_v))}
        \|a\|_{B^{0,\frac{1}{2}}(T)}.
    $$
This ends the proof of Proposition \ref{ANS-P4.2}.   {\hfill
$\square$\medskip}

\section{The proof of the uniqueness}\label{ANS-Sec4}
The first step to prove the uniqueness part of Theorem
\ref{ANS-T1} is the proof of the following regularity theorem.
\begin{thm}\label{ANS-T4.1}
  Let $u\in B^{-1+\frac{2}{p},\frac{1}{2}}_p(T)$ be a solution of
  (\ref{ANS-E1.1}) with initial data $u_0\in
  B^{-1+\frac{2}{p},\frac{1}{2}}_p$, $[u_0]_{E^p_T}<\infty$, $p\geq2$. Then, there exists a
  $T_1\in(0,T]$ such that
    $$w=u-u_F\in B^{0,\frac{1}{2}}(T_1).$$
\end{thm}
\noindent\textbf{Proof.}
  We already observe at the beginning of Section \ref{ANS-Sec3}
  that the vector field $w$ is the solution of the linear problem,
  which is
    \begin{equation}\label{ANS-E4.2}
      \left\{\begin{array}{l}
            w_t-\nu_h\Delta_hw-\nu_3\partial_3^2w=-\nabla P-u\cdot\nabla
            w-w\cdot\nabla u_F-u_F\cdot\nabla u_F,\\
            \mathrm{div}w=0,            \\
            w|_{t=0}=u_{0ll}.
      \end{array}
      \right.
    \end{equation}
Let us apply the operator $\Delta_j^v$ to the system
(\ref{ANS-E4.2}), and set
    $
    w_j=\Delta^v_jw.
    $
By the $L^2$ energy estimate, we have
    \begin{eqnarray*}
      &&\|w_j(t)\|^2_{L^2}+2\nu_h\int^t_0\|\nabla_h
      w_j(s)\|^2_{L^2}ds+2\nu_3\int^t_0\|\partial_3
      w_j(s)\|^2_{L^2}ds\\
        &\leq&\|\Delta^v_ju_{0ll}\|^2_{L^2}
        +2\int^t_0\left|\int_{\mathbb{R}^3}\Delta^v_j(u\cdot\nabla
        w)w_jdx
        \right|ds\\
                &&        +2\int^t_0\left|\int_{\mathbb{R}^3}\Delta^v_j(w\cdot\nabla
        u_F)w_jdx
        \right|ds
        +2\int^t_0\left|\int_{\mathbb{R}^3}\Delta^v_j(u_F\cdot\nabla
        u_F)w_jdx
        \right|ds.
    \end{eqnarray*}
From Propositions \ref{ANS-P4.1}-\ref{ANS-P4.2}, we obtain,
$t\in[0,T]$,
        \begin{eqnarray*}
      &&\|w_j(t)\|^2_{L^2}+\nu_h\int^t_0\|\nabla_h
      w_j(s)\|^2_{L^2}ds+\nu_3\int^t_0\|\partial_3
      w_j(s)\|^2_{L^2}ds\\
        &\lesssim&\|\Delta^v_ju_{0ll}\|^2_{L^2}
        +d_j^2\nu_h^{-\frac{1}{2}-\frac{1}{2p}}2^{-j}\|w\|^2_{B^{0,\frac{1}{2}}(t)}
        \left(\|u
        \|_{\widetilde{L}^\frac{2p}{p-1}_t(\widetilde{L}^{2p}_h(B^\frac{1}{2}_v))}+
            \|u_F\|_{\widetilde{L}^\frac{2p}{p-1}_t(\widetilde{L}^{2p}_h(B^\frac{1}{2}_v))}\right)\\
                &&+d_j^22^{-j}\|w\|^2_{B^{0,\frac{1}{2}}(t)}
            \|u_F\|_{\widetilde{L}^1_t(\widetilde{L}^\infty_h(B^\frac{3
}{2}_v))}+d_j^22^{-j}\|w\|_{B^{0,\frac{1}{2}}(t)}
                \|u_F\cdot\nabla
                u_F\|_{ {L}^1_t(B^{0,\frac{1}{2}})}
    \end{eqnarray*}
and
     \begin{eqnarray*}
    &&\|w\|_{B^{0,\frac{1}{2}}(t)}\\
            &\lesssim& \|u_{0ll}\|_{B^{0,\frac{1}{2}}}
        +\|u_F\cdot\nabla u_F\|_{ {L}^1_t(B^{0,\frac{1}{2}})} +(1+\nu_h^{-\frac{1}{4}-\frac{1}{4p}})\|w\|_{B^{0,\frac{1}{2}}(t)} \\
            &&   \times\left(\|u
        \|_{\widetilde{L}^\frac{2p}{p-1}_t(\widetilde{L}^{2p}_h(B^\frac{1}{2}_v))}+
            \|u_F\|_{\widetilde{L}^\frac{2p}{p-1}_t(\widetilde{L}^{2p}_h(B^\frac{1}{2}_v))}
+\|u_F\|_{\widetilde{L}^1_t(\widetilde{L}^\infty_h(B^\frac{3
}{2}_v))}\right)^\frac{1}{2}.
    \end{eqnarray*}
Thus, we can choose a small $T_1\in(0,T]$, such that
        $\|u
        \|_{\widetilde{L}^\frac{2p}{p-1}_{T_1}(\widetilde{L}^{2p}_h(B^\frac{1}{2}_v))}+
            \|u_F\|_{\widetilde{L}^\frac{2p}{p-1}_{T_1}(\widetilde{L}^{2p}_h(B^\frac{1}{2}_v))}
+\|u_F\|_{\widetilde{L}^1_{T_1}(\widetilde{L}^\infty_h(B^\frac{3
}{2}_v))}$
            is small enough and
                $$\|w\|_{B^{0,\frac{1}{2}}(T_1)}\\
            \lesssim \|u_{0ll}\|_{B^{0,\frac{1}{2}}}
            +\|u_F\cdot\nabla u_F\|_{ {L}^1_{T_1}(B^{0,\frac{1}{2}})}.$$
This concludes the proof of Theorem \ref{ANS-T4.1}. {\hfill
$\square$\medskip}

The above theorem implies that, if $u_i$, $i=1,2$, are  two
solutions of (\ref{ANS-E1.1}) in the space
$B^{-1+\frac{2}{p},\frac{1}{2}}_p(T)$ associated with the same
initial data, then there exists a $T_1\in(0,T]$ such that the
difference $\delta:=u_2-u_1$ belongs to $B^{0,\frac{1}{2}}(T_1)$.
Moreover, $\delta$ satisfies the following system:
    \begin{equation}
      \left\{
      \begin{array}{l}
            \delta_t-\nu_h\Delta_h\delta-\nu_3\partial_3^2\delta=L\delta-\nabla
            P,\\
                \mathrm{div}\delta=0,\\
                    \delta|_{t=0}=0,
      \end{array}
      \right.\label{ANS-E4.2-0}
    \end{equation}
where $L$ is the following linear operator
    $$
    L\delta:=-\delta\cdot\nabla u_1-u_2\cdot\nabla \delta.
    $$
In order to prove the uniqueness, we have to prove that
$\delta\equiv0$.

As in \cite{Chemin07}, we give the following definitions.
\begin{defn}
  Let $s\in \mathbb{R}$, and let us define the following
  semi-norm:
    $$
        \|a\|_{H^{0,s}}:=\left(
        \sum_{j\in\mathbb{Z}}2^{2js}\|\Delta^v_ja\|_{L^2}^2
        \right)^\frac{1}{2}.
    $$
\end{defn}

\begin{defn}
  We denote by $\mathcal{H}$ the space of distributions, which
  is the completion of $\mathcal{S}(\mathbb{R}^3)$ by the following
  norm:
    $$
    \|a\|^2_{\mathcal{H}}:=
    \sum_{j\in\mathbb{Z}}2^{-j}\|\Delta^{vi}_ja\|_{L^2(\mathbb{R}^3)}^2<\infty,
    $$
  where
    $$
    \Delta^{vi}_j=\Delta^v_j,\ \textrm{ if }\ j\geq0,
    \ \textrm{ and }\ \Delta^{vi}_j=0,\ \textrm{ if }\ j\leq-2,
    \ \textrm{ and }\ \Delta^{vi}_{-1}=S^v_0.
    $$
\end{defn}

\begin{defn}
Let us denote by $\mathcal{B}$ the following semi-norm:
    $$
    \|a\|^2_{\mathcal{B}}:=\sumetage{k\in\mathbb{Z}}{j\in\mathbb{N}}2^{j-k(2-\frac{4}{p})}
    \|\Delta^h_k\Delta^v_ja\|^2_{L^p_h(L^2_v)}.
    $$
\end{defn}
\begin{rem}
  It is obvious that
    \begin{equation}
      \|a\|^2_{L^\infty_T(H^{0,\frac{1}{2}})}
      +\nu_h\|\nabla_ha\|^2_{L^2_T(H^{0,\frac{1}{2}})}\lesssim\|a\|^2_{B^{0,\frac{1}{2}}(T)}\label{ANS-E4.5-1}
    \end{equation}
  and
    \begin{equation}
      \|a\|^2_{L^\infty_T(\mathcal{B})}
      +\nu_h\|\nabla_ha\|^2_{L^2_T(\mathcal{B})}\lesssim\|a\|^2_{B^{-1+\frac{2}{p},\frac{1}{2}}_p(T)}.
      \label{ANS-E4.6-1}
    \end{equation}
\end{rem} Let us state the following variation of Lemma 4.2 of
\cite{Chemin07}.

\begin{lem}\label{ANS-L4.4}
  A constant $C$ exists such that, for any $p'\in[2p,\infty)$, $p\geq2$, we
  have
    $$
        \|\Delta^v_j b\|_{L^{p'}_h(L^2_v)}
        \leq Cc_j\sqrt{p'}2^{-\frac{j}{2}}\|b\|_{\mathcal{B}}^\frac{2}{p'}
        \|\nabla_hb\|_{\mathcal{B}}^{1-\frac{2}{p'}}, \ j\geq0.
    $$
\end{lem}
\noindent\textbf{Proof.}
  From Lemma \ref{ANS-L2.1}, we get
    \begin{eqnarray*}
      2^\frac{j}{2}\|\Delta^v_j b\|_{L^{p'}_h(L^2_v)}
      &\leq&C\sum_{k\leq N}2^{\frac{j}{2}+k(-1+\frac{2}{p})}2^{k(1-\frac{2}{p'})}
      \|\Delta^h_k\Delta^v_jb\|_{L^p_h(L^2_v)}\\
        &&+C\sum_{k\geq N}2^{\frac{j}{2}+k(-1+\frac{2}{p})}2^{-\frac{2}{p'}k}
      \|\Delta^h_k\Delta^v_j\nabla_hb\|_{L^p_h(L^2_v)}\\
                &\leq&C\|b\|_{\mathcal{B}}\sum_{k\leq
                N}2^{k(1-\frac{2}{p'})}c_{k,j}+C\|\nabla_hb\|_{\mathcal{B}}\sum_{k\geq
                N}2^{-\frac{2}{p'}k}c'_{k,j}.
    \end{eqnarray*}
  Using the Cauchy-Schwarz inequality, we get
        \begin{eqnarray*}
      2^\frac{j}{2}\|\Delta^v_j b\|_{L^{p'}_h(L^2_v)}
                &\leq&Cc_j\|b\|_{\mathcal{B}}\left(\sum_{k\leq
                N}2^{2k(1-\frac{2}{p'})}\right)^\frac{1}{2}+Cc_j\|\nabla_hb\|_{\mathcal{B}}
                \left(\sum_{k\geq
                N}2^{-\frac{4}{p'}k}\right)^\frac{1}{2}\\
                &\leq&Cc_j\|b\|_{\mathcal{B}}2^{N(1-\frac{2}{p'})}+C\sqrt{p'}c_j\|\nabla_hb\|_{\mathcal{B}}
                2^{-\frac{2}{p'}N}.
    \end{eqnarray*}
  Choosing $2^N\simeq \frac{\|\nabla_hb\|_{\mathcal{B}}}{\|b\|_{\mathcal{B}}}$ gives the
  lemma.
{\hfill $\square$\medskip}

Let us state the following variation of Lemma 4.1 of
\cite{Chemin07}.
\begin{lem}\label{ANS-L4.5}
  Let $a$ and $b$ be two divergence free vector fields such that
  $a$,  $\nabla_ha\in H^{0,\frac{1}{2}}\cap\mathcal{H}$, $b\in
 \mathcal{B}\cap \widetilde{L}^{2p}_h(B^\frac{1}{2}_v)$ and $\nabla_hb\in
  \mathcal{B}$. Let us assume in addition that
    $\|a\|_{\mathcal{H}}^2\leq 2^{-2^{2p}}$. Then, we have
        $$
        |(b\cdot\nabla a|a)_{\mathcal{H}}|
        + |(a\cdot\nabla b|a)_{\mathcal{H}}|
        \leq\frac{\nu_h}{10}\|\nabla_ha\|_{\mathcal{H}}^2+C(a,b)\mu(\|a\|_{\mathcal{H}}^2),
        $$
  where $ (f|g
  )_{\mathcal{H}}:=\sum_{j\in\mathbb{Z}}2^{-j}\int_{\mathbb{R}^3}\Delta^{vi}_jf\Delta^{vi}_jgdx$,
  $\mu(r):=r(1-\log_2r)\log_2(1-\log_2r)$ and
        \begin{eqnarray*}
          C(a,b)&:=&C_{\nu_h}\big(\|b\|^{\frac{2p}{p-1}}_{L^{2p}_h(L^\infty_v)}
                    +\|b\|^{\frac{2p}{2p-1}}_{L^{2p}_h(L^\infty_v)}+\|b\|^2_{\widetilde{L}^{2p}_h(B^\frac{1}{2}_v)}
            \|a\|^{2-\frac{2}{p}}_{H^{0,\frac{1}{2}}}\|\nabla_ha\|^{\frac{2}{p}}_{H^{0,\frac{1}{2}}}\\
                    &&+(1
             +\|b\|^2_{\mathcal{B}})
            \|\nabla_hb\|^{2}_{\mathcal{B}}
            \big).
        \end{eqnarray*}
\end{lem}
\noindent\textbf{Proof.}
  \textbf{The estimate of the term $(b\cdot\nabla a
  |a)_{\mathcal{H}}$.} Using  Bony's paradifferential decomposition
  in the vertical variable and in the inhomogeneous context, we
  have
        $$
        b\cdot \nabla a =T_b\nabla a+\widetilde{R}(b,\nabla a)
        $$
  with
    $$
    T_b\nabla a:=\sum_lS^{vi}_{l-1}b\cdot \nabla \Delta^{vi}_la,
    \ \widetilde{R}(b,\nabla a):=\sum_{l}\Delta^{vi}_lb\cdot\nabla
    S^{vi}_{l+2}a\ \textrm{ and }
    \ S^{vi}_l =\sum_{l'\leq l-1}\Delta^{vi}_{l'}.
    $$

    \noindent\textit{Step 1.} The estimate of $(T_b\nabla
    a|a)_{\mathcal{H}}$. As usual, we shall treat terms
    involving vertical derivatives in a different way from terms
    involving horizontal derivatives. This leads to
            $$
        \Delta^{vi}_j(T_b\nabla a)=T^h_j+T^v_j,
            $$
    with
        $$
        T^h_j:=\Delta^{vi}_j\sum_{|j-l|\leq5}S^{vi}_{l-1}b^h\cdot\nabla_h\Delta^{vi}_la,
        \ \textrm{ and }\
        T^v_j:=\Delta^{vi}_j\sum_{|j-l|\leq5}S^{vi}_{l-1}b^3\partial_3\Delta^{vi}_la.
        $$
Using H\"{o}lder's inequality, we obtain
        $$
          \|T^h_j\|_{L^\frac{2p}{p+1}_h(L^2_v)}\lesssim \|b\|_{L^{2p}_h(L^\infty_v)}
          \sum_{|j-l|\leq5}\|\nabla_h\Delta^{vi}_la\|_{L^2}
                    \lesssim c_j2^{\frac{j}{2}}\|b\|_{L^{2p}_h(L^\infty_v)}
                    \|\nabla_ha\|_{\mathcal{H}}
        $$
and
    $$
    |(T^h_j|\Delta^{vi}_ja)_{L^2}|
    \lesssim c_j2^\frac{j}{2}\|b\|_{L^{2p}_h(L^\infty_v)}
    \|\Delta^{vi}_ja\|_{L^\frac{2p}{p-1}_h(L^2_v)}
                    \|\nabla_ha\|_{\mathcal{H}}.
    $$

Using Minkowski's inequality and the Gagliardo-Nirenberg
inequality, we have
      $$
      \|f\|_{L^4_h(L^2_v)}\leq C\|f\|_{L^2_v(L^4_h)}
        \leq C\left\|\|f(\cdot,x_3)\|_{L^2_h}^{\frac{1}{2}}
      \|\nabla_hf(\cdot,x_3)\|_{L^2_h}^{\frac{1}{2}}\right\|_{L^2_v}
      \leq C\|f\|_{L^2}^{\frac{1}{2}}
      \|\nabla_hf\|_{L^2}^{\frac{1}{2}}.
    $$
By interpolation, we have
    \begin{equation}
      \|f\|_{L^q_h(L^2_v)}\leq C
      \|f\|_{L^2_h(L^2_v)}^{\frac{2}{q}}
      \|\nabla_hf\|_{L^2_h(L^2_v)}^{1-\frac{2}{q}},
      \ q\in[2,4]. \label{ANS-E4.4}
    \end{equation}
 Then, we get
     $$
      \|\Delta^{vi}_ja\|_{L^\frac{2p}{p-1}_h(L^2_v)}\lesssim
      \|\Delta^{vi}_ja\|_{L^2_h(L^2_v)}^{1-\frac{1}{p}}
      \|\nabla_h\Delta^{vi}_ja\|_{L^2_h(L^2_v)}^{\frac{1}{p}}
    $$
and
    \begin{equation}
      \sum_j2^{-j}|(T^h_j|\Delta^{vi}_ja)_{L^2}|
    \lesssim \|b\|_{L^{2p}_h(L^\infty_v)}
    \|a\|_{\mathcal{H}}^{1-\frac{1}{p}}
                    \|\nabla_ha\|_{\mathcal{H}}^{1+\frac{1}{p}}.
                    \label{ANS-E4.5}
    \end{equation}
The estimate of $(T^v_j|\Delta^{vi}_ja)_{L^2}$ is more delicate.
Let us write that
    $$
    T^v_j=\sum^3_{n=1}T^{v,n}_j
    $$
with
    $$
    T^{v,1}_j:=S^{vi}_{j-1}b^3\partial_3\Delta^{vi}_ja,
     \   T^{v,2}_j:=\sum_{|j-l|\leq5}[\Delta^{vi}_j;S^{vi}_{l-1}b^3]\partial_3\Delta^{vi}_la,
        $$
            $$
        T^{v,3}_j:=\sum_{|j-l|\leq5}(S^{vi}_{l-1}b^3-S^{vi}_{j-1}b^3)\partial_3\Delta^{vi}_j\Delta^{vi}_la.
            $$

\noindent\textit{Step 1a.} The estimate of
$\sum_j2^{-j}|(T^{v,1}_j|\Delta^{vi}_ja)_{L^2}|$. To do this, we
use the tricks from \cite{Chemin07,Chemin00} once again. Using
integration by parts and div$b=0$, we get
    \begin{eqnarray*}
    (T^{v,1}_j|\Delta^{vi}_ja)_{L^2}&=&-\frac{1}{2}
    \int_{\mathbb{R}^3}S^{vi}_{j-1}\partial_3b^3\Delta^{vi}_ja\Delta^{vi}_jadx\\
            &=&\frac{1}{2}
    \int_{\mathbb{R}^3}S^{vi}_{j-1}\mathrm{div}_hb^h\Delta^{vi}_ja\Delta^{vi}_jadx\\
            &=&-
    \int_{\mathbb{R}^3}S^{vi}_{j-1}b^h\cdot\nabla_h\Delta^{vi}_ja\Delta^{vi}_jadx.
    \end{eqnarray*}
Similar to the proof of (\ref{ANS-E4.5}), we have
    \begin{equation}
      \sum_j2^{-j}|(T^{v,1}_j|\Delta^{vi}_ja)_{L^2}|
      \lesssim \|b\|_{L^{2p}_h(L^\infty_v)}
    \|a\|_{\mathcal{H}}^{1-\frac{1}{p}}
                    \|\nabla_ha\|_{\mathcal{H}}^{1+\frac{1}{p}}.
    \end{equation}

\noindent\textit{Step 1b.} The estimate of
$\sum_j2^{-j}|(T^{v,2}_j|\Delta^{vi}_ja)_{L^2}|$. In order to
estimate the commutator, let us use the Taylor formula (as in
\cite{Chemin07}). For a function $f$ on $\mathbb{R}^3$, we define
the function $\tilde{f}$ on $\mathbb{R}^4$ by
    $$
    \tilde{f}(x,y_3):=\int^1_0f(x_h,x_3+\tau(y_3-x_3))d\tau.
    $$
Then, denoting $\bar{h}(x_3):=x_3h(x_3)$, we have
    $$
    T^{v,2}_j=\sum_{|j-l|\leq5}
    \int_{\mathbb{R}}\bar{h}(2^j(x_3-y_3))(\widetilde{S^{vi}_{l-1}\partial_3b^3})(x,y_3)
    \partial_3\Delta^{vi}_la(x_h,y_3)dy_3.
    $$
Using div$b=0$ and $\partial_h\tilde{f}=\widetilde{\partial_hf}$,
we obtain
    $$
    T^{v,2}_j=-\sum_{|j-l|\leq5}
    \int_{\mathbb{R}}\bar{h}(2^j(x_3-y_3))\mathrm{div}_h(\widetilde{S^{vi}_{l-1}b^h})(x,y_3)
    \partial_3\Delta^{vi}_la(x_h,y_3)dy_3.
    $$
Using integration by parts with respect to the horizontal
variable, we have
    \begin{eqnarray*}
            &&(T^{v,2}_j|\Delta^{vi}_ja)_{L^2}\\
    &=&\sum_{|j-l|\leq5}\left\{
    \int_{\mathbb{R}^4}\bar{h}(2^j(x_3-y_3))(\widetilde{S^{vi}_{l-1}b^h})(x,y_3)
    \cdot\nabla_h\partial_3\Delta^{vi}_la(x_h,y_3)\Delta^{vi}_ja(x)dxdy_3
    \right.\\
           && +\left.
     \int_{\mathbb{R}^4}\bar{h}(2^j(x_3-y_3))(\widetilde{S^{vi}_{l-1}b^h})(x,y_3)
    \cdot\nabla_h\Delta^{vi}_ja(x)\partial_3\Delta^{vi}_la(x_h,y_3)dxdy_3
            \right\}.
    \end{eqnarray*}
From $\|\tilde{f}(x_h,\cdot,\cdot)\|_{L^\infty}\leq
\|f(x_h,\cdot)\|_{L^\infty_v}$, Young's inequality and
H\"{o}lder's inequality, we obtain
        \begin{eqnarray*}
       && |(T^{v,2}_j|\Delta^{vi}_ja)_{L^2}|\\
        &\lesssim&
        \|b\|_{L^{2p}_h(L^\infty_v)}\sum_{|j-l|\leq5}\left(
    \|\Delta^{vi}_l\nabla_ha\|_{L^2}
                    \|\Delta^{vi}_ja\|_{L^\frac{2p}{p-1}_h(L^2_v)}
                    +\|\Delta^{vi}_la\|_{L^\frac{2p}{p-1}_h(L^2_v)}
                    \|\nabla_h\Delta^{vi}_ja\|_{L^2}
                    \right).
        \end{eqnarray*}
Similar to the proof of (\ref{ANS-E4.5}), we have
    \begin{equation}
      \sum_j2^{-j}|(T^{v,2}_j|\Delta^{vi}_ja)_{L^2}|
      \lesssim \|b\|_{L^{2p}_h(L^\infty_v)}
    \|a\|_{\mathcal{H}}^{1-\frac{1}{p}}
                    \|\nabla_ha\|_{\mathcal{H}}^{1+\frac{1}{p}}.
                    \label{ANS-E4.8}
    \end{equation}

\noindent\textit{Step 1c.} The estimate of
$\sum_j2^{-j}|(T^{v,3}_j|\Delta^{vi}_ja)_{L^2}|$. For any
divergence free vector field $u$, from (\ref{ANS-E4.1}), we have,
$l\geq0$,
    \begin{eqnarray}
        \Delta^v_lu^3(x)&=&\int_{\mathbb{R}}g^v(2^l(x_3-y_3))\Delta^v_l
        \partial_3u^3(x_h,y_3)dy_3\nonumber\\
                &=&-\mathrm{div}_h\int_{\mathbb{R}}g^v(2^l(x_3-y_3))\Delta^v_lu^h(x_h,y_3)dy_3
                \nonumber\\
        &=&-2^{-l}\mathrm{div}_h\widetilde{\Delta}^v_lu^h.\label{ANS-E4.9}
    \end{eqnarray}
If $j\geq7$, then the terms $S^{vi}_{l-1}b^3-S^{vi}_{j-1}b^3$ that
appear in $T^{v,3}_j$ are a sum of the terms $\Delta^{vi}_{l'}$
with $l'\geq0$. If $j\geq7$, using (\ref{ANS-E4.9}) and
integration by parts in the horizontal variable, we obtain
    \begin{eqnarray*}
      (T^{v,3}_j|\Delta^{vi}_ja)_{L^2}
      &=&\sumetage{|l'-l|\leq5}{|l-j|\leq5}
      2^{-l'}\left(\int_{\mathbb{R}^3}\Delta^{vi}_j(\widetilde{\Delta}^{vi}_{l'}b^h\cdot\nabla_h\Delta^{vi}_l\partial_3a)
      \Delta^{vi}_ja
      dx\right.\\
           && +\left.\int_{\mathbb{R}^3}\Delta^{vi}_j(\Delta^{vi}_l\partial_3a\widetilde{\Delta}^{vi}_{l'}b^h)
            \cdot\nabla_h\Delta^{vi}_ja
      dx      \right).
    \end{eqnarray*}
Similar to the proof of (\ref{ANS-E4.8}), we have
    \begin{equation}
      \sum_{j\geq7}2^{-j}|(T^{v,3}_j|\Delta^{vi}_ja)_{L^2}|
      \lesssim \|b\|_{L^{2p}_h(L^\infty_v)}
    \|a\|_{\mathcal{H}}^{1-\frac{1}{p}}
                    \|\nabla_ha\|_{\mathcal{H}}^{1+\frac{1}{p}}.\label{ANS-E4.10}
    \end{equation}
If $j\leq 7$, we can easily get
    $$
      |(T^{v,3}_j|\Delta^{vi}_ja)_{L^2}|
      \lesssim \|b\|_{L^{2p}_h(L^\infty_v)}
    \|a\|_{\mathcal{H}}^{2-\frac{1}{p}}
                    \|\nabla_ha\|_{\mathcal{H}}^{\frac{1}{p}}.
    $$
Plugging this inequality with inequalities
(\ref{ANS-E4.5})-(\ref{ANS-E4.8}) and (\ref{ANS-E4.10}), using
Young's inequality, we have
    \begin{eqnarray}
    |(T_b\nabla a|a)_{\mathcal{H}}|&\lesssim&
    \|b\|_{L^{2p}_h(L^\infty_v)}\left(
    \|a\|_{\mathcal{H}}^{1-\frac{1}{p}}
                    \|\nabla_ha\|_{\mathcal{H}}^{1+\frac{1}{p}}+
                    \|a\|_{\mathcal{H}}^{2-\frac{1}{p}}
                    \|\nabla_ha\|_{\mathcal{H}}^{\frac{1}{p}}
                    \right)\nonumber\\
                        &\leq&
                    \frac{\nu_h}{100}\|\nabla_h
                    a\|_{\mathcal{H}}^2
                    +C_{\nu_h}
                    \|a\|_{\mathcal{H}}^2\left(\|b\|^{\frac{2p}{p-1}}_{L^{2p}_h(L^\infty_v)}
                    +\|b\|^{\frac{2p}{2p-1}}_{L^{2p}_h(L^\infty_v)}\right).
    \end{eqnarray}

\noindent\textit{Step 2.} The estimate of $(\widetilde{R}(b,\nabla
a)|a)_{\mathcal{H}}$. Again, let us treat terms involving vertical
derivatives in a different way from terms involving horizontal
derivatives. This leads to
        $$
    \Delta^{vi}_j\widetilde{R}(b,\nabla a)=
    R^h_j+R^v_j+R^0_j
        $$
    with
        $$
    R^h_j:=\Delta^{vi}_j\sum_{l\geq(j-N_0)^+}\Delta^v_lb^h\cdot\nabla_hS^v_{l+2}a,
        $$
        $$
    R^v_j:=\Delta^{vi}_j\sum_{l\geq(j-N_0)^+}\Delta^v_lb^3\partial_3S^v_{l+2}a,
        $$
        $$
    R^0_j:=\Delta^{vi}_j(S^v_0b\cdot\nabla S^v_{2}a).
        $$
    Let us first estimate $R^0_j$. It is obvious that if $j$ is
    large enough, this term is 0. Thus, if $j\leq N_1$, we obtain
   \begin{eqnarray}
   \sum_{j\leq N_1} |(R^0_j|\Delta^{vi}_ja)_{L^2}|&\lesssim&
    \|b\|_{L^{2p}_h(L^\infty_v)}\left(
    \|a\|_{\mathcal{H}}^{1-\frac{1}{p}}
                    \|\nabla_ha\|_{\mathcal{H}}^{1+\frac{1}{p}}+
                    \|a\|_{\mathcal{H}}^{2-\frac{1}{p}}
                    \|\nabla_ha\|_{\mathcal{H}}^{\frac{1}{p}}
                    \right)\nonumber\\
                        &\leq&
                    \frac{\nu_h}{100}\|\nabla_h
                    a\|_{\mathcal{H}}^2
                    +C_{\nu_h}
                    \|a\|_{\mathcal{H}}^2\left(\|b\|^{\frac{2p}{p-1}}_{L^{2p}_h(L^\infty_v)}
                    +\|b\|^{\frac{2p}{2p-1}}_{L^{2p}_h(L^\infty_v)}\right).
    \end{eqnarray}

\noindent\textit{Step 2a.} The estimate of
$\sum_j2^{-j}|(R^h_j|\Delta^{vi}_ja)_{L^2}|$. First, we estimate
$R^h_j$ in high (vertical) frequencies. From Lemma \ref{ANS-L2.1}
and H\"{o}lder's inequality, we have
    \begin{eqnarray*}
      \|R^h_j\|_{L^\frac{2p}{p+1}_h(L^2_v)}&\lesssim&
      2^{\frac{j}{2}}\sum_{l\geq(j-N_0)^+}\|\Delta^v_lb^h\cdot\nabla_hS^v_{l+2}a\|_{
      L^\frac{2p}{p+1}_h(L^1_v)}\\
            &\lesssim&2^{\frac{j}{2}}\sum_{l\geq(j-N_0)^+}\|\Delta^v_lb^h\|_{L^{2p}_h(L^2_v)}
            \|\nabla_hS^v_{l+2}a\|_{
      L^2_h(L^2_v)}\\
            &\lesssim&2^{\frac{j}{2}}\sum_{l\geq(j-N_0)^+}\|\Delta^v_lb^h\|_{L^{2p}_h(L^2_v)}
            \sum^{l+1}_{l'=-1}\|\nabla_h\Delta^{vi}_{l'}a\|_{
      L^2_h(L^2_v)}\\
            &\lesssim&\|b^h\|_{\widetilde{L}^{2p}_h(B^\frac{1}{2}_v)}
            \|\nabla_ha\|_{\mathcal{H}}\sum_{l\geq(j-N_0)^+}d_l2^{\frac{1}{2}(j-l)}
            \left(
            \sum^{l+1}_{l'=-1}2^{l'}\right)^\frac{1}{2}\\
            &\lesssim&2^\frac{j}{2}\|b^h\|_{\widetilde{L}^{2p}_h(B^\frac{1}{2}_v)}
            \|\nabla_ha\|_{\mathcal{H}}.
    \end{eqnarray*}
Combining it with (\ref{ANS-E4.4}), we have
    \begin{eqnarray}
      |(R^h_j|\Delta^{vi}_ja)_{L^2}|&\lesssim&
      2^\frac{j}{2}\|b^h\|_{\widetilde{L}^{2p}_h(B^\frac{1}{2}_v)}
            \|\nabla_ha\|_{\mathcal{H}}\|\Delta^{vi}_ja\|_{L^\frac{2p}{p-1}_h(L^2_v)}\nonumber\\
                &\lesssim&
      \|b \|_{\widetilde{L}^{2p}_h(B^\frac{1}{2}_v)}
            \|\nabla_ha\|_{\mathcal{H}}\|a\|^{1-\frac{1}{p}}_{H^{0,\frac{1}{2}}}
            \|\nabla_ha\|^{\frac{1}{p}}_{H^{0,\frac{1}{2}}}.\label{ANS-E4.13}
    \end{eqnarray}

     Then, we estimate
$R^h_j$ in low (vertical) frequencies. Following the idea of
\cite{Chemin07,Paicu05}, using Lemmas \ref{ANS-L2.1} and
\ref{ANS-L4.4}, we obtain
    \begin{eqnarray*}
      \|R^h_j\|_{L^\frac{2p'}{p'+2}_h(L^2_v)}
      &\lesssim&2^{\frac{j}{2}}\sum_{l\geq(j-N_0)^+}
      \|\Delta^v_lb^h\cdot\nabla_hS^v_{l+2}a\|_{L^\frac{2p'}{p'+2}_h(L^1_v)}\\
        &\lesssim&2^{\frac{j}{2}}\sum_{l\geq(j-N_0)^+}
      \|\Delta^v_lb^h\|_{L^{p'}_h(L^2_v)}\|
      \nabla_hS^v_{l+2}a\|_{L^2_h(L^2_v)}\\
            &\lesssim&2^{\frac{j}{2}}\sum_{l\geq(j-N_0)^+}c_l2^{-\frac{l}{2}}\sqrt{p'}
      \|b\|^\frac{2}{p'}_{\mathcal{B}}\|\nabla_hb\|^{1-\frac{2}{p'}}_{\mathcal{B}}
      \sum_{k\leq l+1}c_k2^{\frac{k}{2}}\|
      \nabla_h a\|_{\mathcal{H}}\\
            &\lesssim&2^{\frac{j}{2}}\sqrt{p'}
      \|b\|^\frac{2}{p'}_{\mathcal{B}}\|\nabla_hb\|^{1-\frac{2}{p'}}_{\mathcal{B}}\|
      \nabla_h a\|_{\mathcal{H}}.
    \end{eqnarray*}
By (\ref{ANS-E4.4}), a constant $C$ exists (independent of $p'$)
such that,
    $$
    \|\Delta^{vi}_ja\|_{L^\frac{2p'}{p'-2}_h(L^2_v)}
    \leq C \|\Delta^{vi}_ja\|^{1-\frac{2}{p'}}_{L^2}
            \|\nabla_h\Delta^{vi}_ja\|^{\frac{2}{p'}}_{L^2}.
    $$
Thus we get
    \begin{equation}
      |(R^h_j|\Delta^{vi}_ja)_{L^2}| \lesssim
        c_j2^j\sqrt{p'}\|b\|^\frac{2}{p'}_{\mathcal{B}}\|\nabla_hb\|^{1-\frac{2}{p'}}_{\mathcal{B}}
     \|a\|_{\mathcal{H}}^{1-\frac{2}{p'}}
    \|\nabla_ha\|_{\mathcal{H}}^{1+\frac{2}{p'}}.\label{ANS-E4.14}
    \end{equation}
From (\ref{ANS-E4.13})-(\ref{ANS-E4.14}), the Cauchy-Schwarz
inequality and Young's inequality, we have
    \begin{eqnarray*}
      \sum_j2^{-j}|(R^h_j|\Delta^{vi}_ja)_{L^2}|
      &\lesssim&\sum_{j>M}2^{-j}\|b\|_{\widetilde{L}^{2p}_h(B^\frac{1}{2}_v)}
            \|\nabla_ha\|_{\mathcal{H}}\|a\|^{1-\frac{1}{p}}_{H^{0,\frac{1}{2}}}
            \|\nabla_ha\|^{\frac{1}{p}}_{H^{0,\frac{1}{2}}}\\
            &&+\sum_{-1\leq j\leq M} c_j\sqrt{p'}\|b\|^\frac{2}{p'}_{\mathcal{B}}
            \|\nabla_hb\|^{1-\frac{2}{p'}}_{\mathcal{B}}
     \|a\|_{\mathcal{H}}^{1-\frac{2}{p'}}
    \|\nabla_ha\|_{\mathcal{H}}^{1+\frac{2}{p'}}\\
        &\lesssim&2^{-M}\|b\|_{\widetilde{L}^{2p}_h(B^\frac{1}{2}_v)}
            \|\nabla_ha\|_{\mathcal{H}}\|a\|^{1-\frac{1}{p}}_{H^{0,\frac{1}{2}}}
            \|\nabla_ha\|^{\frac{1}{p}}_{H^{0,\frac{1}{2}}}\\
            &&+\sqrt{p'M}\|b\|^\frac{2}{p'}_{\mathcal{B}}
            \|\nabla_hb\|^{1-\frac{2}{p'}}_{\mathcal{B}}
     \|a\|_{\mathcal{H}}^{1-\frac{2}{p'}}
    \|\nabla_ha\|_{\mathcal{H}}^{1+\frac{2}{p'}}\\
        &\leq&
        \frac{\nu_h}{100}\|\nabla_ha\|_{\mathcal{H}}^2
        +\frac{C}{\nu_h}2^{-2M}\|b\|^2_{\widetilde{L}^{2p}_h(B^\frac{1}{2}_v)}
            \|a\|^{2-\frac{2}{p}}_{H^{0,\frac{1}{2}}}
            \|\nabla_ha\|^{\frac{2}{p}}_{H^{0,\frac{1}{2}}}\\
            &&+C\nu_h^{-\frac{p'+2}{p'-2}}(p'M)^\frac{p'}{p'-2}\|b\|^\frac{4}{p'-2}_{\mathcal{B}}
            \|\nabla_hb\|^{2}_{\mathcal{B}}
     \|a\|_{\mathcal{H}}^{2}.
    \end{eqnarray*}
Let us assume that $M\geq 2^{2p}$. Choosing $p'=\log_2M$, we get
        \begin{eqnarray*}
      \sum_j2^{-j}|(R^h_j|\Delta^{vi}_ja)_{L^2}|
        &\leq&
        \frac{\nu_h}{100}\|\nabla_ha\|_{\mathcal{H}}^2
        +\frac{C}{\nu_h}2^{-2M}\|b \|^2_{\widetilde{L}^{2p}_h(B^\frac{1}{2}_v)}
            \|a\|^{2-\frac{2}{p}}_{H^{0,\frac{1}{2}}}
            \|\nabla_ha\|^{\frac{2}{p}}_{H^{0,\frac{1}{2}}}\\
            &&+C_{\nu_h}(1+\|b\|^{2}_{\mathcal{B}})
            \|\nabla_hb\|^{2}_{\mathcal{B}}
     \|a\|_{\mathcal{H}}^{2}M\log_2M.
    \end{eqnarray*}
If $\|a\|_{\mathcal{H}}\leq 2^{-2^{2p}}$, then we can choose $M$
such that $2^{-M}\simeq \|a\|_{\mathcal{H}}$, and get
    \begin{equation}
        \sum_j2^{-j}|(R^h_j|\Delta^{vi}_ja)_{L^2}|
        \leq \frac{\nu_h}{100}\|\nabla_ha\|_{\mathcal{H}}^2
        +C_1(a,b)\mu(\|a\|_{\mathcal{H}}),\label{ANS-E4.15}
    \end{equation}
with
    $$
    C_1(a,b)=\frac{C}{\nu_h}\|b\|^2_{\widetilde{L}^{2p}_h(B^\frac{1}{2}_v)}
            \|a\|^{2-\frac{2}{p}}_{H^{0,\frac{1}{2}}}
            \|\nabla_ha\|^{\frac{2}{p}}_{H^{0,\frac{1}{2}}}
            +C_{\nu_h}(1+\|b\|^{2}_{\mathcal{B}})
            \|\nabla_hb\|^{2}_{\mathcal{B}}.
    $$

\noindent\textit{Step 2b.} The estimate of
$\sum_j2^{-j}|(R^v_j|\Delta^{vi}_ja)_{L^2}|$. From
(\ref{ANS-E4.9}) and integration by parts in the horizontal
variable, we have
    $$
    (R^v_j|\Delta^{vi}_ja)_{L^2}=R^{v,1}_j+R^{v,2}_j
    $$
with
    $$
    R^{v,1}_j:=\sum_{l\geq(j-N_0)^+}
    2^{-l}\int_{\mathbb{R}^3}\Delta^{vi}_j(\widetilde{\Delta}^v_lb^h\cdot\nabla_h\partial_3S^v_{l+2}a)
    \Delta^{vi}_jadx,
    $$
    $$
    R^{v,2}_j:=\sum_{l\geq(j-N_0)^+}
    2^{-l}\int_{\mathbb{R}^3}\Delta^{vi}_j(\widetilde{\Delta}^v_lb^h\partial_3S^v_{l+2}a)\cdot\nabla_h
    \Delta^{vi}_jadx.
    $$
Since $a\in H^{0,\frac{1}{2}}\cap\mathcal{H}$, we have
    \begin{equation}\label{ANS-E4.16}
    \|\partial_3S^v_la\|_{L^2}\lesssim c_l2^\frac{3l}{2}\|a\|_{\mathcal{H}}
    \ \textrm{ and }\ \|\partial_3S^v_la\|_{L^2}\lesssim
    c_l2^\frac{l}{2}\|a\|_{H^{0,\frac{1}{2}}},
    \end{equation}
using the similar argument to that in the proof of
(\ref{ANS-E4.15}), we have
    \begin{equation}
      \sum_j2^{-j}|R^{v,1}_j|\leq\frac{\nu_h}{100}\|\nabla_ha\|^2_{\mathcal{H}}
      +C_1(a,b)\mu(\|a\|^2_{\mathcal{H}}).
    \end{equation}

    Now let us estimate $R^{v,2}_j$  in high (vertical) frequencies by using that $a$ and
    $\nabla_ha$ are in $H^{0,\frac{1}{2}}$. Using Lemma
    \ref{ANS-L2.1}, we obtain
        \begin{eqnarray*}
          \|\Delta^{vi}_j(\widetilde{\Delta}^v_lb^h\partial_3S^v_{l+2}a)\|_{L^2}
          &\lesssim& 2^{\frac{j}{2}}\|\widetilde{\Delta}^v_lb^h\partial_3S^v_{l+2}a\|_{L^2_h(L^1_v)}
          \\
            &\lesssim&2^{\frac{j}{2}}\|\widetilde{\Delta}^v_lb^h
            \|_{L^{2p}_h(L^2_v)}\|\partial_3S^v_{l+2}a\|_{L^\frac{2p}{p-1}_h(L^2_v)}.
        \end{eqnarray*}
Using (\ref{ANS-E4.4}) and (\ref{ANS-E4.16}), we obtain
    $$
    \|\partial_3S^v_{l+2}a\|_{L^\frac{2p}{p-1}_h(L^2_v)}
    \lesssim c_l2^\frac{l}{2}\|\nabla_ha\|^{\frac{1}{p}}_{H^{0,\frac{1}{2}}}
    \|a\|^{1-\frac{1}{p}}_{H^{0,\frac{1}{2}}}
    $$
and
    \begin{equation}
      |R^{v,2}_j|\lesssim \|b^h\|_{\widetilde{L}^{2p}_h(B^\frac{1}{2}_v)}
            \|\nabla_ha\|_{\mathcal{H}}\|a\|^{1-\frac{1}{p}}_{H^{0,\frac{1}{2}}}
            \|\nabla_ha\|^{\frac{1}{p}}_{H^{0,\frac{1}{2}}}.
            \label{ANS-E4.18}
    \end{equation}

Then, let us estimate $R^{v,2}_j$  in low (vertical) frequencies
by using that $a$ and
    $\nabla_ha$ are in $\mathcal{H}$. Using Lemma
    \ref{ANS-L4.4}, (\ref{ANS-E4.4}) and (\ref{ANS-E4.16}), we have
        \begin{eqnarray*}
          \|\Delta^{vi}_j(\widetilde{\Delta}^v_lb^h\partial_3S^v_{l+2}a)\|_{L^2}
          &\lesssim&2^\frac{j}{2}\|\widetilde{\Delta}^v_lb\|_{L^{p'}_h(L^2_v)}
          \|\partial_3S^v_{l+2}a\|_{L^\frac{2p'}{p'-2}_h(L^2_v)}\\
                &\lesssim&2^\frac{j}{2}2^ld_l\sqrt{p'}\|b\|^\frac{2}{p'}_{\mathcal{B}}
                \|\nabla_hb\|^{1-\frac{2}{p'}}_{\mathcal{B}}
          \|a\|_{\mathcal{H}}^{1-\frac{2}{p'}}
                    \|\nabla_ha\|_{\mathcal{H}}^{\frac{2}{p'}}.
        \end{eqnarray*}
Thus, we deduce that
    $$
    |R^{v,2}_j|\lesssim c_j2^j\sqrt{p'}
    \|b\|^\frac{2}{p'}_{\mathcal{B}}
                \|\nabla_hb\|^{1-\frac{2}{p'}}_{\mathcal{B}}
          \|a\|_{\mathcal{H}}^{1-\frac{2}{p'}}
                    \|\nabla_ha\|_{\mathcal{H}}^{1+\frac{2}{p'}}.
    $$
Combining it with (\ref{ANS-E4.18}), using  the similar argument
to that in the proof of (\ref{ANS-E4.15}), we have
\begin{equation}
        \sum_j2^{-j}|(R^h_j|\Delta^{vi}_ja)_{L^2}|
        \leq \frac{\nu_h}{100}\|\nabla_ha\|_{\mathcal{H}}^2
        +C_1(a,b)\mu(\|a\|_{\mathcal{H}}).
    \end{equation}
This proves the estimate of the term $(b\cdot\nabla a
  |a)_{\mathcal{H}}$.

 \textbf{The estimate of the term $(a\cdot\nabla b
  |a)_{\mathcal{H}}$.} Using  Bony's paradifferential decomposition
  in the vertical variable and in the inhomogeneous context, we
  have
        $$
        a\cdot \nabla b =T_a\nabla b+{R}(a,\nabla b)+T_{\nabla b}a
        $$
  with
    $$
    T_a\nabla b:=\sum_{l\geq1}S^{vi}_{l-1}a\cdot \nabla \Delta^{vi}_lb,
    \ {R}(a,\nabla b):=\sumetage{l\geq-1}{i=-1,0,1}\Delta^{vi}_la\cdot\nabla
    \Delta^{vi}_{l+i}b
    $$
        $$\textrm{and }
    \ T_{\nabla b}a:=\sum_{l\geq1}\Delta^{vi}_la\cdot \nabla S^{vi}_{l-1}
    b.
    $$
   \noindent\textit{Step 3.} The estimate of $(T_a\nabla
    b|a)_{\mathcal{H}}$. Using the similar argument to that in Step 1, we have
            $$
        \Delta^{vi}_j(T_a\nabla b)=\overline{T}^h_j+\overline{T}^v_j,
            $$
    with
        $$
        \overline{T}^h_j:=\Delta^{vi}_j\sum_{|j-l|\leq5}S^{vi}_{l-1}a^h\cdot\nabla_h\Delta^{vi}_lb
      \ \textrm{ and }\  \overline{T}^v_j:=\Delta^{vi}_j\sum_{|j-l|\leq5}S^{vi}_{l-1}a^3\partial_3\Delta^{vi}_lb.
        $$
  \noindent\textit{Step 3a.} The estimate of $\sum_j2^{-j}(\overline{T}^h_j|\Delta^{vi}_ja)_{L^2}$.
        Then, we have
            $$
        (\overline{T}^h_j|\Delta^{vi}_ja)_{L^2}=\overline{T}^{h1}_j+\overline{T}^{h2}_j,
            $$
        with
         $$
        \overline{T}^{h1}_j:=-\int_{\mathbb{R}^3}
        \Delta^{vi}_j(\sum_{|j-l|\leq5}S^{vi}_{l-1}\mathrm{div}_ha^h\Delta^{vi}_lb)
        \Delta^{vi}_jadx,
        $$
        $$
        \overline{T}^{h2}_j:=-\int_{\mathbb{R}^3}
        \Delta^{vi}_j(\sum_{|j-l|\leq5}S^{vi}_{l-1}a^h\Delta^{vi}_lb)\cdot\nabla_h
        \Delta^{vi}_jadx.
        $$
Using H\"{o}lder's inequality and (\ref{ANS-E4.4}), we obtain
        \begin{eqnarray*}
          |\overline{T}^{h1}_j|&\lesssim& \sum_{|j-l|\leq5}\|\nabla_hS^{vi}_{l-1}a\|_{L^2}\|b\|_{L^{2p}_h(L^\infty_v)}
          \|\Delta^{vi}_ja\|_{L^\frac{2p}{p-1}_h(L^2_v)}\\
                    &\lesssim& \sum_{|j-l|\leq5}c_l2^{\frac{l}{2}}
                        \|\nabla_ha\|_{\mathcal{H}}
                    \|b\|_{L^{2p}_h(L^\infty_v)}
          c_j2^\frac{j}{2}\|a\|_{\mathcal{H}}^{1-\frac{1}{p}}\|\nabla_ha\|_{\mathcal{H}}^{\frac{1}{p}}\\
             &\lesssim& d_j2^j
                    \|a\|_{\mathcal{H}}^{1-\frac{1}{p}}\|\nabla_ha\|_{\mathcal{H}}^{1+\frac{1}{p}}
                    \|b\|_{L^{2p}_h(L^\infty_v)},
        \end{eqnarray*}
            \begin{eqnarray*}
          |\overline{T}^{h2}_j|&\lesssim& \sum_{|j-l|\leq5}\|S^{vi}_{l-1}a\|_{L^\frac{2p}{p-1}_h(L^2_v)}\|b\|_{L^{2p}_h(L^\infty_v)}
          \|\nabla_h\Delta^{vi}_ja\|_{L^2}\\
                    &\lesssim& \sum_{|j-l|\leq5}c_l2^{\frac{l}{2}}
                    \|a\|_{\mathcal{H}}^{1-\frac{1}{p}}\|\nabla_ha\|_{\mathcal{H}}^{\frac{1}{p}}
                    \|b\|_{L^{2p}_h(L^\infty_v)}
          c_j2^\frac{j}{2}\|\nabla_ha\|_{\mathcal{H}}\\
             &\lesssim& d_j2^j
                    \|a\|_{\mathcal{H}}^{1-\frac{1}{p}}\|\nabla_ha\|_{\mathcal{H}}^{1+\frac{1}{p}}
                    \|b\|_{L^{2p}_h(L^\infty_v)}.
        \end{eqnarray*}
\noindent\textit{Step 3b.} The estimate of
$\sum_j2^{-j}(\overline{T}^v_j|\Delta^{vi}_ja)_{L^2}$. Let
    $$
    \overline{T}^v_j:=\overline{T}^{v0}_j+\overline{T}^{v1}_j,
    $$
with
    $$
\overline{T}^{v0}_j:=\Delta^{vi}_j\sum_{|j-l|\leq5}S^{vi}_{0}a^3\partial_3\Delta^{vi}_lb,
    $$
      $$
\overline{T}^{v1}_j:=\Delta^{vi}_j\sumetage{|j-l|\leq5}{l'\in[0,l-2]}\Delta^{vi}_{l'}a^3\partial_3\Delta^{vi}_lb.
    $$
Using H\"{o}lder's inequality and (\ref{ANS-E4.4}), we obtain
        \begin{eqnarray*}
          |(\overline{T}^{v0}_j|\Delta^{vi}_ja)_{L^2}|&\lesssim& \sum_{|j-l|\leq5}\|S^{vi}_{
          0}a\|_{L^2_h(L^\infty_v)}
          \|\partial_3\Delta^{vi}_lb\|_{L^{2p}_h(L^2_v)}
          \|\Delta^{vi}_ja\|_{L^\frac{2p}{p-1}_h(L^2_v)}\\
                    &\lesssim& \sum_{|j-l|\leq5}
                        \|a\|_{\mathcal{H}}
                    d_l2^\frac{l}{2}\|b\|_{\widetilde{L}^{2p}_h(B^\frac{1}{2}_v)}
          c_j2^\frac{j}{2}\|a\|_{\mathcal{H}}^{1-\frac{1}{p}}\|\nabla_ha\|_{\mathcal{H}}^{\frac{1}{p}}\\
             &\lesssim& d_j2^j
                    \|a\|_{\mathcal{H}}^{2-\frac{1}{p}}\|\nabla_ha\|_{\mathcal{H}}^{\frac{1}{p}}
                    \|b\|_{\widetilde{L}^{2p}_h(B^\frac{1}{2}_v)}.
        \end{eqnarray*}
From (\ref{ANS-E4.9}), we have
        $$
        |\sum_{j}2^{-j}(\overline{T}^{v1}_j|\Delta^{vi}_ja)_{L^2}|
        \leq \overline{T}^{v1,N}+\overline{T}^{v1}_N,
        $$
with
    $$
    \overline{T}^{v1,N}=\sum_{j\geq N}\sumetage{|j-l|\leq5}{l'\geq N}2^{-j}\left|\int_{\mathbb{R}^3}\Delta^{vi}_j(
            2^{-l'}\widetilde{\Delta}^{vi}_{l'}\mathrm{div}_ha^h\partial_3\Delta^{vi}_lb)\Delta^{vi}_jadx
            \right|,
    $$
    $$
    \overline{T}^{v1}_N=\sum_{j}\sumetage{|j-l|\leq5}{l'\in[0,N]}2^{-j}\left|\int_{\mathbb{R}^3}\Delta^{vi}_j(
            2^{-l'}\widetilde{\Delta}^{vi}_{l'}\mathrm{div}_ha^h\partial_3\Delta^{vi}_lb)\Delta^{vi}_jadx
            \right|.
    $$
From (\ref{ANS-E4.4}), Lemma \ref{ANS-L4.4} and H\"{o}lder's
inequality, we obtain
        \begin{eqnarray*}
          \overline{T}^{v1}_N
            &\lesssim&
            \sum_{j}2^{-j}\sumetage{|j-l|\leq5}{l'\in[0,N]}2^{-l'}\|\nabla_h\Delta^{vi}_{l'}
          a\|_{L^2_h(L^\infty_v)}
          \|\partial_3\Delta^{vi}_lb\|_{L^{p'}_h(L^2_v)}
          \|\Delta^{vi}_ja\|_{L^\frac{2p'}{p'-2}_h(L^2_v)}\\
          &\lesssim&\sum_{j}2^{-j}
            \sumetage{|j-l|\leq5}{l'\in[0,N]}c_{l'}\|\nabla_h
          a\|_{\mathcal{H}}c_l2^\frac{l}{2}\sqrt{p'}
          \|b\|^\frac{2}{p'}_{\mathcal{B}}\|\nabla_hb\|^{1-\frac{2}{p'}}_{\mathcal{B}}
         c_j
         2^\frac{j}{2}\|a\|^{1-\frac{2}{p'}}_{\mathcal{H}}\|\nabla_ha\|^{\frac{2}{p'}}_{\mathcal{H}}\\
                    &\lesssim&\sqrt{p'N}
          \|b\|^\frac{2}{p'}_{\mathcal{B}}\|\nabla_hb\|^{1-\frac{2}{p'}}_{\mathcal{B}}
         \|a\|^{1-\frac{2}{p'}}_{\mathcal{H}}\|\nabla_ha\|^{1+\frac{2}{p'}}_{\mathcal{H}},
         \ \textrm{ (low vertical frequencies) }
        \end{eqnarray*}and
        \begin{eqnarray*}
          \overline{T}^{v1,N}
            &\lesssim&
            \sum_{j\geq N}2^{-\frac{j}{2}}\sumetage{|j-l|\leq5}{l'\geq N}2^{-l'}\|\nabla_h\Delta^{vi}_{l'}
          a\|_{L^2_h(L^2_v)}
          \|\partial_3\Delta^{vi}_lb\|_{L^{2p}_h(L^2_v)}
          \|\Delta^{vi}_ja\|_{L^\frac{2p}{p-1}_h(L^2_v)}\\
          &\lesssim&\sum_{j\geq N}2^{-\frac{j}{2}}
            \sumetage{|j-l|\leq5}{l'\geq N}2^{-\frac{l'}{2}}c_{l'}\|\nabla_h
          a\|_{\mathcal{H}}d_l2^{\frac{l}{2}}
          \|b\|_{\widetilde{L}^{2p}_h(B^\frac{1}{2}_v)}
         2^{-\frac{j}{2}}\|a\|^{1-\frac{1}{p}}_{H^{0,\frac{1}{2}}}
         \|\nabla_ha\|^{\frac{1}{p}}_{H^{0,\frac{1}{2}}}\\
                    &\lesssim&2^{-N}\|\nabla_h
          a\|_{\mathcal{H}}
          \|b\|_{\widetilde{L}^{2p}_h(B^\frac{1}{2}_v)}
         \|a\|^{1-\frac{1}{p}}_{H^{0,\frac{1}{2}}}
         \|\nabla_ha\|^{\frac{1}{p}}_{H^{0,\frac{1}{2}}},
         \ \textrm{ (high vertical frequencies). }
        \end{eqnarray*}
Using the similar argument to that in the proof of
(\ref{ANS-E4.15}), we get
        $$
        |\sum_{j}2^{-j}(\overline{T}^{v1}_j|\Delta^{vi}_ja)_{L^2}|
        \leq \frac{\nu_h}{200}\|\nabla_ha\|_{\mathcal{H}}^2
        +C_1(a,b)\mu(\|a\|_{\mathcal{H}}).
    $$
Thus, from above estimates, we get
    $$
    (T_a\nabla b|a)_{\mathcal{H}} \leq\frac{\nu_h}{100}\|\nabla_ha\|_{\mathcal{H}}^2+C(a,b)\mu(\|a\|_{\mathcal{H}}^2).
    $$

\noindent\textit{Step 4.} The estimate of $(R(a,\nabla
b)|a)_{\mathcal{H}}$. Obviously, we have
    $$
    (\Delta^{vi}_jR(a,\nabla
b)|\Delta^{vi}_ja)_{L^2}=\overline{R}^h_j+\overline{R}^v_j+\overline{R}^0_j,
    $$
with
    $$
    \overline{R}^h_j:=(\Delta^{vi}_j(\sumetage{l\geq 1+(j-N_0)^+}{i=-1,0,1}
    \Delta^{vi}_la^h\cdot\nabla_h\Delta^{vi}_{l+i}b)
    |\Delta^{vi}_ja)_{L^2},
    $$
    $$
    \overline{R}^v_j:=(\Delta^{vi}_j(\sumetage{l\geq 1+(j-N_0)^+}{i=-1,0,1}
    \Delta^{vi}_la^3\partial_3\Delta^{vi}_{l+i}b)
    |\Delta^{vi}_ja)_{L^2},
    $$
    $$
    \overline{R}^0_j:=(\Delta^{vi}_j(
    S^{vi}_{0}a\cdot\nabla S^{vi}_{1}b)
    |\Delta^{vi}_ja)_{L^2}+(\Delta^{vi}_j(
    \Delta^{vi}_{0}a\cdot\nabla S^{vi}_{2}b)
    |\Delta^{vi}_ja)_{L^2}.
    $$
It is obvious that if $j$ is large enough, $\overline{R}^0_j$ is
0. Thus, if $j\leq N_1$, we obtain
    $$
      |\overline{R}^0_j|
            \lesssim\|b\|_{L^{2p}_h(L^\infty_v)}\left(\|a\|^{2-\frac{1}{p}}_{\mathcal{H}}
            \|\nabla_ha\|^{\frac{1}{p}}_{\mathcal{H}}+
            \|a\|^{1-\frac{1}{p}}_{\mathcal{H}}
            \|\nabla_ha\|^{1+\frac{1}{p}}_{\mathcal{H}}\right).
    $$
  \noindent\textit{Step 4a.} The estimate of $\sum_j2^{-j}\overline{R}^h_j$.
        Then, we have
            $$
        \overline{R}^h_j=\overline{R}^{h1}_j+\overline{R}^{h2}_j,
            $$
        with
         $$
        \overline{R}^{h1}_j:=-\int_{\mathbb{R}^3}
        \Delta^{vi}_j(\sumetage{l\geq 1+(j-N_0)^+}{i=-1,0,1}\Delta^{vi}_{l}\mathrm{div}_ha^h\Delta^{vi}_{l+i}b)
        \Delta^{vi}_jadx,
        $$
        $$
        \overline{R}^{h2}_j:=-\int_{\mathbb{R}^3}
        \Delta^{vi}_j(\sumetage{l\geq 1+(j-N_0)^+}{i=-1,0,1}\Delta^{vi}_{l}a^h\Delta^{vi}_{l+i}b)\cdot\nabla_h
        \Delta^{vi}_jadx.
        $$
Using H\"{o}lder's inequality and (\ref{ANS-E4.4}), we obtain
        \begin{eqnarray*}
          |\overline{R}^{h1}_j|&\lesssim& 2^\frac{j}{2}\sumetage{l\geq 1+(j-N_0)^+}{i=-1,0,1}
          \|\nabla_h\Delta^{vi}_{l}a\|_{L^2}\|\Delta^{vi}_{l+i}b\|_{L^{p'}_h(L^2_v)}
          \|\Delta^{vi}_ja\|_{L^\frac{2p'}{p'-2}_h(L^2_v)}\\
                    &\lesssim& 2^\frac{j}{2}\sumetage{l\geq 1+(j-N_0)^+}{i=-1,0,1}c_l2^{\frac{l}{2}}
                        \|\nabla_ha\|_{\mathcal{H}}
                    c_{l+i}2^{-\frac{l}{2}}\sqrt{p'}\|b\|^{\frac{2}{p'}}_{\mathcal{B}}
                    \|\nabla_hb\|^{1-\frac{2}{p'}}_{\mathcal{B}}
          c_j2^\frac{j}{2}\|a\|_{\mathcal{H}}^{1-\frac{2}{p'}}\|\nabla_ha\|_{\mathcal{H}}^{\frac{2}{p'}}\\
             &\lesssim& c_j\sqrt{p'}2^j
                    \|a\|_{\mathcal{H}}^{1-\frac{1}{p'}}\|\nabla_ha\|_{\mathcal{H}}^{1+\frac{1}{p'}}
                   \|b\|^{\frac{2}{p'}}_{\mathcal{B}}
                    \|\nabla_hb\|^{1-\frac{2}{p'}}_{\mathcal{B}}, \textrm{
                    (low vertical
                    frequencies)},
        \end{eqnarray*}
        \begin{eqnarray*}
          |\overline{R}^{h1}_j|&\lesssim& 2^\frac{j}{2}\sumetage{l\geq 1+(j-N_0)^+}{i=-1,0,1}
          \|\nabla_h\Delta^{vi}_{l}a\|_{L^2}\|\Delta^{vi}_{l+i}b\|_{L^{2p}_h(L^2_v)}
          \|\Delta^{vi}_ja\|_{L^\frac{2p}{p-1}_h(L^2_v)}\\
                    &\lesssim& 2^\frac{j}{2}\sumetage{l\geq 1+(j-N_0)^+}{i=-1,0,1}c_l2^{\frac{l}{2}}
                        \|\nabla_ha\|_{\mathcal{H}}
                    d_{l+i}2^{-\frac{l}{2}}\|b\|_{\widetilde{L}^{2p}_h(B^\frac{1}{2}_v)}
          c_j2^{-\frac{j}{2}}\|a\|_{H^{0,\frac{1}{2}}}^{1-\frac{1}{p}}\|\nabla_ha\|_{H^{0,\frac{1}{2}}}^{\frac{1}{p}}\\
             &\lesssim&
                        \|\nabla_ha\|_{\mathcal{H}}
                    \|b\|_{\widetilde{L}^{2p}_h(B^\frac{1}{2}_v)}
          \|a\|_{H^{0,\frac{1}{2}}}^{1-\frac{1}{p}}\|\nabla_ha\|_{H^{0,\frac{1}{2}}}^{\frac{1}{p}}. \textrm{
                    (high vertical
                    frequencies)}.
        \end{eqnarray*}
Using the similar argument to that in the proof of
(\ref{ANS-E4.15}), we get
     \begin{equation}
        |\overline{R}^{h1}_j|
        \leq \frac{\nu_h}{100}\|\nabla_ha\|_{\mathcal{H}}^2
        +C_1(a,b)\mu(\|a\|_{\mathcal{H}}).\label{ANS-E4.22-1}
     \end{equation}
Similarly, we have
    $$
        |\overline{R}^{h2}_j|
        \leq \frac{\nu_h}{100}\|\nabla_ha\|_{\mathcal{H}}^2
        +C_1(a,b)\mu(\|a\|_{\mathcal{H}}).
    $$
From (\ref{ANS-E4.9}), we have
    $$
    \overline{R}^{v}_j=-(\Delta^{vi}_j(\sumetage{l\geq 1+(j-N_0)^+}{i=-1,0,1}
    2^{-l}\widetilde{\Delta}^{vi}_l\mathrm{div}_ha^h\partial_3\Delta^{vi}_{l+i}b)
    |\Delta^{vi}_ja)_{L^2},
    $$
and using the similar argument to that in the proof of
(\ref{ANS-E4.22-1}), we obtain
    $$
        |\overline{R}^{v}_j|
        \leq \frac{\nu_h}{100}\|\nabla_ha\|_{\mathcal{H}}^2
        +C_1(a,b)\mu(\|a\|_{\mathcal{H}}).
    $$
     \noindent\textit{Step 5.} The estimate of $(T_{\nabla
     b}a|a)_{\mathcal{H}}$. Using  similar arguments to that in the proof
     of Steps 3-4, we get
        $$
        (T_{\nabla b}a|a)_{\mathcal{H}}
        \leq\frac{\nu_h}{100}\|\nabla_ha\|_{\mathcal{H}}^2+C(a,b)\mu(\|a\|_{\mathcal{H}}^2).
        $$
 This proves Lemma \ref{ANS-L4.5}.
{\hfill $\square$\medskip}

Then, we will prove that $\delta\in L^\infty_{T_1}(\mathcal{H})$
and $\nabla_h \delta \in L^2_{T_1}(\mathcal{H})$ in the following
lemma.

\begin{lem}\label{ANS-L4.6} We have
$\delta\in L^\infty_{T_1}(\mathcal{H})$ and $\nabla_h \delta \in
L^2_{T_1}(\mathcal{H})$.
\end{lem}
\noindent\textbf{Proof.}
 Since $\delta\in B^{0,\frac{1}{2}}$, we only need to prove that
    $S^v_0\delta\in L^\infty_{T_1}(L^2)$ and $\nabla_hS^v_0\delta\in
    L^2_{T_1}(L^2)$.

  Let us write that $S^v_0\delta$ is a solution (with initial value
  0) of
    $$
    \partial_tS^v_0\delta-\nu_h\Delta_hS^v_0\delta-\nu_3\partial_3^2S^v_0\delta
    =\sum_{i=1}^6g_i-\nabla S^v_0P,
    $$
  with
    $$
    g_1:=-S^v_0\partial_3(\delta^3u_1),
    $$
        $$
        g_2:=-S^v_0\mathrm{div}_h(u_1(Id-S^v_0)\delta^h),
        $$
            $$
            g_3:=-S^v_0\mathrm{div}_h(u_1S^v_0\delta^h),
            $$
    $$
    g_4:=-S^v_0\partial_3(u_2^3\delta),
    $$
        $$
        g_5:=-S^v_0\mathrm{div}_h(u_2^h(Id-S^v_0)\delta),
        $$
            $$
            g_6:=-S^v_0\mathrm{div}_h(u_2^hS^v_0\delta).
            $$
  Using Lemmas \ref{ANS-L2.1} and \ref{ANS-L4.2}, we have
        \begin{eqnarray}
          \|g_1\|_{L^2_{T_1}(L^2)}&\lesssim&
          \sum_{j\leq-1}2^j\|\Delta^v_j(u_1
          \delta^3)\|_{L^2_{T_1}(L^2(\mathbb{R}^3))}
          \nonumber\\
        &\lesssim&\sum_{j\leq-1}2^j
        d_j\nu_h^{-\frac{1}{2p}}2^{-\frac{j}{2}}
        \|\delta\|_{B^{0,\frac{1}{2}}(T_1)}
        \|u_1\|_{\widetilde{L}^\frac{2p}{p-1}_{T_1}(\widetilde{L}^{2p}_h(B^\frac{1}{2}_v))}
        \nonumber\\
                        &\lesssim&\nu_h^{-\frac{1}{2p}}
        \|\delta\|_{B^{0,\frac{1}{2}}(T_1)}
        \|u_1\|_{\widetilde{L}^\frac{2p}{p-1}_{T_1}(\widetilde{L}^{2p}_h(B^\frac{1}{2}_v))}.
        \label{ANS-E4.20}
        \end{eqnarray}
  Similarly, we have
        \begin{equation}
        \|g_4\|_{L^2_{T_1}(L^2)}\lesssim \nu_h^{-\frac{1}{2p}}
        \|\delta\|_{B^{0,\frac{1}{2}}(T_1)}
        \|u_2\|_{\widetilde{L}^\frac{2p}{p-1}_{T_1}(\widetilde{L}^{2p}_h(B^\frac{1}{2}_v))}.
        \label{ANS-E4.21}
        \end{equation}
By (\ref{ANS-E4.4}) and Young's inequality, we obtain
        \begin{eqnarray*}
          \|(Id-S^v_0)\delta\|_{L^{2p}_{T_1}(L^\frac{2p}{p-1}_h(L^2_v))}
          &\lesssim&\sum_{j\geq0}\|\Delta^v_j\delta\|_{L^{2p}_{T_1}(L^\frac{2p}{p-1}_h(L^2_v))}\\
            &\lesssim&\sum_{j\geq0}\|\Delta^v_j\delta\|_{L^\infty_{T_1}(L^2_h(L^2_v))}^{1-\frac{1}{p}}
            \|\Delta^v_j\nabla_h\delta\|_{L^2_{T_1}(L^2_h(L^2_v))}^{\frac{1}{p}}\\
                &\lesssim&\nu_h^{-\frac{1}{2p}}\sum_{j\geq0}\left(\|\Delta^v_j\delta\|_{L^\infty_{T_1}(L^2_h(L^2_v))}+
           \nu_h^\frac{1}{2} \|\Delta^v_j\nabla_h\delta\|_{L^2_{T_1}(L^2_h(L^2_v))}\right)\\
            &\lesssim&\nu_h^{-\frac{1}{2p}}\|\delta\|_{B^{0,\frac{1}{2}}(T_1)}.
        \end{eqnarray*}
Then, we have
    \begin{equation}
      g_2=\mathrm{div}_h\widetilde{g_2}
      \ \textrm{ and }\ g_5=\mathrm{div}_h\widetilde{g_5},
      \label{ANS-E4.22}
    \end{equation}
with
    $$
    \| \widetilde{g_2}\|_{L^2_t(L^2)}\lesssim
    \nu_h^{-\frac{1}{2p}}
        \|\delta\|_{B^{0,\frac{1}{2}}(t)}
        \|u_1\|_{\widetilde{L}^\frac{2p}{p-1}_{t}(\widetilde{L}^{2p}_h(B^\frac{1}{2}_v))},
    $$
    $$
   \| \widetilde{g_5}\|_{L^2_t(L^2)}\lesssim
    \nu_h^{-\frac{1}{2p}}
        \|\delta\|_{B^{0,\frac{1}{2}}(t)}
        \|u_2\|_{\widetilde{L}^\frac{2p}{p-1}_{t}(\widetilde{L}^{2p}_h(B^\frac{1}{2}_v))}.
    $$
    The terms $g_3$ and $g_6$ must be treated with a commutator
    argument based on the following estimate, which is proved in Lemma 4.3 of
    \cite{Chemin07}:
    Let $\chi$ be a function of $\mathcal{S}(\mathbb{R})$. A
    constant $C$ exists such that, for any function $a$ in
    $L^2_h(L^\infty_v)$, we have
        \begin{equation}\label{ANS-E4.23}
        \|[\chi(\varepsilon x_3);S^v_0]a\|_{L^2}\leq
        C\varepsilon^\frac{1}{2}\|a\|_{L^2_h(L^\infty_v)}.
        \end{equation}

    Now let us choose $\chi\in\mathcal{D}(\mathbb{R})$ with value
    1 near 0 and let us state
    $$
    S^v_{0,\varepsilon}a:=\chi(\varepsilon x_3)S^v_0a.
    $$
Using a classical $L^2$ energy estimate and Young's inequality, we
    have
        \begin{eqnarray*}
        &&\|S^v_{0,\varepsilon}\delta(t)\|_{L^2}^2
        +\nu_h\int^t_0\|\nabla_hS^v_{0,\varepsilon}\delta(s)\|_{L^2}^2ds
        +2\nu_3\int^t_0\|\partial_3S^v_{0,\varepsilon}\delta(s)\|_{L^2}^2ds\\
            &\lesssim& \int^t_0\|S^v_{0,\varepsilon}\delta(s)\|_{L^2
            }\left(\|g_1(s)\|_{L^2}
            +\|g_4(s)\|_{L^2}\right)ds+\frac{1}{\nu_h}\int^t_0
            \left(\|\widetilde{g_2}(s)\|^2_{L^2}
            +\|\widetilde{g_5}(s)\|_{L^2}^2\right)ds\\
                    &&
                    +\int^t_0\int_{\mathbb{R}^3}\chi(\varepsilon
                    x_3)(g_3+g_6)(s)S^v_{0,\varepsilon}\delta(s)dxds
            +\nu_3\varepsilon^2\int^t_0\|\chi'(\varepsilon x_3)S^v_0\delta(s)\|_{L^2}^2ds.
        \end{eqnarray*}
By the definition of $g_3$, we have
    $$
        \int^t_0\int_{\mathbb{R}^3}\chi(\varepsilon
                    x_3)g_3(s)S^v_{0,\varepsilon}\delta(s)dxds
                    = \int^t_0D_1(s)ds+ \int^t_0D_2(s)ds,
    $$
with
    $$
    D_1:=\int_{\mathbb{R}^3}[\chi(\varepsilon
    x_3);S^v_0](u_1S^v_0\delta^h)\cdot
    \nabla_hS^v_{0,\varepsilon}\delta dx,
    $$
        $$
    D_2:=\int_{\mathbb{R}^3}S^v_0(u_1S^v_{0,\varepsilon}\delta^h)\cdot
    \nabla_hS^v_{0,\varepsilon}\delta dx.
    $$
From Lemma \ref{ANS-L4.0}, (\ref{ANS-E4.4}), (\ref{ANS-E4.23}),
H\"{o}lder's inequality and Young's inequality, we obtain
    \begin{eqnarray*}
      \int^t_0|D_1(s)|ds&\lesssim&\varepsilon^\frac{1}{2}\|u_1S^v_0\delta^h\|_{L^2_t(L^2_h(L^\infty_v))}
      \|\nabla_hS^v_{0,\varepsilon}\delta\|_{L^2_t(L^2)}\\
            &\lesssim&\varepsilon^\frac{1}{2}\|u_1\|_{L^\frac{2p}{p-1}_t(L^{2p}_h(L^\infty_v
            ))}\|S^v_0\delta^h\|_{L^{2p}_t(L^\frac{2p}{p-1}_h(L^\infty_v))}
      \|\nabla_hS^v_{0,\varepsilon}\delta\|_{L^2_t(L^2)}\\
        &\leq&\frac{\nu_h}{10}\|\nabla_hS^v_{0,\varepsilon}\delta\|_{L^2_t(L^2)}^2
        +C\nu_h^{-1-\frac{1}{p}}\varepsilon\|u_1\|_{\widetilde{L}^\frac{2p}{p-1}_t(\widetilde{L}^{2p}_h(B^\frac{1}{2}_v
            ))}^2\|S^v_0\delta^h\|_{B^{0,\frac{1}{2}}(t)}^2,
    \end{eqnarray*}
        \begin{eqnarray*}
          |D_2|&\lesssim&\|u_1\|_{L^{2p}_h(L^\infty_v)}
          \|S^v_{0,\varepsilon}\delta\|_{L^\frac{2p}{p-1}_h(L^2_v)}
        \|\nabla_hS^v_{0,\varepsilon}\delta\|_{L^2}\\
            &\lesssim&\|u_1\|_{L^{2p}_h(L^\infty_v)}
          \|S^v_{0,\varepsilon}\delta\|^{1-\frac{1}{p}}_{L^2}
        \|\nabla_hS^v_{0,\varepsilon}\delta\|^{1+\frac{1}{p}}_{L^2}\\
                  &\leq&\frac{\nu_h}{10}\|\nabla_hS^v_{0,\varepsilon}\delta\|_{L^2}^2
        +C\nu_h^{-\frac{p+1}{p-1}}
        \|u_1\|_{L^{2p}_h(L^\infty_v)}^\frac{2p}{p-1}\|S^v_{0,\varepsilon}\delta^h\|_{L^2}^2
        \end{eqnarray*}
and
    \begin{eqnarray}
       &&\int^t_0\int_{\mathbb{R}^3}\chi(\varepsilon
                    x_3)g_3(s)S^v_{0,\varepsilon}\delta(s)dxds\nonumber\\
       &\leq&\frac{\nu_h}{5}\|\nabla_hS^v_{0,\varepsilon}\delta\|_{L^2_t(L^2)}^2
        +C\nu_h^{-1-\frac{1}{p}}\varepsilon
        \|u_1\|_{\widetilde{L}^\frac{2p}{p-1}_t(\widetilde{L}^{2p}_h(B^\frac{1}{2}_v
            ))}^2\|S^v_0\delta\|_{B^{0,\frac{1}{2}}(t)}^2\nonumber\\
             &&
        +C\nu_h^{-\frac{p+1}{p-1}}\int^t_0
        \|u_1\|_{L^{2p}_h(L^\infty_v)}^\frac{2p}{p-1}\|S^v_{0,\varepsilon}\delta\|_{L^2}^2ds.
        \label{ANS-E4.24}
    \end{eqnarray}
Similarly, we have
        \begin{eqnarray}
       &&\int^t_0\int_{\mathbb{R}^3}\chi(\varepsilon
                    x_3)g_6(s)S^v_{0,\varepsilon}\delta(s)dxds\nonumber\\
       &\leq&\frac{\nu_h}{5}\|\nabla_hS^v_{0,\varepsilon}\delta\|_{L^2_t(L^2)}^2
        +C\nu_h^{-1-\frac{1}{p}}\varepsilon\|u_2\|_{\widetilde{L}^\frac{2p}{p-1}_t(\widetilde{L}^{2p}_h(B^\frac{1}{2}_v
            ))}^2\|S^v_0\delta\|_{B^{0,\frac{1}{2}}(t)}^2\nonumber\\
             &&
        +C\nu_h^{-\frac{p+1}{p-1}}\int^t_0
        \|u_2\|_{L^{2p}_h(L^\infty_v)}^\frac{2p}{p-1}\|S^v_{0,\varepsilon}\delta\|_{L^2}^2ds.
 \label{ANS-E4.25}
    \end{eqnarray}

From (\ref{ANS-E4.20})-(\ref{ANS-E4.22}),
(\ref{ANS-E4.24})-(\ref{ANS-E4.25}) and Lemma \ref{ANS-L2.3}, we
have, $t\in[0,T_1]$,
\begin{eqnarray*}
        &&\|S^v_{0,\varepsilon}\delta(t)\|_{L^2}^2
        +\nu_h\int^t_0\|\nabla_hS^v_{0,\varepsilon}\delta(s)\|_{L^2}^2ds
        \\
            &\leq& C_{\nu_h}(1+\varepsilon)C^4_{12}(T_1)+C_{\nu_h}\int^t_0
            \left(1+\|u_1\|_{L^{2p}_h(L^\infty_v)}^\frac{2p}{p-1}+
            \|u_2\|_{L^{2p}_h(L^\infty_v)}^\frac{2p}{p-1}
            \right)\|S^v_{0,\varepsilon}\delta\|_{L^2}^2ds\\
                    &&
            +\nu_3\varepsilon^2\int^t_0\|\chi'(\varepsilon
            x_3)\|_{L^2_v}^2\|S^v_0\delta(s)\|_{L^2_h(L^\infty_v)}^2ds\\
             &\leq& C_{\nu_h}(1+\varepsilon)C^4_{12}(T_1)+C_{\nu_h}\int^t_0
            \left(1+\|u_1\|_{L^{2p}_h(L^\infty_v)}^\frac{2p}{p-1}+
            \|u_2\|_{L^{2p}_h(L^\infty_v)}^\frac{2p}{p-1}
            \right)\|S^v_{0,\varepsilon}\delta \|_{L^2}^2ds\\
    &&+C\nu_3\varepsilon T_1C_{12}^2(T_1),
        \end{eqnarray*}
where
$C_{12}(T_1)=\|u_1\|_{B^{-1+\frac{2}{p},\frac{1}{2}}_p(T_1)}+\|u_2\|_{B^{-1+\frac{2}{p},\frac{1}{2}}_p(T_1)}
+\|\delta\|_{B^{0,\frac{1}{2}}(T_1)}$. Using Gronwall's inequality
and Lemma \ref{ANS-L2.3}, we obtain, $t\in[0,T_1]$,
    \begin{eqnarray*}
    &&\|S^v_{0,\varepsilon}\delta(t)\|_{L^2}^2
        +\nu_h\int^t_0\|\nabla_hS^v_{0,\varepsilon}\delta(s)\|_{L^2}^2ds\\
                &   \leq&
                    \left(C_{\nu_h}(1+\varepsilon)C^4_{12}(T_1)+C\nu_3\varepsilon T_1C_{12}^2(T_1)
        \right)\exp\{C_{\nu_h}(T_1+C_{12}^{\frac{2p}{p-1}}(T_1))\}.
    \end{eqnarray*}
Passing to the limit when $\varepsilon$ tends to 0 allows to
conclude the proof of this lemma. {\hfill $\square$\medskip}

\noindent\textbf{Conclusion of the proof of the uniqueness.} From
Lemmas \ref{ANS-L4.5}-\ref{ANS-L4.6}, we have
    $$
    \|\delta(t)\|^2_{\mathcal{H}}\leq
    \int^t_0f(s)\mu(\|\delta(s)\|^2_{\mathcal{H}})ds,
    \ t\in[0,T_1].
    $$
with
    $f(t):=C(\delta(t),u_1(t))+C(\delta(t),u_2(t))$. Lemma
    \ref{ANS-L2.3} and (\ref{ANS-E4.5-1})-(\ref{ANS-E4.6-1})
imply that $f\in L^1([0,T_1])$. Then, the uniqueness on $[0,T_1]$
follows from the Osgood Lemma (see for instance \cite{Fleet80}).
Since $u\in C([0,T];B^{-1+\frac{2}{p},\frac{1}{2}}_p)$, one can
    easily obtain the uniqueness of the solution $u$  on $[0,T]$.{\hfill $\square$\medskip}

\section{Continuous dependence}\label{ANS-Sec5-0}
\textbf{Proof of (\ref{ANS-E1.5}).} Here, we give a sketch proof
of (\ref{ANS-E1.5}).

 From (\ref{ANS-E4.4}), we have
    \begin{equation}
    \|a\|_{L^\frac{2p}{p-1}_h(L^2_v)}\lesssim
   \|a\|_{L^2}^{1-\frac{1}{p}}
      \|\nabla_ha\|_{L^2}^{\frac{1}{p}},\ p\geq2.\label{ANS-E-5.2}
    \end{equation}

Similar to (\ref{ANS-E4.2-0}), we obtain
    \begin{equation}
      \left\{
      \begin{array}{l}
            \delta_t-\nu_h\Delta_h\delta-\nu_3\partial_3^2\delta=-\delta\cdot\nabla u_1-u_2\cdot\nabla \delta
            -\nabla
            P,\\
                \mathrm{div}\delta=0,\\
                    \delta|_{t=0}=\delta_0:=u_{02}-u_{01},
      \end{array}
      \right.
    \end{equation}
where $\delta:=u_2-u_1$. By the $L^2$ energy estimate,
(\ref{ANS-E-5.2})  and the Cauchy-Schwarz inequality, we have, for
$p\geq2$
    \begin{eqnarray*}
      &&\frac{d}{dt}\|\delta\|^2_{L^2}+2\nu_h\|\nabla_h
      \delta\|^2_{L^2}+2\nu_3\|\partial_3
      \delta\|^2_{L^2}\\
        &\leq&2\left|\int_{\mathbb{R}^3}(
        \delta\otimes u_1): \nabla \delta dx
        \right|        +2\left|\int_{\mathbb{R}^3}(u_2\otimes
        \delta): \nabla \delta dx
        \right|\\
            &\leq&C\|\nabla
      \delta\|_{L^2}\|
        \delta\|_{L^\frac{2p}{p-1}_h(L^2_v)}\left(\| u_1\|_{L^{2p}_h(L^\infty)}
        +\| u_2\|_{L^{2p}_h(L^\infty)}\right)\\
            &\leq&C\|\nabla
      \delta\|_{L^2}\|\delta\|_{L^2}^{1-\frac{1}{p}}
      \|\nabla_h\delta\|_{L^2}^{\frac{1}{p}}\left(\| u_1\|_{L^{2p}_h(L^\infty)}
        +\| u_2\|_{L^{2p}_h(L^\infty)}\right)\\
         &\leq&C\|\nabla
      \delta\|_{L^2}^{1+\frac{1}{p}}\|\delta\|_{L^2}^{1-\frac{1}{p}}
      \left(\| u_1\|_{L^{2p}_h(L^\infty)}
        +\| u_2\|_{L^{2p}_h(L^\infty)}\right)\\
            &\leq&\nu_h\|\nabla_h
      \delta\|^2_{L^2}+\nu_3\|\partial_3
      \delta\|^2_{L^2}
       +C(\nu_h^{-\frac{p+1}{p-1}}+\nu_3^{-\frac{p+1}{p-1}})\|\delta\|_{L^2}^{2}
      \left(\sum_{i=1}^2\| u_i\|_{L^{2p}_h(L^\infty)}^\frac{2p}{p-1}
      \right).
    \end{eqnarray*}
Then, we have
    $$
    \frac{d}{dt}\|\delta\|^2_{L^2}\leq
       C(\nu_h^{-\frac{p+1}{p-1}}+\nu_3^{-\frac{p+1}{p-1}})\|\delta\|_{L^2}^{2}
      \left(\sum_{i=1}^2\| u_i\|_{L^{2p}_h(L^\infty)}^\frac{2p}{p-1}\right),\ p\geq2.
    $$
Using Gronwall's inequality and Lemma \ref{ANS-L2.3}, we obtain,
for $p\geq2$,
    \begin{eqnarray*}
    \|\delta\|^2_{L^2}&\leq&\|\delta_0\|^2_{L^2}\exp\left\{
       C(\nu_h^{-\frac{p+1}{p-1}}+\nu_3^{-\frac{p+1}{p-1}})\int^T_0
      \left(\sum_{i=1}^2\| u_i\|_{L^{2p}_h(L^\infty)}^\frac{2p}{p-1}
        \right)ds\right\}\\
        &\leq& \|\delta_0\|^2_{L^2}\exp\left
        \{C\nu_h^{-1}(\nu_h^{-\frac{p+1}{p-1}}+\nu_3^{-\frac{p+1}{p-1}})\left(
        \sum_{i=1}^2\| u_i\|_{B^{-1+\frac{2}{p},\frac{1}{2}}_p(T)
        }\right)^{\frac{2p}{p-1}}\right\}.
    \end{eqnarray*}
This finishes the proof of (\ref{ANS-E1.5}) and Theorem
\ref{ANS-T1}. {\hfill $\square$\medskip}

\section{Proof of Proposition \ref{ANS-P2}}\label{ANS-Sec5}

From Lemmas \ref{ANS-L2.1} and \ref{ANS-L2.4}, we have
            \begin{equation}
        \|\Delta^h_k\Delta^v_l a_F\|_{L^q_T(L^4_h(L^2_v))}\lesssim
        \left\{\begin{array}{ll}
               d_{k,l}\nu_h^{-\frac{1}{q}}2^{(\frac{1}{2}-\frac{2}{q})k}
                2^{-\frac{l}{2}}\|a\|_{B^{-1+\frac{2}{p},\frac{1}{2}}_p},&\textrm{for}
                \ k\geq l-1,\\
            0,&\textrm{else,}
        \end{array}
        \right.\label{ANS-E5.1}
            \end{equation}
            where  $1\leq q\leq \infty$ and $p\in[2,4]$. Then, we
            obtain,
        \begin{equation}
    \|\Delta^h_k a_F\|_{L^q(\mathbb{R}^+;L^{4}_h(L^\infty_v))}
    \lesssim \nu_h^{-\frac{1}{q}}c_k2^{-k(\frac{2}{q}-\frac{1}{2})}
    \|a\|_{B^{-1+\frac{2}{p},\frac{1}{2}}_p},  \ \textrm{ for }\
    q\in[1,\infty].\label{ANS-E5.2}
     \end{equation}

 Using Bony's decomposition in the vertical variable, we obtain
        $$
        \Delta^v_j(a_F^3\partial_3a_F)=\sum_{|j-j'|\leq5}\Delta^v_j(S^v_{j'-1}a_F^3\partial_3\Delta^v_{j'} a_F)
        +\sum_{j'\geq j-N_0}\Delta^v_j(\Delta^v_{j'}a_F^3\partial_3
        S^v_{j'+2}a_F
        ).
        $$
The two terms of the above sum are estimated exactly along the
same lines. Using Bony's decomposition in the horizontal variable,
we obtain
    $$
    S^v_{j'-1}a_F^3\partial_3\Delta^v_{j'} a_F
    =\sum_{k\geq j'-N_0}
    \left\{S^h_{k-1}S^v_{j'-1}a_F^3\partial_3\Delta^h_k\Delta^v_{j'} a_F
    +\Delta^h_kS^v_{j'-1}a_F^3\partial_3S^h_{k+2}\Delta^v_{j'} a_F
    \right\}.
    $$
 Using H\"{o}lder's inequality, Young's inequality,
Lemma \ref{ANS-L2.1} and (\ref{ANS-E5.1})-(\ref{ANS-E5.2}), we get
    \begin{eqnarray*}
      &&\|S^h_{k-1}S^v_{j'-1}a_F^3\partial_3\Delta^h_k\Delta^v_{j'}
      a_F\|_{L^1_T(L^2(\mathbb{R}^3))}\\
        &\leq& \|S^h_{k-1}S^v_{j'-1}a_F^3\|_{L^\infty_T(L^4_h(L^\infty_v))}
        \|\partial_3\Delta^h_k\Delta^v_{j'}
      a_F\|_{L^1_T(L^4_h(L^2_v))}\\
        &\lesssim&
        \sum_{k'\leq k-2}c_{k'}2^{\frac{k'}{2}}\|a\|_{B^{-1+\frac{2}{p},\frac{1}{2}}_p}
        d_{k,j'}\nu_h^{-1}2^{\frac{j'}{2}}2^{-\frac{3}{2}k}\|a\|_{B^{-1+\frac{2}{p},\frac{1}{2}}_p}\\
        &\lesssim&
        c_{k}2^{\frac{k}{2}}\|a\|_{B^{-1+\frac{2}{p},\frac{1}{2}}_p}
        d_{k,j'}\nu_h^{-1}2^{\frac{j'}{2}}2^{-\frac{3}{2}k}\|a\|_{B^{-1+\frac{2}{p},\frac{1}{2}}_p}\\
        &\lesssim&
        c_{k}d_{k,j'}
        \nu_h^{-1}2^{\frac{j'}{2}}2^{-k}\|a\|_{B^{-1+\frac{2}{p},\frac{1}{2}}_p}^2,
    \end{eqnarray*}
    \begin{eqnarray*}
      &&\|\Delta^h_{k}S^v_{j'-1}a_F^3\partial_3S^h_{k+2}\Delta^v_{j'}
      a_F\|_{L^1_T(L^2(\mathbb{R}^3))}\\
        &\leq& \|\Delta^h_{k}S^v_{j'-1}a_F^3\|_{L^1_T(L^4_h(L^\infty_v))}
        \|\partial_3S^h_{k+2}\Delta^v_{j'}
      a_F\|_{L^\infty_T(L^4_h(L^2_v))}\\
        &\lesssim&
    c_{k}\nu_h^{-1}2^{-\frac{3}{2}k}\|a\|_{B^{-1+\frac{2}{p},\frac{1}{2}}_p}
        \sum_{k'\leq k}d_{k',j'}2^{\frac{j'}{2}}2^{\frac{k'}{2}}\|a\|_{B^{-1+\frac{2}{p},\frac{1}{2}}_p}\\
        &\lesssim&
    c_{k}\nu_h^{-1}2^{-\frac{3}{2}k}\|a\|_{B^{-1+\frac{2}{p},\frac{1}{2}}_p}
        d_{k,j'}2^{\frac{j'}{2}}2^{\frac{k}{2}}\|a\|_{B^{-1+\frac{2}{p},\frac{1}{2}}_p}\\
        &\lesssim&
        c_{k}d_{k,j'}
        \nu_h^{-1}2^{\frac{j'}{2}}2^{-k}\|a\|_{B^{-1+\frac{2}{p},\frac{1}{2}}_p}^2
    \end{eqnarray*}
and
    \begin{eqnarray*}
      &&\|S^v_{j'-1}a_F^3\partial_3\Delta^v_{j'}
      a_F\|_{L^1_T(L^2(\mathbb{R}^3))}\\
        &\lesssim& \sum_{k\geq j'-N_0} c_{k}d_{k,j'}
        \nu_h^{-1}2^{\frac{j'}{2}}2^{-k}\|a\|_{B^{-1+\frac{2}{p},\frac{1}{2}}_p}^2\\
            &\lesssim&
            d_{j'}2^{-\frac{j'}{2}}\nu_h^{-1}\|a\|_{B^{-1+\frac{2}{p},\frac{1}{2}}_p}^2.
    \end{eqnarray*}
Thus, we can easily obtain
    \begin{eqnarray*}
      &&\|\Delta^v_j(a_F^3\partial_3
      a_F)\|_{L^1_T(L^2(\mathbb{R}^3))}\\
        &\lesssim&
         d_{j}2^{-\frac{j}{2}}\nu_h^{-1}\|a\|_{B^{-1+\frac{2}{p},\frac{1}{2}}_p}^2+\sum_{j'\geq j-N_0}
            d_{j'}2^{-\frac{j'}{2}}\nu_h^{-1}\|a\|_{B^{-1+\frac{2}{p},\frac{1}{2}}_p}^2\\
            &\leq&
            d_{j}2^{-\frac{j}{2}}\nu_h^{-1}\|a\|_{B^{-1+\frac{2}{p},\frac{1}{2}}_p}^2
    \end{eqnarray*}
and
    $$\sum_{j\in\mathbb{Z}}2^{\frac{j}{2}}\int^\infty_0\|\Delta^v_j(a_F^3\partial_3 a_F)\|_{L^2(\mathbb{R}^3)}dt
         \lesssim \nu_h^{-1}\|a\|_{B^{-1+\frac{2}{p},\frac{1}{2}}_p}^2.$$
Using the similar argument, we can easily estimate the term
$\sum2^{\frac{j}{2}}\int^\infty_0\|\Delta^v_j(a_F^h\cdot\nabla_h
a_F)\|_{L^2}$, and finish the proof of Proposition \ref{ANS-P2}.
{\hfill $\square$\medskip}

\section{Proof of Proposition \ref{ANS-P1}}\label{ANS-Sec6-0}
Using methods in \cite{Chemin07-0,Chemin07}, we can prove
Proposition \ref{ANS-P1} as follows.

We shall start by estimating the high frequencies. Defining a
threshold $k_0\geq0$ to be determined later on, we have
    $$
    \sum_{k=k_0}^\infty
    2^{-\sigma k}\|\Delta^h_k\phi_\varepsilon\|_{L^q(\mathbb{R}^2)}\leq
    C
    2^{-\sigma k_0}\|\phi_\varepsilon\|_{L^q(\mathbb{R}^2)}=
      C 2^{-\sigma k_0}\|\phi\|_{L^q(\mathbb{R}^2)}
        $$
On the other hand, noting that
$e^{i\frac{x_1}{\varepsilon}}=(-i\varepsilon\partial_1)^N(e^{i\frac{x_1}{\varepsilon}})$,
we get, for any $N\in\mathbb{N}$,
    $$
        \Delta^h_k
        \phi_\varepsilon=(i\varepsilon)^N2^{2k}\sum^N_{l=0}\int_{\mathbb{R}^2}
        C^l_Ne^{i\frac{y_1}{\varepsilon}}\partial_{x_1}^l({g}(2^k(x_h-y_h)))\partial_{x_1}^{N-l}
        \phi(y_h)dy_h,
    $$
where ${g}(x_h)\in \mathcal{S}(\mathbb{R}^2)$ satisfying
$\mathcal{F}{g}(\xi_h)={\varphi}(|\xi_h|)$. Young's inequality
enables us to infer that
    $$
        \|\Delta^h_k
        \phi_\varepsilon\|_{L^q}
        \leq C_\phi\varepsilon^N 2^{2k}\min\left(
        \sum^N_{l=0}2^{k(l-2)},
        \sum^N_{l=0}2^{k(l-\frac{2}{q})}
        \right).
    $$
So, choosing $N$ large enough, we obtain
    $$
        \sum_{0\leq k\leq k_0}2^{-\sigma k}\|\Delta^h_k
        \phi_\varepsilon\|_{L^q}
                    \leq C_\phi\sum_{0\leq k\leq k_0}
                    2^{k(N-\sigma)}\varepsilon^N
        \leq C_\phi 2^{k_0(N-\sigma)}\varepsilon^N
    $$
and
        $$
  \|S^h_0\phi_\varepsilon\|_{L^2(\mathbb{R}^2)}\leq
\sum_{k\leq-1}\|\Delta^h_k\phi_\varepsilon\|_{L^2} \leq
\sum_{k\leq-1}C_\phi\varepsilon^N2^{2k(1-\frac{1}{2})}\leq C_\phi
\varepsilon^N.
    $$
Choosing the best $k_0$, we have
               $$
                              \|\phi_\varepsilon\|_{\tilde{B}^{-\sigma}_{q,1}}\leq C_\phi
               \varepsilon^{\sigma}.
               $$

Similarly, since $\alpha<2(1-\frac{1}{q})$, we obtain
    $$
\sum_{k\leq-1}2^{-\alpha k}\|\Delta^h_k\phi_\varepsilon\|_{L^q}
\leq
\sum_{k\leq-1}C_\phi\varepsilon^N2^{k(2(1-\frac{1}{q})-\alpha)}\leq
C_\phi \varepsilon^N
    $$
and
    $$
                              \|\phi_\varepsilon\|_{\dot{B}^{-\alpha}_{q,1}}\leq C_\phi
               \varepsilon^{\alpha}.
               $$

 From (1.1) in \cite{Chemin07-0} (or Definition 1.1 in \cite{Chemin07}), we
have
    $$
    \|f\|_{\dot{B}^{-\sigma}_{q,p}(\mathbb{R}^2)}\simeq \left\|
    t^{\frac{\sigma}{2}}\|e^{t\Delta_h}f\|_{L^q}
    \right\|_{L^p(\mathbb{R}^+,\frac{dt}{t})}
    $$
and
    $$
        \|\phi_\varepsilon\|_{\dot{B}^{-\sigma}_{q,\infty}(\mathbb{R}^2)}\simeq
        \sup_{t>0}t^{\frac{\sigma}{2}}\|e^{t\Delta_h}\phi_\varepsilon\|_{L^q}\geq
        C\varepsilon^{\sigma}\|e^{\varepsilon^2\Delta_h}\phi_\varepsilon\|_{L^q},
        \ \sigma>0.
    $$
For any function $g$ satisfying
supp$\hat{g}\in\varepsilon^{-1}\mathcal{C}_h$, we have
    $$
    \|\mathcal{F}^{-1}(e^{\varepsilon^2|\xi_h|^2}\hat{g})\|_{L^q}\leq
    C\|g\|_{L^q}.
    $$
Since the support of $\mathcal{F}\phi_\varepsilon$ is included in
$\varepsilon^{-1}\mathcal{C}_h$ for some ring $\mathcal{C}_h$,
applied with $g=e^{\varepsilon^2\Delta}\phi_\varepsilon$, this
inequality gives
    $$
    \|\phi_\varepsilon\|_{L^q}\leq
    C\|e^{\varepsilon^2\Delta}\phi_\varepsilon\|_{L^q}
    \ \textrm{ and }
    \    \|\phi_\varepsilon\|_{\dot{B}^{-\sigma}_{q,\infty}}\geq
    C^{-1}\varepsilon^{\sigma}\|\phi\|_{L^q}.
    $$
From (\ref{ANS-E1.4}), we have
    \begin{eqnarray*}
      &&\int^\infty_0\|\phi_{\varepsilon,F'}\psi_{\varepsilon,F'}\|_{L^2(\mathbb{R}^2)}dt\\
      &\leq&\left(
      \int^\infty_0\|\phi_{\varepsilon,F'}\|_{L^4}^2dt\right)^\frac{1}{2}
      \left(\int^\infty_0\|\psi_{\varepsilon,F'}\|_{L^4}^2dt\right)^\frac{1}{2}\\
        &=&\nu_h^{-1}\left(\int^\infty_0\|e^{t\Delta_h}\sum_{k\geq 0}\Delta^h_k
        \phi_{\varepsilon}\|_{L^4}^2dt\right)^\frac{1}{2}
            \left(\int^\infty_0\|e^{t\Delta_h}\sum_{k\geq 0}\Delta^h_k
        \psi_{\varepsilon}\|_{L^4}^2dt\right)^\frac{1}{2}\\
        &\simeq& \nu_h^{-1}\|\sum_{k\geq 0}\Delta^h_k
        \phi_{\varepsilon}\|_{\dot{B}^{-1}_{4,2}}\|\sum_{k\geq 0}\Delta^h_k
        \psi_{\varepsilon}\|_{\dot{B}^{-1}_{4,2}}
        \\
        &\leq&\nu_h^{-1}\|\sum_{k\geq 0}\Delta^h_k
        \phi_{\varepsilon}\|_{\dot{B}^{-1}_{4,1}}\|\sum_{k\geq 0}\Delta^h_k
        \psi_{\varepsilon}\|_{\dot{B}^{-1}_{4,1}}
        \leq C_{\phi,\psi,\nu_h}\varepsilon^2.
    \end{eqnarray*}
This concludes the proof of Proposition \ref{ANS-P1}.  {\hfill
$\square$\medskip}

\section{An imbedding result}\label{ANS-Sec6}
\textbf{Proof of Proposition \ref{ANS-P3}.} It is easy to obtain
that $  \dot{B}^{-1}_{\infty,2}\subset BMO^{-1}\subset
  \dot{B}^{-1}_{\infty,\infty}=C^{-1},$
  (See \cite{Chemin07-0}). Thus, we only need to prove that
    $B^{-1+\frac{2}{p},\frac{1}{2}}_p\subset \dot{B}^{-1}_{\infty,2}.$
From Lemma \ref{ANS-L2.1} and Young's inequality, we have
    \begin{eqnarray*}
            &&\sum_{q\in\mathbb{Z}}2^{-2q}\|\Delta_q u\|^2_{L^\infty}\\
      &=&
           \sum_{q\in\mathbb{Z}}2^{-2q}\|\sum_{l\in\mathbb{Z}}\Delta_q(\sum_{k\geq l-1}\Delta^h_k\Delta^v_l u
           +S^h_{l-1}\Delta^v_lu)\|_{L^\infty}^2\\
                &\leq&           \sum_{q\in\mathbb{Z}}\sum_{l\leq q+N_0}\sum_{|k-q|\leq N_0}2^{-2q}
                \|
                \Delta_q\Delta^h_k\Delta^v_l u\|_{L^\infty}^2
                +\sum_{q\in\mathbb{Z}}\sum_{l\geq q-N_0}2^{-2q}\|\Delta_qS^h_{l-1}\Delta^v_lu\|_{L^\infty}^2\\
           &\lesssim&    \sum_{q\in\mathbb{Z}}\sum_{l\leq q+N_0}\sum_{|k-q|\leq
           N_0}2^{-2q}2^l2^\frac{4k}{p}
                \|
                \Delta^h_k\Delta^v_l u\|_{L^p_h(L^2_v)}^2
                +\sum_{q\in\mathbb{Z}}\sum_{l\geq q-N_0}2^{-2q}2^{3q}
                \|S^h_{l-1}\Delta^v_lu\|_{L^2}^2\\
            &\lesssim& \left(   \sum_{q\in\mathbb{Z}}\sum_{l\leq q+N_0}\sum_{|k-q|\leq
           N_0}2^{-2q+2k}d_{k,l}^2
                +\sum_{q\in\mathbb{Z}}\sum_{l\geq
                q-N_0}2^{q-l}d_l^2\right)
                \|u\|_{B^{-1+\frac{2}{p},\frac{1}{2}}_p}^2\\
            &\lesssim&
                \|u\|_{B^{-1+\frac{2}{p},\frac{1}{2}}_p}^2,
                \ p\geq2.\\
    \end{eqnarray*}
Then, we finish the proof of Proposition \ref{ANS-P3}.
 {\hfill $\square$\medskip}

\end{document}